\documentclass[10pt,twoside]{amsart}

\usepackage{psfrag}

\usepackage{array,amscd,amsmath,amsthm}
\usepackage{amssymb}
\usepackage{graphicx}
\usepackage[all]{xy}
\allowdisplaybreaks

\makeatletter

\newenvironment{figurehere}
  {\def\@captype{figure}}
  {}
\makeatother

\newcommand\myCaption[1]{\small\refstepcounter{figure}%
   \centering\figurename\ \thefigure :\ #1}

\newenvironment{theorem2}[2][Theorem]{\begin{trivlist}
\item[\hskip \labelsep {\bfseries #1}\hskip
\labelsep {\bfseries #2}{\bfseries{.}}]\it}{\end{trivlist}}

\newenvironment{corollary2}[2][Corollary]{\begin{trivlist}
\item[\hskip \labelsep {\bfseries #1}\hskip
\labelsep {\bfseries #2}{\bfseries{.}}]\it}
{\end{trivlist}}

\newtheorem*{theorem*}{Theorem}

\newtheorem{theorem}{Theorem}[section]
\newtheorem{proposition}[theorem]{Proposition}
\newtheorem{lemma}[theorem]{Lemma}
\newtheorem{lemmadef}[theorem]{Lemma-Definition}
\newtheorem{corollary}[theorem]{Corollary}

\newtheorem{remark}[theorem]{Remark}
\newtheorem{definition}[theorem]{Definition}
\newtheorem{example}[theorem]{Example}
\newtheorem*{notation}{Notation}

\newcommand{\C}{\mathbb{C}}
\newcommand{\Z}{\mathbb{Z}}

\newcommand{\R}{\mathbb{R}}
\newcommand{\N}{\mathbb{N}}

\def\Cc{\mathcal{C}}
\def\Aa{\mathcal{A}}

\def\Hh{\mathcal{H}}
\def\extended#1{\check{#1}}
\def\Teich{\mathcal{T}}
\def\eTeich{\extended{\Teich}}
\def\Tbar{\overline{\Teich}}
\def\That{\widehat{\Teich}}
\def\Ttil{\widetilde{\Teich}}
\def\M{\mathcal{M}}

\def\Mbar{\overline{\M}}
\def\Mhat{\widehat{\M}}
\def\Mtil{\widetilde{\M}}
\def\pa{\partial}
\def\Pr{\mathbb{P}}
\def\Aut{\mathrm{Aut}}
\def\Iso{\mathrm{Iso}}

\def\rar{\rightarrow}
\def\arr#1#2{\stackrel{#1}{#2}}
\def\hra{\hookrightarrow}
\def\a{\alpha}
\def\Af{\mathfrak{A}}
\def\Cf{\mathfrak{C}}
\def\b{\beta}

\def\Ao{\Af^\circ}
\def\A8{\Af_{\infty}}

\def\l{\lambda}

\def\i{\iota}
\def\g{\gamma}
\def\G{\Gamma}

\def\d{\delta}
\def\lra{\longrightarrow}

\def\Ll{\mathcal{L}}

\def\s{\sigma}
\def\Si{\Sigma}
\def\e{\varepsilon}

\def\ul#1{\underline{#1}}
\def\ol#1{\overline{#1}}
\def\wh#1{\widehat{#1}}
\def\dis{\displaystyle}
\def\ti#1{\tilde{#1}}

\def\Sp{\mathrm{Sp}}
\def\Pr{\mathbb{P}}
\def\Qq{\mathcal{Q}}
\def\ora#1{\overrightarrow{#1}}
\def\ola#1{\overleftarrow{#1}}

\def\up{\ul{p}}

\def\tw{\mathrm{tw}}

\def\ML{\mathcal{ML}}
\def\th{\vartheta}
\def\bm#1{\text{\boldmath${#1}$}}
\def\dev{\mathrm{dev}}
\def\Pp{\mathcal{P}}
\def\pp{\mathfrak{p}}
\def\Sch{\text{\boldmath$S$}}
\def\gr8{\mathrm{gr}_{\infty}}
\def\To{\mathbb{T}}
\def\Ho{\bm{H}}
\def\Bel{\mathcal{B}}
\def\Pj{\mathcal{P}}
\def\Jac{\mathcal{J}}

\begin{document}

\title[Riemann surfaces with boundary]{Riemann surfaces with boundary
and natural triangulations of the {T}eichm\"uller space}

\author{Gabriele Mondello}

\address{Imperial College of London\\
Department of Mathematics, Huxley Building\\
South Kensington Campus\\
SW7 2AZ - London, UK}

\email{g.mondello@imperial.ac.uk}


\begin{abstract}
We compare some natural triangulations of the Teichm\"uller space
of hyperbolic surfaces with geodesic boundary and of some bordifications.
We adapt Scannell-Wolf's proof to show that grafting semi-infinite
cylinders at the ends of hyperbolic surfaces with fixed boundary lengths
is a homeomorphism. This way, we construct a family of
equivariant triangulations
of the Teichm\"uller space of punctured surfaces that interpolates
between Penner-Bowditch-Epstein's (using the spine construction)
and Harer-Mumford-Thurston's (using Strebel's differentials).
Finally, we show (adapting arguments of Dumas)
that on a fixed punctured surface,
when the triangulation approaches HMT's,
the associated Strebel differential is well-approximated
by the Schwarzian of the
associated projective structure and by the Hopf differential of
the collapsing map.
\end{abstract}


\maketitle

%
\begin{section}{Introduction}
The aim of this paper is to compare two different ways of
triangulating the Teichm\"uller space $\Teich(R,x)$ of
conformal structures on a compact oriented surface
$R$ with distinct ordered marked points $x=(x_1,\dots,x_n)$.
Starting with $[f:R\rar R']\in\Teich(R,x)$ and a collection of
weights $\up=(p_1,\dots,p_n)\in\Delta^{n-1}$,
both constructions produce a ribbon graph $G$
embedded inside the punctured surface
$\dot{R}=R\setminus x$ as a deformation retract,
together with a positive weight for each edge.
The space of such weighted graphs can be identified
to the topological realization of the arc complex $\Af(R,x)$
(via Poincar\'e-Lefschetz duality on $(R,x)$, see for instance
\cite{mondello:survey}),
which is the simplicial complex of (isotopy classes of)
systems of (homotopically nontrivial, pairwise nonhomotopic)
arcs that join couples of marked points and that admit
representatives with disjoint interior (\cite{harer:virtual},
\cite{bowditch-epstein:natural}, \cite{looijenga:cellular}).

Thus, both constructions provide
a $\G(R,x)$-equivariant homeomorphism
$\Teich(R,x)\times\Delta^{n-1}\rar |\Ao(R,x)|$,
where $\G(R,x)=\pi_0\mathrm{Diff}_+(R,x)$ is the mapping class group
of $(R,x)$ and $\Ao(R,x)\subset\Af(R,x)$ consists of proper
systems of arcs $\bm{A}=\{\a_0,\dots,\a_k\}$, namely such that
$\dot{R}\setminus(\a_0\cup\dots\cup\a_k)$ is a disjoint union of discs
and pointed discs. In fact, properness of $\bm{A}$ is exactly
equivalent to its dual ribbon graph being a deformation retract
of $\dot{R}$.

The HMT construction (due to Harer, Mumford and Thurston)
appears in \cite{harer:virtual}. It uses Strebel's result
\cite{strebel:67} on existence and uniqueness of meromorphic
quadratic differential $\varphi$ on a Riemann surface
$R$ with prescribed residues $\up$ at $x$ to decompose $\dot{R}$ into a
disjoint union of semi-infinite $|\varphi|$-flat cylinders (one for each
puncture $x_i$ with $p_i>0$), that are identified along a
critical graph $G$ which inherits this way a metric.
The length of each edge of $G$ will be its weight.

The PBE construction (due Penner \cite{penner:decorated}
and Bowditch-Epstein \cite{bowditch-epstein:natural}) uses
the unique hyperbolic metric on the punctured Riemann
surface $\dot{R}$. Given a {\it (projectively) decorated surface},
that is a hyperbolic surface $\dot{R}$ with cusps plus
a weight $\up\in\Delta^{n-1}$, there are
disjoint embedded horoballs of circumference $p_1,\dots,p_n$ at
the $n$ cusps of $\dot{R}$. Removing the horoballs, we obtain
a truncated surface $R^{tr}$ with boundary, on which the function
``distance from the boundary'' is well-defined. The critical locus
of this function is a spine $G$ embedded in $R^{tr}\subset\dot{R}$
as a deformation retract and with geodesic edges, whose horocyclic
lengths provide the associated weights.

Both constructions share similar properties
of homogeneity and real-analiticity (see \cite{hubbard-masur:foliations}
and \cite{penner:decorated}) and they also enjoy some
good compatibility with the Weil-Petersson symplectic structure
on $\Teich(R,x)$, as explained later.

In this paper, we will interpolate these two constructions using the
Teichm\"uller space $\Teich(S)$ of {\it hyperbolic
surfaces with geodesic boundary} (see also
\cite{luo:decomposition}), where $S$ is a surface with boundary
endowed with a homotopy equivalence $S\hra\dot{R}$.
The spine construction works perfectly on such surfaces and it can
be easily seen to reduce to the PBE case
as the boundary lengths $(p_1,\dots,p_n)\in\R_+^n\hra\Delta^{n-1}\times
(0,\infty)$ become infinitesimal
(see also \cite{mondello:poisson}). Also, the Weil-Petersson
Poisson structure can be explicitly determined, thus obtaining
a generalization of Penner's formula \cite{penner:wp}.

Thus, the limit $\bm{p}:=p_1+\dots+p_n\rar 0$
is completely understood and it behaves
as the Weil-Petersson completion (or Bers's augmentation
\cite{bers:degenerating}) of the Teichm\"uller space.

Instead, the limit $\bm{p}\rar\infty$ behaves more like Thurston's
compactification \cite{FLP} of the Teichm\"uller space; in fact,
the arc complex $|\Af(S)|$ naturally embeds inside the space of
projective measured laminations. From a symplectic point of view,
the Weil-Petersson structure admits a precise limit as $\bm{p}\rar\infty$,
after a suitable normalization, which agrees with Kontsevich's
piecewise-linear symplectic form on $|\Af(S)|$
defined in \cite{kontsevich:intersection} (see \cite{mondello:poisson}).

To give a more geometric framework to these limiting considerations,
we produce a few different bordifications of the Teichm\"uller
space $\Teich(S)$ of a surface $S$ with boundary, whose quotients by
the mapping class group $\G(S)$
give different compactifications of the moduli space.
A convenient bordification from the point
of view of the Weil-Petersson Poisson structure is the {\it extended
Teichm\"uller space $\Ttil(S)$}; whereas the most suitable one for
triangulations and spine constructions is the
{\it bordification of arcs} $\Tbar^a(S)$, whose definition looks
a bit like Thurston's but with some relevant differences
(for instance, we use $t$-lengths related to hyperbolic collars instead
of hyperbolic lengths).
It is reasonable to believe that careful iterated
blow-ups of $\Tbar^a(S)$ along
its singular locus would produce finer bordifications of $\Teich(S)$
in the spirit of \cite{looijenga:cellular}
(see also \cite{mcshane-penner:screens}).

In order to explicitly
link the HMT and PBE constructions, we construct an isotopic
family of triangulations of $\Teich(R,x)\times\Delta^{n-1}$, parametrized
by $t\in[0,\infty]$, that coincides with PBE for $t=0$ and with
HMT for $t=\infty$. In particular, we prove that, for every complex
structure on $\dot{R}$ and every $(\up,t)\in\Delta^{n-1}\times[0,\infty]$,
there exists a unique projective structure
$\mathcal{P}(\dot{R},t\up)$ on $\dot{R}$,
whose associated Thurston metric has flat cylindrical ends
(with circumferences $t\up$) and a hyperbolic core.
Rescaling the lengths by a factor $1/t$, we recognize that at $t=\infty$
the hyperbolic core shrinks to a graph $G$ and the metric is of the type
$|\varphi|$, where $\varphi$ is a Strebel differential.
This result can be restated in term of infinite grafting at the ends
of a hyperbolic surface with geodesic boundary and the proof adapts
arguments of Scannell-Wolf \cite{scannell-wolf:grafting}.

Finally, we show that, for large $t$, two results of Dumas
\cite{dumas:grafting} \cite{dumas:schwarzian} for compact surfaces
still hold. The first one says that, for $t$ large, the Strebel
differential $\varphi$ is well-approximated
in $L^1_{loc}(\dot{R})$ by the Hopf differential of the collapsing map
associated to $\mathcal{P}(\dot{R},t\up)$,
that is the quadratic differential which writes $dz^2$ on the
flat cylinders $S^1\times[0,\infty)$ and is zero on the hyperbolic part.
The second result says that $\varphi$ is also well-approximated
by the Schwarzian derivative of the projective structure
$\mathcal{P}(\dot{R},t\up)$.
%
%
\begin{subsection}{Content of the paper}
In Section~\ref{sec:preliminaries}, we recall basic concepts like
Teichm\"uller space $\Teich(R)$, measured laminations $\ML(R)$
and Thurston's compactification $\Tbar^{Th}(R)=\Teich(R)\cup\ML(R)$,
when $R$ is an oriented compact surface with $\chi(R)<0$.

We also extend these concepts to the case of an oriented
surface $S$ with boundary and $\chi(S)<0$,
using the doubling construction $S\rightsquigarrow dS$. We also remark
that the arc complex $|\Af(S)|$ embeds in
$\ML(S)=\ML(dS)^\s$ (where $\s$ is the natural anti-holomorphic
involution of $dS$) and, even though its image is neither open nor closed,
the subspace topology coincides with the metric topology.

Next, we introduce the Weil-Petersson pairing on a closed surface and
on a surface with boundary, we describe the augmentation $\Tbar^{WP}$
of the Teichm\"uller space and we restate Wolpert's
formula \cite{wolpert:fenchel-nielsen},
which expresses the WP symplectic structure
in Fenchel-Nielsen coordinates. Finally, we recall the definition of the
mapping class group $\G(S)=\pi_0\mathrm{Diff}_+(S)$, the moduli space
$\M(S)=\Teich(S)/\G(S)$ and the Deligne-Mumford compactification.\\

We begin Section~\ref{sec:triangulations} by defining the geometrical
quantities that are associated to an arc in a hyperbolic surface $S$
with boundary: the hyperbolic
length $a_i=\ell_{\a_i}$ of (the geodesic representative of)
$\a_i\in\Aa(S)=\{\text{isotopy classes of arcs in $S$}\}$, its associated
$s$-length $s_{\a_i}=\cosh(a_i/2)$ and $t$-length $t_{\a_i}=T(\ell_{\a_i})$,
where $T$ is defined by $\sinh(T(x)/2)\sinh(x/2)=1$.
The $t$-lengths give a continuous embedding
\[
j:\xymatrix@R=0in{
\Teich(S)\ar@{^(->}[r] &  \Pr L^\infty(\Aa(S))\times[0,\infty]\\
[f:S\rar\Si] \ar@{|->}[r] & ([t_\bullet(f)],\|t_\bullet(f)\|_\infty)
}
\]
and we call {\it bordification of arcs}
the closure of its image $\Tbar^a(S)$.

Then we define the {\it spine} $\mathrm{Sp}(\Si)$
(of a hyperbolic surface $\Si$)
as the critical locus of the function ``distance from the boundary''
and we produce its dual {\it spinal arc system} $\bm{A}_{sp}\in\Af(\Si)$
and a system of weights (the {\it widths})
$w_{sp}$ so that $w_{sp}(\a)$ is the length
of either of the two projections of the edge $\a^*$ of the spine (dual to $\a$)
to the boundary.
We also define the width of an arc $\a$ (and of an oriented
arc $\ora{\a}$) associated to a maximal system of arcs $\bm{A}$
and we show that the two concepts agree \cite{ushijima:decomposition}
(see also \cite{mondello:poisson}).

We recall the PBE and Luo's result on the cellularization of $\Teich(S)$
using the spine construction.

\begin{theorem*}
Let $S$ be a compact oriented surface with $n\geq 1$ boundary components
and $\chi(S)<0$ and let $(R,x)$ be a pointed surface such that $S\hra
\dot{R}$ is a homotopy equivalence.
\begin{itemize}
\item[(a)]
If $\Teich(R,x)\times\Delta^{n-1}$
is the Teichm\"uller space of (projectively) decorated surfaces,
then the spine construction
\[
\bm{W}_{PBE}:
\xymatrix@R=0in{
\Teich(R,x)\times\Delta^{n-1}\ar[r]      & |\Ao(R,x)| \\
([f:R\rar R'],\up) \ar@{|->}[r] &     f^* \ti{w}_{sp,R',\up}
}
\]
induces a $\G(R,x)$-equivariant homeomorphism (\cite{penner:decorated},
\cite{bowditch-epstein:natural}).
\item[(b)]
The spine construction applied to hyperbolic surfaces
with geodesic boundary
\[
\bm{W}:
\xymatrix@R=0in{
\Teich(S)\ar[r]      & |\Ao(S)|\times \R_+ \\
[f:S\rar\Si] \ar@{|->}[r] &     f^* w_{sp,\Si}
}
\]
gives a $\G(S)$-equivariant homeomorphism (\cite{luo:decomposition}).
\end{itemize}
\end{theorem*}

To deal with stable surfaces, we first define
$\That(S)$ as the real blow-up
of $\Tbar^{WP}(S):=\bigcup_{\up\in\Delta^{n-1}\times[0,\infty)}
\Tbar^{WP}(S)(\up)$
along the locus $\Tbar^{WP}(S)(0)$
of surfaces with $n$ boundary cusps
and we identify the exceptional locus $\That(S)(0)$
with the space of projectively decorated surfaces (that is, of surfaces
with $n$ boundary cusps and weights $(p_1,\dots,p_n)\in\Delta^{n-1}$).
Then, we call {\it visible} the subsurface $\Si_+\subset\Si$
consisting of the components of $\Si$ which have
positive boundary length (or some positively weighted cusp)
and we declare that $[f_1:S\rar\Si_1]$
and $[f_2:S\rar\Si_2]$ in $\That(S)$
are {\it visibly equivalent}
if there exists a third $[f:S\rar\Si]$ and maps
$h_i:\Si\rar\Si_i$ for $i=1,2$ that are isomorphisms
on the visible components and such that
$h_i\circ f\simeq f_i$ for $i=1,2$.

\begin{theorem2}{\ref{prop:w}}
Let $S$ be a compact oriented surface with $n\geq 1$ boundary components
and $\chi(S)<0$.
The spine construction gives a $\G(S)$-equivariant
homeomorphism
\[
\wh{\bm{W}}:
\xymatrix@R=0in{
\That^{vis}(S)\ar[r]      & |\Af(S)|\times[0,\infty) \\
[f:S\rar\Si] \ar@{|->}[r] &     f^* w_{sp,\Si}
}
\]
where $\That^{vis}(S)$ is obtained from $\That(S)$ by identifying
visibly equivalent surfaces. Moreover, $\wh{\bm{W}}$ extends
Penner's and Luo's constructions.
\end{theorem2}

We then consider the bordification of arcs.

\begin{theorem2}{\ref{thm:phi}}
The map $\Phi:|\Af(S)|\times[0,\infty]\lra\Tbar^a(S)$
defined as
\[
\Phi(w,\bm{p})=
\begin{cases}
([\l^{-1}_\bullet(\wh{\bm{W}}^{-1}(w,0))],0) &
\text{if $\bm{p}=0$} \\
j(\wh{\bm{W}}^{-1}(w,\bm{p})) & \text{if $0<\bm{p}<\infty$}\\
([w],\infty) & \text{if $\bm{p}=\infty$}
\end{cases}
\]
is a $\G(S)$-equivariant homeomorphism, where $\l_\a$ is Penner's
$\l$-length of $\a$.
\end{theorem2}
The situation is illustrated in the following
$\G(S)$-equivariant commutative diagram
\[
\xymatrix{
\That(S) \ar@{^(->}[d] \ar@{->>}[r] &  \That^{vis}(S)
\ar@{^(->}[d] \ar[rr]^{\wh{\bm{W}}\qquad}_{\sim\qquad} &&
|\Af(S)|\times[0,\infty) \ar@{^(->}[d] \\
\Ttil(S) \ar@{->>}[r]  &    \Tbar^a(S)  
&& |\Af(S)|\times[0,\infty] \ar[ll]_{\Phi\qquad}^{\sim\qquad}
}
\]
in which
$\Tbar^a(S)$ is exhibited as a quotient of the {\it extended
Teichm\"uller space} $\Ttil(S):=\That(S)\cup|\Af(S)|_\infty$
(endowed with a suitable topology, where $|\Af(S)|_\infty$ is
just a copy of $|\Af(S)|$) by visible equivalence.

Section~\ref{sec:wp} describes how to extend the previous
triangulations to the case of a surface with boundary $S$
and a marked point $v_i$ on each boundary component $C_i$
(with $i=1,\dots,n$),
so that we obtain a commutative diagram
\[
\xymatrix{
\That^{vis}(S,v) \ar[rr]^{\wh{\bm{W}}_v\quad}\ar[d] &&
|\Af(S,v)|\times[0,\infty) \ar[d] \\
\That^{vis}(S) \ar[rr]^{\wh{\bm{W}}\quad} && |\Af(S)|\times[0,\infty)
}
\]
in which the horizontal arrows are $\G(S,v)$-equivariant
homeomorphisms and
the vertical arrows are $\R^n$-fibrations on the smooth locus
(with some possible degenerations on the stable surfaces).
After passing to the associated moduli spaces,
the vertical arrows become $(S^1)^n$-bundles,
which are actually products of the circle bundles
$L_1,\dots,L_n$ associated to the respective boundary components
$C_1,\dots,C_n$.
This $(S^1)^n$-action is Hamiltonian for the Weil-Petersson
structure with moment map $\mu=(p_1^2/2,\dots,p_n^2/2)$ and this
shows that $[\omega_{\up}]=[\omega_0]+\sum_i p_i^2/2 [c_1(L_i)]$
in cohomology (\cite{mirzakhani:volumes}), where $\omega_{\up}$
is the restriction of the Weil-Petersson form to the symplectic leaf
$\Mhat(S)(\up)$, that is the moduli space
of surfaces with boundary lengths $\up$.
Pointwise, the Poisson structure $\eta$ on $\Mhat(S)$
can be described as follows.

\begin{theorem2}{\ref{thm:poisson} {\normalfont{(\cite{mondello:poisson})}}}
Let $\bm{A}$ be a maximal system of arcs on $S$.
Then
\[
\eta=\frac{1}{4}\sum_{k=1}^n
\sum_{\substack{y_i\in\a_i\cap C_k \\ y_j\in\a_j\cap C_k}}
\frac{\sinh(p_k/2-d_{C_k}(y_i,y_j))}{\sinh(p_k/2)}
\frac{\pa}{\pa a_i}\wedge\frac{\pa}{\pa a_j}
\]
on $\Teich(S)$,
where $d_{C_k}(y_i,y_j)$ is the length
of the geodesic running from $y_i$ to $y_j$ along $C_k$ in the positive
direction. Moreover, if we normalize $\ti{w}_i=
(\bm{p}/2)^{-1}w_i$ and $\ti{\eta}=(1+\bm{p}/2)^2\eta$,
then $\ti{\eta}$ extends to
$\Ttil(S)$ and
\[
\tilde{\eta}_\infty=\frac{1}{2}\sum_{r}\left(
\frac{\pa}{\pa \tilde{w}_{r_1}}\wedge\frac{\pa}{\pa\tilde{w}_{r_2}}+
\frac{\pa}{\pa \tilde{w}_{r_2}}\wedge\frac{\pa}{\pa\tilde{w}_{r_3}}+
\frac{\pa}{\pa \tilde{w}_{r_3}}\wedge\frac{\pa}{\pa\tilde{w}_{r_1}}
\right)
\]
where $r$ ranges over all (trivalent) vertices of the ribbon graph
representing a point in $|\Ao(S)|$ and $(r_1,r_2,r_3)$ is the (cyclically)
ordered triple of edges incident at $r$.
\end{theorem2}

The result can be seen to reduce to Penner's formula
\cite{penner:wp} as $\bm{p}=0$.\\

Finally, Section~\ref{sec:grafting} relates hyperbolic surfaces
with boundary homeomorphic to $S$ to punctured surfaces
homeomorphic to $R\setminus x=\dot{R}\simeq S$.
We describe first Strebel's result and the HMT construction
and its extension to $\That^{vis}(R,x)$
(see \cite{mondello:survey}),
which provides a $\G(R,x)$-equivariant homeomorphism
\[
\bm{W}_{HMT}:\That^{vis}(R,x)\lra |\Af(R,x)|
\]
Then, we define
$\gr8(\Si)\in\Tbar^{WP}(R,x)\times\Delta^{n-1}\times[0,\infty)$ to be the Riemann
surface obtained from the hyperbolic surface $\Si\in\That(S)$ with geodesic
boundary by {\it grafting semi-infinite flat cylinders at its ends}.
Moreover, for every $w \in |\Af(S)|_\infty\cong|\Af(R,x)|$,
we let $\gr8(w):=\bm{W}_{HMT}^{-1}(w)$.

The key result is the following.

\begin{theorem2}{\ref{thm:grafting}}
The map
$(\gr8,\Ll):\Ttil(S)\rar \Tbar^{WP}(R,x)\times\Delta^{n-1}
\times[0,\infty]/\!\!\sim$
is a homeomorphism, where
$\Ll$ is the boundary length map and $\sim$ identifies
$([f_1:R\rar R_1],\up,\infty)$ and $([f_2:R\rar R_2],\up,\infty)$
if and only if $([f_1],\up)$ and $([f_2],\up)$ are visibly equivalent.
\end{theorem2}

The continuity at infinity requires some explicit computations,
whereas the proof of the injectivity
simply adapts arguments of Scannell-Wolf \cite{scannell-wolf:grafting}
to our situation.

We can summarize our results in the following commutative
diagram
\[
\xymatrix{
\That^{vis}(R,x)\times[0,\infty]
\ar[rrd]_{\Psi} &&
\Tbar^a(S) \ar[ll]_{\qquad\qquad (\gr8,\Ll)} \\
&& |\Af(S)|\times[0,\infty] \ar[u]_{\Phi}
}
\]
in which $\Psi=\Phi^{-1}\circ(\gr8,\Ll)^{-1}$.

\begin{corollary2}{\ref{cor:grafting}}
The maps $\Psi_t:\That^{vis}(R,x)\rar|\Af(R,x)|$
obtained by restricting $\Psi$ to
$\That^{vis}(R,x)\times\{t\}$
form a continuous family of $\G(S)$-equivariant triangulations,
which specializes to PBE for $t=0$ and to HMT to $t=\infty$.
\end{corollary2}

The last result concerns the degeneration of
the projective structure $\mathrm{Gr}_\infty(\Si)$
on the Riemann surface $\mathrm{gr}_\infty(\Si)$.
It adapts arguments of Dumas \cite{dumas:grafting}
\cite{dumas:schwarzian} to our case.

\begin{theorem2}{\ref{thm:mydumas}}
Let $\{f_m:S\rar \Si_m\}\subset\Teich(S)$
be a sequence such that
$(\mathrm{gr}_\infty,\Ll)(f_m)=([f:R\rar R'],\up_m)\in\Teich(R,x)
\times\R_+^n$.
The following are equivalent:
\begin{enumerate}
\item
$\up_m\rar(\up,\infty)$ in $\Delta^{n-1}\times(0,\infty]$
\item
$[f_m]\rar [w]$ in $\Tbar^a(S)$,
where $[w]$ is the projective multi-arc associated to the
vertical foliation of the Jenkins-Strebel differential $\varphi_{JS}$
on $R'$ with weights $\up$ at $x'=f(x)$.
\end{enumerate}
When this happens, we also have
\begin{itemize}
\item[(a)]
$4\bm{p_m}\!\!\!^{-2}\Ho(\kappa_m)\rar
\varphi_{JS}$ in $L^1_{loc}(R',K(x')^{\otimes 2})$, where
$\Ho(\kappa_m)$ is the Hopf differential of the
collapsing map $\kappa_m:R'\rar \Si_m$
\item[(b)]
$2\bm{p_m}\!\!\!^{-2}\Sch(\mathrm{Gr}_\infty([f_m]))\rar -\varphi_{JS}$
in $H^0(R',K(x')^{\otimes 2})$, where
$\Sch$ is the Schwarzian derivative
with respect to the Poincar\'e projective structure.
\end{itemize}
\end{theorem2}
\end{subsection}
%
%
\begin{subsection}{Acknowledgments}
I would like to thank Enrico Arbarello, Curtis McMullen,
Tomasz Mrowka, Mike Wolf and Scott Wolpert for very fruitful
and helpful discussions.
\end{subsection}
\end{section}
%
%
\begin{section}{Preliminaries}\label{sec:preliminaries}
%
%
\begin{subsection}{Double of a surface with boundary}
By a {\it surface} with boundary and/or marked points we will always mean
a compact oriented surface $S$ possibly
with boundary and/or distinct ordered marked points $x=(x_1,\dots,x_n)$,
with $x_i\in S^\circ$.
By a {\it nodal surface} with boundary and marked points
we mean a compact, Hausdorff topological space $S$ with countable basis
in which every $q\in S$ has an open neighbourhood $U_q$ such that
$(U_q,q)$ is homeomorphic to: either
$(\C,0)$ and $q$ is called {\it smooth} point; or
$(\{z\in\C\,|\,\mathrm{Im}(z)\geq 0\},0)$ and $q$ is called {\it boundary}
point; or $(\{(z,w)\in\C^2\,|\,zw=0\},0)$ and $q$ is called {\it node}.

We will say that a (nodal) surface
$S$ is {\it closed} if it has no boundary and
no marked points.

A {\it hyperbolic metric} on $S$ is a complete metric $g$
of finite volume on the smooth locus $\dot{S}_{sm}$ of the {\it punctured
surface} $\dot{S}:=S\setminus x$ of constant
curvature $-1$, such that $\pa S$ is geodesic.
Clearly, such a $g$ acquires cusps at the marked points and at the nodes.

Given a (possibly nodal)
surface $S$ with boundary and/or marked points, we can construct
its double $dS$ in the following way. Let $S'$ be another copy of $S$,
with opposite orientation, and call $q'\in S'$ the point corresponding
to $q\in S$.
Define $dS$ to be $S\coprod S'/\!\!\sim$,
where $\sim$ is the equivalence relation generated by $q\sim q'$ for
every $q\in\pa S$ and every $x_i$.
Clearly, $dS$ is closed and it is smooth whenever
$S$ has no nodes and no marked points.

$dS$ can be oriented
so that the natural embedding $\i:S\hra dS$ is
orientation-preserving.
Moreover, $dS$ comes naturally equipped with an
orientation-reversing involution
$\sigma$ that fixes the boundary and the cusps of $\i(S)$
and such that $dS/\sigma\cong S$. If $S$ is hyperbolic,
then $dS$ can be given a hyperbolic metric such that $\i$
and $\sigma$ are isometries.

Clearly, on $dS$ there is a correspondence between complex structures
and hyperbolic metrics and, in fact, $\s$-invariant hyperbolic
metrics correspond to complex structures such that $\s$ is
anti-holomorphic. Thus, the datum of a
hyperbolic metric with geodesic boundary
on $S$ is equivalent to that of a complex structure on $S$,
such that $\pa S$ is totally real.
\end{subsection}
%
%
\begin{subsection}{Teichm\"uller space}
Let $S$ be a hyperbolic surface with $n\geq 0$ boundary
components $C_1,\dots,C_n$ and no cusps.
\begin{definition}
An {\it $S$-marked} hyperbolic surface is an orientation
preserving map $f:S\lra\Si$ of (smooth) hyperbolic
surfaces that may shrink boundary components of $S$
to cusps of $\Si$ and that is a diffeomorphism
everywhere else.
\end{definition}

Two $S$-marked
surfaces $f_1:S\lra \Si_1$ and $f_2:S\lra\Si_2$
are {\it equivalent} if there exists an isometry $h:\Si_1\lra \Si_2$
such that $h\circ f_1$ is homotopic to $f_2$.

\begin{definition}
Call $\eTeich(S)$ the space of equivalence classes of
$S$-marked hyperbolic surfaces. The {\it Teichm\"uller
space} $\Teich(S)\subset\eTeich(S)$ is the locus of
surfaces $\Sigma$ with no cusps.
\end{definition}

The space $\mathfrak{Met}(S)$ of smooth metrics on $\dot{S}$ has
the structure of an open convex subset of a Fr\'echet space.
Consider the map $\mathfrak{Met}(S)\rar\eTeich(S)$ that associates
to $g\in\mathfrak{Met}(S)$ the unique hyperbolic
metric with geodesic boundary in the conformal class of $g$.
Endow $\eTeich(S)$ with the quotient topology.

Let $\bm{\g}=\{C_1,\dots,C_n,\g_1,\dots,\g_{3g-3+n}\}$
be a maximal system of disjoint
simple closed curves of $S$ such that no $\g_i$ is contractible
and no couple $\{\g_i,\g_j\}$ or $\{\g_i,C_j\}$ bounds a cylinder.
The system $\bm{\g}$ induces a {\it pair of pants decomposition} of $S$,
that is $S^\circ\setminus\bigcup_i \g_i=P_1\cup\dots\cup
P_{2g-2+n}$,
and each $P_i$ is a pair of pants
(i.e. a surface homeomorphic to $\C\setminus\{0,1\}$).

Given $[f:S\rar\Si]\in\Teich(S)$, we can define $\ell_i(f)$ to
be the length of the unique geodesic curve isotopic to $f(\g_i)$.
Let $\tau_i(f)$ be the associated twist parameter (whose definition
depends on some choices).

The {\it Fenchel-Nielsen coordinates} $(p_j,\ell_i,\tau_i)$
exhibit a homeomorphism $\eTeich(S)\arr{\sim}{\lra} \R_{\geq 0}^n
\times(\R_+\times\R)^{3g-3+n}$.
In particular, the {\it boundary length map}
$\Ll:\eTeich(S)\lra \R_{\geq 0}^n$ is
defined as $\Ll([f])=(p_1,\dots,p_n)$
and we write $\Teich(S)(\up):=\Ll^{-1}(\up)$
for $\up\in\R_{\geq 0}^n$.
Thus, $\Teich(S)=\eTeich(S)(\R_+^n)$.
\end{subsection}
%
%
%
\begin{subsection}{Measured laminations}\label{sec:laminations}
A {\it geodesic lamination} on a closed smooth
hyperbolic surface $(R,g)$
is a closed subset $\l\subset R$ which
is foliated in complete simple geodesics.
A {\it transverse measure} $\mu$ for $\l$ is a function
$\mu:\Lambda(\l)\lra\R_{\geq 0}$,
where $\Lambda(\l)$ is the collection of compact smooth arcs imbedded in $R$
with endpoints in $R\setminus\l$, such that
\begin{enumerate}
\item
$\mu(\a)=\mu(\b)$ if $\a$ is isotopic to $\b$ through elements of $\Lambda(\l)$
({\it homotopy invariance})
\item
$\dis \mu(\a)=\sum_{i\in I}\mu(\a_i)$ if $\a=\bigcup_{i\in I}\a_i$,
if $\a_i\in\Lambda(\l)$ for all $i$ in a countable set $I$
and distinct $\a_i,\a_j$ meet at most at their endpoints
({\it $\sigma$-additivity})
\item
for every $\a\in\Lambda(\l)$, $\mu(\a)>0$ if and only
if $\a\cap\l\neq \emptyset$
({\it the support of $\mu$ is $\l$}).
\end{enumerate}
In this case, the couple $(\l,\mu)$ is called a
{\it measured geodesic lamination}
on $(R,g)$ (often denoted just by $\mu$).

\begin{lemmadef}
If $g$ and $g'$ are hyperbolic metrics on $R$, then
there is a canonical identification between measured
$g$-geodesic laminations and measured $g'$-geodesic laminations.
Thus, we call the set $\ML(R)$ of such $(\l,\mu)$'s just
the space of {\it measured laminations} on $R$
(see \cite{FLP} and \cite{penner-harer:traintracks} for
more details).
\end{lemmadef}

Given a measured lamination $(\l,\mu)$ and a simple
closed curve $\g$ on $R$, one can decompose $\g$ as
a union of geodesic arcs $\g=\g_1\cup\dots\cup \g_k$
with $\g_i\in\Lambda(\l)$. The {\it intersection}
$\i(\mu,\g)$ is defined to be $\mu(\g_1)+\dots+\mu(\g_k)$.
Clearly, if $\g\simeq\g'$, then $\i(\mu,\g)=\i(\mu,\g')$.

Call $\Cc(R)$ the
set of nontrivial isotopy classes of simple
closed curves $\g$ contained in $R$.

\begin{remark}\label{rmk:topology-ML}
There is a map $\ML(R)\times\Cc(R)\lra\R_{\geq 0}$ given
by $(\mu,\g)\mapsto \i(\mu,\g)$. The induced
$\ML(R)\lra (\R_{\geq 0})^{\Cc(R)}$ is injective: identifying
$\ML(R)$ with its image, we can induce a
topology on $\ML(R)$ which is independent of the
hyperbolic structure on $R$ (see \cite{FLP}).
\end{remark}

A {\it $k$-system of curves}
$\bm{\g}=\{\g_1,\dots,\g_k\}\subset\Cc(R)$
is a subset of curves of $R$,
which admit disjoint representatives.

\begin{definition}
The {\it complex of curves} $\Cf(R)$ is the simplicial
complex whose $k$-simplices are $(k+1)$-systems of curves on $R$
(see \cite{harvey:bordification}).
\end{definition}

\begin{notation}
Given a simplicial complex $\mathfrak{X}$,
denote by $|\mathfrak{X}|$ its geometric realization.
It comes endowed with two natural topologies.
The {\it coherent topology} is the finest topology
that makes the realization of all simplicial maps continuous.
The {\it metric topology} is induced by the {\it path metric},
for which every $k$-simplex is isometric to the standard
$\Delta^k\subset\R^{k+1}$.
Where $|\mathfrak{X}|$ is not locally finite,
the metric topology is coarser than the coherent one.
Now on, we will endow all realizations with the metric topology,
unless differently specified.
\end{notation}

Clearly, there are continuous injective maps
$|\Cf(R)|\lra \Pr\ML(R)$ and $|\Cf(R)|\times\R_+\lra \ML(R)$.
Points in the image of the latter map are
called {\it multi-curves}.

\begin{definition}
Let $S$ is a compact hyperbolic surface with boundary.
A {\it geodesic lamination $\l$ on $S$} is a closed subset of $S$
foliated in geodesics that can meet $\pa S$ only perpendicularly;
equivalently, a $\s$-invariant geodesic lamination on its double $dS$.
A {\it measured lamination on $S$}
is a $\s$-invariant measured lamination $\l$ on $dS$.
%
\end{definition}

If $S$ has at least a boundary component or a marked point,
call $\Aa(S)$ the set of
all nontrivial isotopy classes of
simple arcs $\a\subset S$ with $\a^\circ\subset S^\circ$ and
endpoints at $\pa S$
or at the marked points of $S$.
A {\it $k$-system of arcs}
$\bm{A}=\{\a_1,\dots,\a_k\}\subset\Aa(S)$
is a subset of arcs of $S$,
that admit representatives which can intersect only at the
marked points.
The system $\bm{A}$ {\it fills} (resp. {\it quasi-fills})
$S$ if $S\setminus\bm{A}:=S\setminus\bigcup_{\a_i\in\bm{A}}\a_i$
is a disjoint union
of discs (resp. discs, pointed discs and annuli homotopic
to boundary components);
$\bm{A}$ is also called {\it proper} if
it quasi-fills $S$ (\cite{looijenga:cellular}).

\begin{definition}
The {\it complex of arcs} $\Af(S)$
of a surface $S$ with boundary and/or cusps
is the simplicial
complex whose
$k$-simplices are $(k+1)$-systems of arcs on $S$
(see \cite{harer:virtual}).
\end{definition}
We will denote by $\Ao(S)\subset\Af(S)$ the subset of
proper systems of arcs, which is the complement of
a lower-dimensional simplicial subcomplex, and by
$|\Ao(S)|\subset|\Af(S)|$ the locus
of weighted proper systems, which is open and dense.
%
%
\begin{notation}
If $\bm{A}=\{\a_1,\dots,\a_k\}\in\Af(S)$, then a point
$w\in|\bm{A}|\subset|\Af(S)|$ is a formal sum
$w=\sum_i w_i \a_i$ such that $w_i\geq 0$ and
$\sum_i w_i=1$, which can be also seen as a function
$w:\Aa(S)\rar\R$ supported on $\bm{A}$.
\end{notation}
We recall the following simple result.
\begin{lemma}
$|\Af(S)|/\G(S)$ is compact.
\end{lemma}
\begin{proof}
It is sufficient to notice the following facts:
\begin{itemize}
\item
$\G(S)$ acts on $\Af(S)$
\item
the above action may not be simplicial, but
it is on the second baricentric subdivision
$\Af(S)''$
\item
$\Af(S)/\G(S)$ is a finite set and so is
$\Af(S)''/\G(S)$.
\end{itemize}
\end{proof}

Clearly, for a hyperbolic surface $S$ with nonempty boundary
and no marked points, there are continuous injective maps
$|\Af(S)|\times\R_+\hra |\Cf(dS)|^{\s}\times\R_+\hra \ML(S)$
and $|\Af(S)|\hra |\Cf(dS)|^{\s}\hra \Pr\ML(S)$.


Notice that, if $R$ is a smooth compact surface without boundary,
then the multi-curves are dense in $\ML(R)$ and so the metric topology
on $|\Cf(R)|\times\R_+$ is stricly finer than the one coming
from $\ML(R)$. The situation for multi-arcs is different.

\begin{lemma}
If $S$ is a smooth surface with boundary,
the metric topology on $|\Af(S)|$ agrees with
the subspace topology induced by the inclusion
$|\Af(S)|\hra \Pr\ML(S)$. The image of $|\Ao(S)|$ is open but
the image of $|\Af(S)|$ is neither open nor close.
\end{lemma}
\begin{proof}
The image of $|\Ao(S)|$ is open by invariance of domain,
because $|\Ao(S)|$ is a topological manifold (for instance,
see \cite{hubbard-masur:foliations}). But if we consider an arc
$\a$ and a simple closed curve $\b$ disjoint from $\a$, then
the lamination $(1-t)[\a]+t[\b]\rar[\a]$ in $\Pr\ML(S)$ as $t\rar 0$,
which shows that the image of $|\Af(S)|$ is not open.
To show that it is not even closed, consider two
disjoint arcs $\{\a,\b\}\in\Af(S)$ and a simple closed curve $\g$
(possibly, a boundary component of $S$) such that $\a\cap\g=\emptyset$
and $i(\b,\g)=1$. Let $U\subset\ML(S)$ be an $\ML(S)$-neighbourhood
of $[\a]$ that contains the $|\Af(S)|$-ball of radius $2\e$ centered
at $\a$.
Consider the weighted arc systems $w^{(k)}=(1-\e)\a+
\e\tw_{k\g}(\b)$ in $|\Af(S)|$, where $\tw_{k\g}$ is
the $k$-uple Dehn twist along $\g$. Then, the sequence
$\{w^{(k)}\}$ is contained in $U$; moreover, it diverges in
$|\Af(S)|$, but it converges to $[\g]$ in $\Pr\ML(S)$.

To compare the topologies, pick $w=\sum w_i\a_i\in|\Af(S)|$
and $w^{(k)}=\sum v_j^{(k)}\b_j^{(k)}\in|\Af(S)|$
such that $w^{(k)}\rar w$ in $\Pr\ML(S)$.
Complete $\bm{A}=\{\a_i\}$ to a maximal system of arcs $\bm{A}'=\{\a_i\}\cup\{\a'_h\}$
and define $w'=w+\d\sum_h \a'_h$, where $\d=\mathrm{min}_i w_i$.

For every $k$, write $w^{(k)}$ as a sum $\tilde{w}^{(k)}+\hat{w}^{(k)}$
of two nonnegative multi-arcs
in such a way that all arcs in the support of $\hat{w}^{(k)}$
cross $\bm{A}'$ and that $i(\tilde{w}^{(k)},w')=0$.

Let $t_k$ be the sum of the weights on
$\hat{w}^{(k)}$, so that
$d(w,w^{(k)})\leq d(w,\tilde{w}^{(k)})+t_k$,
where $d$ is the path metric on $|\Af(S)|$. Because
\[
t_k\d \leq i(\hat{w}^{(k)},w')=i(w^{(k)},w')\rar
i(w,w')=0
\]
it follows that $t_k\rar 0$.
Moreover, $\tilde{w}^{(k)}$ has support contained in $\bm{A}'$
and the result follows.
\end{proof}
\end{subsection}
%
%
\begin{subsection}{Thurston's compactification}
Let $R$ be a closed hyperbolic surface and let $\Cc(R)$ be the
set of nontrivial isotopy classes of simple closed curves of $R$.
The map $\ell:\Teich(R)\times\Cc(R)\lra \R_+$, that
assigns to $([f:R\rar R'],\g)$ the length of the
geodesic representative for $f(\g)$ in the hyperbolic metric
of $R'$, induces an embedding
$I:\Teich(R)\hra (\R_+)^{\Cc(R)}$, where $(\R_+)^{\Cc(R)}$
is given the product topology (which is the same as the
weak$^*$ topology on $L^{\infty}(\Cc(R))$).

\begin{theorem}[Thurston \cite{FLP}]\label{prop:thurston}
The composition $[I]:\Teich(R)\hra\R_+^{\Cc(R)}\rar
\Pr(\R_+^{\Cc(R)})$ is an embedding with
relatively compact image and
the boundary of $\Teich(R)$ inside $\Pr(\R_+^{\Cc(R)})$
is exactly $\Pr\ML(R)$.
\end{theorem}

Let $S$ be a hyperbolic surface with boundary
and no cusps.
The doubling map
$\Cc(S)\cup\Aa(S)\hra \Cc(dS)$ identifies
$\Cc(S)\cup\Aa(S)$ with $\Cc(dS)^{\s}$.

\begin{corollary}\label{cor:thurston}
$\Teich(S)$ embeds in
$\Pr(\R_+^{\Cc(S)\cup\Aa(S)})$ and
its boundary is $\Pr\ML(S)$.
\end{corollary}

For $S$ a hyperbolic surface with no cusps,
$\Tbar^{Th}(S):=\Teich(S)\cup\Pr\ML(S)$
is called {\it Thurston's compactification} of $\Teich(S)$.
Notice that the doubling map $S\hra dS$ induces
a closed embedding
$D:\Tbar^{Th}(S)\hra\Tbar^{Th}(dS)$.
\end{subsection}
\begin{subsection}{Weil-Petersson metric}
Let $S$ be a hyperbolic surface with (possibly empty)
geodesic boundary $\pa S=C_1\cup\dots\cup C_n$,
and let $[f:S\rar \Si]$ a point of $\Teich(S)$.

Define $\Qq_{\Si}$ to be the {\it real}
vector space of holomorphic quadratic differentials $q(z)dz^2$
whose restriction to $\pa \Si$ is real. Similarly, define
the real vector space of harmonic Beltrami differentials
as $\Bel_{\Si}:=
\{\mu=\mu(z)d\bar{z}/dz=\bar{\varphi}\,ds^{-2}\,|\,\varphi\in\Qq_{\Si}\}$,
where $ds^{2}$ is the hyperbolic metric on $\Si$.

It is well-known that $T_{[f]}\Teich(S)$ can be identified to $\Bel_{\Si}$
and, similarly, $T^*_{[f]}\Teich(S)\cong\Qq_{\Si}$.
The natural coupling is given by
\[
\xymatrix@R=0in{
\Bel_{\Si}\times\Qq_{\Si} \ar[rr] && \C \\
(\mu,\varphi) \ar@{|->}[rr] && \dis\int_\Si \mu\varphi
}
\]

\begin{definition}
The {\it Weil-Petersson pairing} on $T_{[f]}\Teich(S)$ is
defined as
\[
h(\mu,\nu):=\int_{\Si} \mu\bar{\nu}\,ds^2 \qquad\text{with}\ \mu,\nu
\in\Bel_{\Si}
\]
Writing $h=g+i\omega$, we call $g$ the {\it Weil-Petersson
Riemannian metric} and $\omega$ the {\it Weil-Petersson form}.
For $T^*_{[f]}\Teich(S)$, we similarly have
$\dis h^{\vee}(\varphi,\psi):=\int_{\Si} \varphi\bar{\psi}\,ds^{-2}$
with $\varphi,\psi\in\Qq_{\Si}$.
The {\it Weil-Petersson Poisson structure} is $\eta:=\mathrm{Im}(h^{\vee})$.
\end{definition}

It follows from the definition that 
the doubling map
$D:\Teich(S)\lra\Teich(dS)$
is a homothety of factor $2$ onto a real Lagrangian
submanifold of $\Teich(dS)$.

From Wolpert's work \cite{wolpert:symplectic}, we learnt
that $\dis\omega=\sum_{i=1}^N d\ell_i\wedge d\tau_i$, where
$(p_1,\dots,p_n,\ell_1,\tau_1,\dots,\ell_N,\tau_N)$ are
Fenchel-Nielsen coordinates, and so $\omega$
is degenerate whenever $S$ has boundary. In this case,
the symplectic leaves (which we will also denote by $\Teich(S)(\up)$)
are exactly the fibers
$\Ll^{-1}(\up)$ of $\Ll$,
which are not totally geodesic subspaces for $g$ (unless $p_1=\dots=p_n=0$
and the boundary components degenerate to cusps).

Using the Weil-Petersson metric, the cotangent
space to $\Teich(S)(\up)$ at
$[f:S\rar\Si]$ can be identified with $\dis(dp_1\oplus
\dots\oplus dp_n)^{\perp}\subset T^*_{[f]}\Teich(S)$.
It follows from \cite{wolf:harmonic} that
the elements of $\dis(dp_1\oplus
\dots\oplus dp_n)^{\perp}$
are those $\varphi\in\Qq_{\Si}$
such that $\dis\int_{C_i}\varphi |ds|^{-1}=0$ for all
$i=1,\dots,n$. Similarly, the tangent space $T_{[f]}\Teich(S)(\up)$
is the subspace of those $\mu\in\Bel_{\Si}$
such that $\dis\int_{C_i}\mu|\l|=0$ for $i=1,\dots,n$.

\begin{remark}
$\Teich(S)$ is naturally a complex manifold if $S$ is closed.
In this case, $\omega$ and $\eta$ are nondegenerate and
the Weil-Petersson metric is K\"ahler (see \cite{ahlfors:remarks}).
\end{remark}
\end{subsection}
%
%
\begin{subsection}{Augmented Teichm\"uller space}
Let $S$ be a hyperbolic surface with geodesic boundary
$\pa S=C_1\cup\dots\cup C_n$ and no cusps.
An {\it $S$-marked stable surface $\Si$} is 
a hyperbolic surface possibly with geodesic boundary
components, cusps and nodes
plus an isotopy
class of maps $f:S\rar\Si$ that may shrink some boundary components of $S$
to cusps of $\Si$, some loops of $S$ to the nodes of $\Si$
and is an oriented diffeomorphism elsewhere.

We say that $f_1:S\lra \Si_1$ and $f_2:S\lra\Si_2$
are equivalent if there exists an isometry $h:\Si_1\lra \Si_2$
such that $h\circ f_1$ is homotopic to $f_2$.
The {\it augmented Teichm\"uller space} $\Tbar^{WP}(S)$ is
the set of stable $S$-marked surfaces up to equivalence
(see \cite{bers:degenerating}).
Clearly, $\Teich(S)\subset\eTeich(S)\subset\Tbar^{WP}(S)$.

To describe the topology of $\Tbar^{WP}(S)$
around a stable surface $[f:S\rar\Si]$
with $k$ cusps and $d$ nodes,
choose a system of curves $\{C_1,\dots,C_n,\g_1,\dots,\g_N\}$
on $S$ (with $N=3g-3+n$) {\it adapted to $f$}, i.e.
such that $f^{-1}(\nu_j)=\g_j$
for each of the nodes $\nu_1,\dots,\nu_d$ of $\Si$.
Clearly, the Fenchel-Nielsen
coordinates $(p_1,\dots,p_n,\ell_1,\tau_1,\dots,\ell_N,\tau_N)$
extend over $[f]$, with the exception of
$\tau_1(f),\dots,\tau_d(f)$, which are not defined
(see \cite{abikoff:book} for more details on the Fenchel-Nielsen
coordinates).

We declare that the sequence $\{f_m:S\rar\Si_m\}\subset\Tbar^{WP}(S)$
converges to $[f]$ if $p_i(f_m)\rar p_i(f)$ for $1\leq i\leq n$,
$\ell_j(f_m)\rar\ell_j(f)$ for $1\leq j\leq N$ and
$\tau_j(f_m)\rar\tau_j(f)$ for $d+1\leq j\leq N$.
By definition, the boundary length map
extends with continuity to
$\Ll:\Tbar^{WP}(S)\lra\R_{\geq 0}^n$
and we call $\Tbar^{WP}(S)(\up)$
the fiber $\Ll^{-1}(\up)$.
We will write $\bm{p}$ for $p_1+\dots+p_n$
and $\bm{\Ll}(f)$ for the $L^1$ norm of $\Ll(f)$.

The cotangent cone $T^*_{[f:S\rar\Si]}\Tbar^{WP}(S)$
(with the analytic smooth structure)
can be identified
to the space $\Qq_{\Si}$ of holomorphic quadratic differentials
$\varphi$ on $\Si$ that are real at $\pa\Si$
and that have (at worst) double poles at the cusps with negative
quadratic residues and (at worst) double poles
at the nodes with the same quadratic residue
on both branches (see \cite{wolf:harmonic}).
Those $\varphi$ which do not
have a double pole at the cusp $f(C_i)$ (resp. at the
node $f(\g_j)$) are perpendicular
to $dp_i$ (resp. to $d\ell_j$).
Similarly, $T_{[f]}\Tbar^{WP}(S)$ can be identified
to the space $\Bel_{\Si}$ of harmonic Beltrami differentials
$\mu=\bar{\varphi}ds^{-2}$, where $\varphi\in\Qq_{\Si}$
and $ds^2$ is the hyperbolic metric.

Notice that the Weil-Petersson metric diverges
in directions transverse
to $\pa\Teich(S)$.
However, the divergence is so mild that $\pa\Teich(S)$
is at finite distance (see \cite{masur:incomplete}).
In fact, for every $\up\in\R_{\geq 0}^n$ the augmented
$\Tbar^{WP}(S)(\up)$ is the completion of $\Teich(S)(\up)$
with respect to the Weil-Petersson metric (it follows from
the $\G(S)$-invariance of the metric, its compatibility with the
doubling map $D$ and the compactness of the Deligne-Mumford
moduli space \cite{deligne-mumford:irreducibility}).
\begin{remark}
According to our definition, if $S$ has nonempty boundary,
then $\Tbar^{WP}(S)$ is not WP-complete and in fact the image
of $\Tbar^{WP}(S)$ inside $\Tbar^{WP}(dS)$ through the doubling map
is not close because it misses thoses surfaces with boundaries of
infinite length.
\end{remark}
We recall here a criterion of convergence in $\Tbar^{WP}(S)$
that will be useful later.
\begin{proposition}[\cite{mondello:criterion}]\label{prop:convergence}
Let $[f:S\rar(\Si,g)]\in\Tbar^{WP}(S)$ and call $\g_1,\dots,\g_d$
the simple closed curves of $S$ that are contracted to a point by $f$,
and let $\{f_m:S\rar(\Si_m,g_m)\}$ be a sequence of points in $\Tbar^{WP}(S)$.
Denote by $V_{\g_i}(f_m)$ a standard collar (of fixed width) of the hyperbolic
geodesic homotopic to $f_m(\g_i)\subset\Si_m$ and set
$V_i=V_{\g_i}(f)$ and $\Si^\circ:=\Si\setminus(V_1\cup\dots\cup V_k)
\subset\Si_{sm}$.
The following are equivalent:
\begin{enumerate}
\item 
$[f_m]\rar [f]$ in $\Tbar^{WP}(S)$
\item 
$\ell_{\g_i}(f_m)\rar 0$ and
there exist representatives $\tilde{f}_m\in[f_m]$ such that
$\dis(f\circ\tilde{f}_m^{-1})\Big|_{V_{\g_i}(f_m)}$ is standard
and $(\tilde{f}_m\circ f^{-1})^*(g_m)\rar g$ uniformly on
$\Si^\circ$
\item 
$\exists\tilde{f}_m\in[f_m]$ such that
the metrics $(\tilde{f}_m\circ f^{-1})^*(g_m)\rar g$ uniformly
on the compact subsets of $\Si_{sm}$ 
\item 
$\ell_{\g_i}(f_m)\rar 0$ and
$\exists\tilde{f}_m\in[f_m]$ such that
$\dis(f\circ\tilde{f}_m^{-1})\Big|_{V_{\g_i}(f_m)}$ is standard
and
$\dis(\tilde{f}_m\circ f^{-1})\Big|_{\Si^\circ}$ is
$K_m$-quasiconformal with $K_m\rar 1$
\item 
$\exists\tilde{f}_m\in[f_m]$ such that,
for every compact subset $F\subset\Si_{sm}$,
the homeomorphism $\dis(\tilde{f}_m\circ f^{-1})\Big|_F$ is
$K_{m,F}$-quasiconformal
and $K_{m,F}\rar 1$. 
\end{enumerate}
\end{proposition}
We denoted by $g_m$ the hyperbolic metric on $\Si_m$
and by $\Si_{sm}$ the locus of $\Si$ on which $g$ is smooth
(namely, $\Si$ with cusps and nodes removed).
By ``standard collar'' of width $t$ of a boundary component $\g$ (resp.
an internal curve $\g$),
we meant an annulus of the form $A_t(\g)$
(resp. the union of the two annuli isometric to $A_t(\g)$
that bound $\g$), as provided
by the following celebrated result.

\begin{lemma}[Collar lemma, \cite{keen:collar}-
\cite{matelski:collar}]\label{lemma:collar}
For every simple closed geodesic $\g\subset\Si$
in a hyperbolic surface and for every ``side'' of $\g$,
and for every $0<t\leq 1$,
there exists an embedded hypercycle $\g'$ parallel
to $\g$ (on the prescribed side) such that
the area of the annulus $A_t(\g)$ enclosed by $\g$ and $\g'$
is $t\ell/2\sinh(\ell/2)$.
For $\ell=0$, the geodesic $\g$ must be intended
to be a cusp and $\g'$ a horocycle of length $t$.
Furthermore, all such annuli (corresponding
to distinct geodesics and sides) are disjoint.
\end{lemma}
{\it Standard maps} between annuli or between pair of pants
are defined in \cite{mondello:criterion}.
\end{subsection}
\begin{subsection}{The moduli space}
Let $S$ be a hyperbolic surface of genus $g$
with geodesic boundary
components $C_1,\dots,C_n$ and no cusps.
The augmented Teichm\"uller space $\Tbar^{WP}(S)$
(as well as Thurston's compactification $\Tbar^{Th}(S)$)
carries a natural
right action of the
group $\mathrm{Diff}_+(S)$
of orientation-preserving diffeomorphisms of $S$
that send $C_i$ to $C_i$ for every $i=1,\dots,n$.
\[
\xymatrix@R=0in{
\Tbar^{WP}(S)\times\mathrm{Diff}_+(S) \ar[rr] && \Tbar^{WP}(S)\\
([f:S\rar\Si],h) \ar@{|->}[rr] && [f\circ h:S\rar \Si]
}
\]
Clearly, the action is trivial on the connected component
$\mathrm{Diff}_0(S)$ of the identity.

\begin{definition}
The {\it mapping class group} of $S$ is the quotient
\[
\G(S):=\mathrm{Diff}_+(S)/\mathrm{Diff}_0(S)=
\pi_0\mathrm{Diff}_+(S)
\]
The quotient $\Mbar(S):=\Tbar^{WP}(S)/\G(S)$
is the {\it moduli space of stable hyperbolic
surfaces} of genus $g$ with $n$ (ordered) boundary components.
\end{definition}

The quotient map $\pi:\Tbar^{WP}(S)\lra\Mbar(S)$ can be identified
with the forgetful map $[f:S\rar\Si]\mapsto [\Si]$.
Moreover, we can identify the
stabilizer $\mathrm{Stab}_{[\Si]}(\G(S))$ with the group $\Iso_+(\Si)$
of orientation-preserving isometries of $\Si$, which is finite.

$\Mbar(S)$ can be given a natural structure of orbifold (with corners),
called {\it Fenchel-Nielsen smooth structure}.
Let $[f:S\rar\Si]$ be a point of $\Tbar^{WP}(S)$
and let $(p_1,\dots,p_n,\ell_1,\tau_1,\dots,\ell_N,\tau_N)$
be Fenchel-Nielsen coordinates adapted to $f$.
A local chart for (the Fenchel-Nielsen smooth structure of)
$\Mbar(S)$ around $[\Si]$ is given by
\[
\xymatrix@R=0in{
\R_{\geq 0}^n\times\C^{3g-3} \ar[rr] && \Mbar(S) \\
(p,z)
\ar@{|->}[rr] &&
\Si(p,z)
}
\]
where $[f':S\rar\Si(p,z)]$ is the
point of $\Tbar^{WP}(S)$ with coordinates
$(p_1,\dots,p_n,\ell_1,\tau_1,\dots,\ell_N,\tau_N)$
with $\ell_j=|z_j|$ and $\tau_j=|z_j|\arg(z_j)/2\pi$.

\begin{remark}
As shown by Wolpert \cite{wolpert:geometry}, the smooth
structure at $\pa\M(S)=\Mbar(S)\setminus\M(S)$
coming from Fenchel-Nielsen
coordinates and the one coming from algebraic geometry
(see \cite{deligne-mumford:irreducibility}) are not
the same.
\end{remark}

We can identify the (co)tangent space
to $\Teich(S)$
(with the analytic structure)
at $[f:S\rar\Si]$ with the (co)tangent space
to $\M(S)$ at $[\Si]$. It follows by its very definition that
the Weil-Petersson metric and the boundary lengths map
descends to $\M(S)$ and that
$\Mbar(S)(\up)$ is the metric completion of $\M(S)(\up)$
for every $\up\in\R_{\geq 0}^n$.
\end{subsection}
%
%
\end{section}
%
%
\begin{section}{Triangulations}\label{sec:triangulations}
%
%
\begin{subsection}{Systems of arcs and widths}
Let $S$ be a hyperbolic surface of genus $g$ with boundary
components $C_1,\dots,C_n$ and no cusps, and let $\bm{A}=\{\a_1,\dots,
\a_N\}\in\Ao(S)$ a maximal system of arcs on $S$
(so that $N=6g-6+3n$).

Fix a point $[f:S\rar\Si]$ in $\Teich(S)$. For every
$i=1,\dots,N$ there exists a unique geodesic
arc on $\Si$ in the isotopy class of $f(\a_i)$
that meets $\pa\Si$ perpendicularly and
which we will still denote by $f(\a_i)$:
call $a_i=\ell_{\a_i}(f)$ its length
and let $s_i=\cosh(a_i/2)$.
Notice that $\{f(\a_i)\}$ decomposes $\Si$ into a disjoint
union of right-angles hexagons
$\{H_1,\dots,H_{4g-4+2n}\}$, so that the following is immediate
(see also \cite{ushijima:decomposition}, \cite{mondello:poisson}).

\begin{lemma}
The maps $a_{\bm{A}}:\Teich(S)\lra \R_+^{\bm{A}}$
and $s_{\bm{A}}:\Teich(S)\lra\R_+^{\bm{A}}$ given
by $a_{\bm{A}}=(a_1,\dots,a_N)$ 
and $s_{\bm{A}}=(s_1,\dots,s_N)$ are real-analytic
diffeomorphisms.
\end{lemma}

Let $H$ be such a right-angled hexagon and let
$(\ora{\a_i},\ora{\a_j},\ora{\a_k})$ the cyclic set of
oriented arcs that bound $H$, so that $\pa H=\ora{\a_i}\ast
\ora{\a_j}\ast\ora{\a_k}$. If $\ora{\a_x},\ora{\a_y}$ are
oriented arcs with endpoint on the same boundary component $C$,
denote by $d(\ora{\a_x},\ora{\a_y})$ the length of the
portion of $C$ running from the endpoint of $\ora{\a_x}$
to the endpoint of $\ora{\a_y}$ along the positive direction
of $C$.

Define $\dis w_{\bm{A}}(\ora{\a_i})=\frac{1}{2}[d(\ora{\a_i},\ola{\a_j})+
d(\ora{\a_k},\ola{\a_i})-d(\ora{\a_j},\ola{\a_k})]$,
where $\ola{\a_x}$ the oriented
arc obtained from $\ora{\a_x}$ by switching its orientation.

\begin{definition}
The {\it $\bm{A}$-width} of $\a_i$ is $w_{\bm{A}}(\a_i)=
w_{\bm{A}}(\ora{\a_i})+w_{\bm{A}}(\ola{\a_i})$.
\end{definition}
In \cite{luo:decomposition}), Luo calls the width ``E-invariant''.\\
\end{subsection}
%
%
\begin{subsection}{The $t$-coordinates}
Let $S$ be a surface as in the previous section.

\begin{definition}
The {\it t(ransverse)-length} of an arc $\a$ at $[f]$ is
$t_{\a}(f):=T(\ell_{\a}(f))$,
where
$\dis T(x):=2\mathrm{arcsinh}\left(\frac{1}{\sinh(x/2)}\right)$.
\end{definition}

Notice that $T(x):[0,+\infty]\rar[0,+\infty]$ is
decreasing function of $x$
(similar to the width of the collar of a closed curve
of length $x$ provided by Lemma~\ref{lemma:collar}).
Moreover, $T$ is involutive, $T(x)\approx 4e^{-x/2}$
as $x\rar\infty$ and $T(x)\approx 2\log(4/x)$ as $x\rar 0$.

Back to the $t$-length, the following lemma reduces
to a statement about hyperbolic hexagons with right angles.

\begin{lemma}\label{lemma:t}
For every maximal system of arcs $\bm{A}$
\[
t_{\bm{A}}:
\xymatrix@R=0in{
\eTeich(S) \ar[rr] && \R_{\geq 0}^{\bm{A}} \\
[f] \ar@{|->}[rr] && (t_{\a_1}(f),\dots,t_{\a_N}(f))
}
\]
is a continuous map that restricts to a real-analytic diffeomorphism
$\Teich(S)\rar\R_+^{\bm{A}}$.
Moreover, for every $[f]\in\eTeich(S)$ with $\Ll(f)\neq 0$,
there exists an $\bm{A}$ such that $t_{\bm{A}}$ is a system
of coordinates around $[f]$.
\end{lemma}

Consequently,
the $t$-lengths map $\eTeich(S)\times\Aa(S)\lra\R_{\geq 0}$
defined as $(f,\a)\mapsto t_{\a}(f)$
gives an injection
\[
j:\xymatrix@R=0in{
\Teich(S)\ar@{^(->}[r] &
\Pr(\Aa(S))\times[0,\infty] \\
[f] \ar@{|->}[r] &
([t_{\bullet}(f)],\| t_\bullet(f)\|_{\infty})
}
\]
where $L^{\infty}(\Aa(S))$ is the $\R_+$-cone of the
bounded maps $t:\Aa(S)\rar\R_{\geq 0}$ and
$\Pr(\Aa(S))$ is its projectivization.

Notice that $\Pr(\Aa(S))$ has a metric induced by
the unit sphere of $L^{\infty}(\Aa(S))$
and that $\G(S)$
acts on $\Pr(\Aa(S))$ permuting some coordinates.
Thus, $\Pr(\Aa(S))\times[0,\infty]$
has a $\G(S)$-invariant metric.

\begin{lemma}
$j$ is continuous.
\end{lemma}

We proof is included in that of
Proposition~\ref{prop:j-continuous}.

\begin{definition}
Call {\bf bordification of arcs} the closure
$\Tbar^a(S)$ of $\Teich(S)$
inside $\Pr(\Aa(S))\times[0,\infty]$.
By ``finite part'' of $\Tbar^a(S)$ we will
mean $\Tbar^a(S)\cap \Pr(\Aa(S))\times[0,\infty)$.
Call {\bf compactification of arcs} the quotient
$\Mbar^a(S):=\Tbar^a(S)/\G(S)$.
\end{definition}

We will give an explicit description of the boundary
points in $\Tbar^a(S)$ and we will show that
$\Mbar^a(S)$ is Hausdorff and compact.
\end{subsection}
%
%
\begin{subsection}{The spine construction}
Let $S$ be a hyperbolic surface with geodesic boundary
$\pa S=C_1\cup\dots\cup C_n$ and no cusps and let $[f:S\rar\Si]$
be a point in $\Teich(S)$.

The {\it valence} $\mathrm{val}(p)$ of a point $p\in\Si$
is the number
of paths from $p$ to $\pa\Si$ of minimal length.

\begin{definition}
The {\it spine} of $\Si$ is the locus $\mathrm{Sp}(\Si)$
of points of $\Si$ of valence at least $2$.
\end{definition}

One can easily show that
$\mathrm{Sp}(\Si)=V\cup E$ is a one-dimensional
CW-complex embedded in $\Si$, where $V=\mathrm{val}^{-1}([3,\infty))$
is a finite set of points, called {\it vertices}, and
$E=\mathrm{val}^{-1}(2)$ is a disjoint union of finitely many
(open) geodesic arcs, called {\it edges}.

For every edge $E_i\subset E$ of $\mathrm{Sp}(\Si)$, we can define
a {\it dual arc} $\a_i$ in the following way. Pick $p\in E_i$
and call $\g_1$ and $\g_2$ the two paths that join $p$ to $\pa\Si$.
Then $\a_i$ is the shortest arc in the homotopy
class (with endpoints on $\pa \Si$) of $\g_1^{-1}\ast\g_2$.
Let the {\it spinal arc system}
$\bm{A}_{sp}(\Si)$ be the system of arcs dual to the
edges of $\mathrm{Sp}(\Si)$, which is proper because
$\Si$ retracts by deformation onto $\mathrm{Sp}(\Si)$ just
flowing away from the boundary.

Even if the spinal arc system is not maximal,
widths $w_{sp}$ can be associated to $\bm{A}_{sp}(\Si)$
in the following way.
For every oriented arc $\ora{\a_i}\in\bm{A}_{sp}(\Si)$
ending at $y_i\in C_m$,
orient the dual edge $E_i$ in such a way that $(\ora{E_i},\ora{\a_i})$
is positively oriented and call $v$ the starting point of $\ora{E_i}$.
Every point of $E_i$ has exactly two projections,
that is two closest points in $\pa\Si$:
the endpoint of $\ora{\a_i}$ selects only one of these, which belongs to $C_m$.
Call $v'\in C_m$ the projection of $v$ determined by $\ora{\a_i}$.
Define $w_{sp}(\ora{\a_i})$ to be the distance with sign
$d_{C_m}(y_i,v')$ along $C_m$, which is certainly
positive if $\a_i$ and $E_i$ intersect.
Clearly, the sum
$w_{sp}(\a_i)=w_{sp}(\ora{\a_i})+w_{sp}(\ola{\a_i})$ is always positive,
being the length of either of the two projections of $E_i$.

\begin{example}
In Figure~\ref{fig:spine}, we have $\ora{\a_i}=\ora{z_i y_i}$, $v'=f_k$
and $w_{sp}(\ora{\a_i})>0$.
\end{example}

\begin{theorem}[Ushijima \cite{ushijima:decomposition}]
Given a hyperbolic surface with nonempty boundary $\Si$,
let $\Af(\Si)_+$ be the set of all maximal systems of arcs $\bm{A}$ 
such that $w_{\bm{A}}(\a_i)\geq 0$ for all $\a_i\in\bm{A}$.
Then $\Af(\Si)_+$ is nonempty and the intersection of all systems in
$\Af(\Si)_+$ is exactly $\bm{A}_{sp}(\Si)$.
Moreover, $w_{sp}(\a)=w_{\bm{A}}(\a)>0$ for all $\a\in\bm{A}_{sp}(\Si)$
and all $\bm{A}\in \Af(\Si)_+$.
\end{theorem}
\begin{center}
\begin{figurehere}
\psfrag{zi}{$z_i$}
\psfrag{yi}{$y_i$}
\psfrag{mi}{$m_i$}
\psfrag{fi}{$f_i$}
\psfrag{zj}{$z_j$}
\psfrag{yj}{$y_j$}
\psfrag{mj}{$m_j$}
\psfrag{fj}{$f_j$}
\psfrag{zk}{$z_k$}
\psfrag{yk}{$y_k$}
\psfrag{mk}{$m_k$}
\psfrag{fk}{$f_k$}
\psfrag{v}{$v$}
\psfrag{ai}{$\ora{\a_i}$}
\psfrag{Ei}{$\ora{E_i}$}
\psfrag{gi}{$\g_i$}
\psfrag{gj}{$\g_j$}
\psfrag{gk}{$\g_k$}
\includegraphics[width=0.9\textwidth]{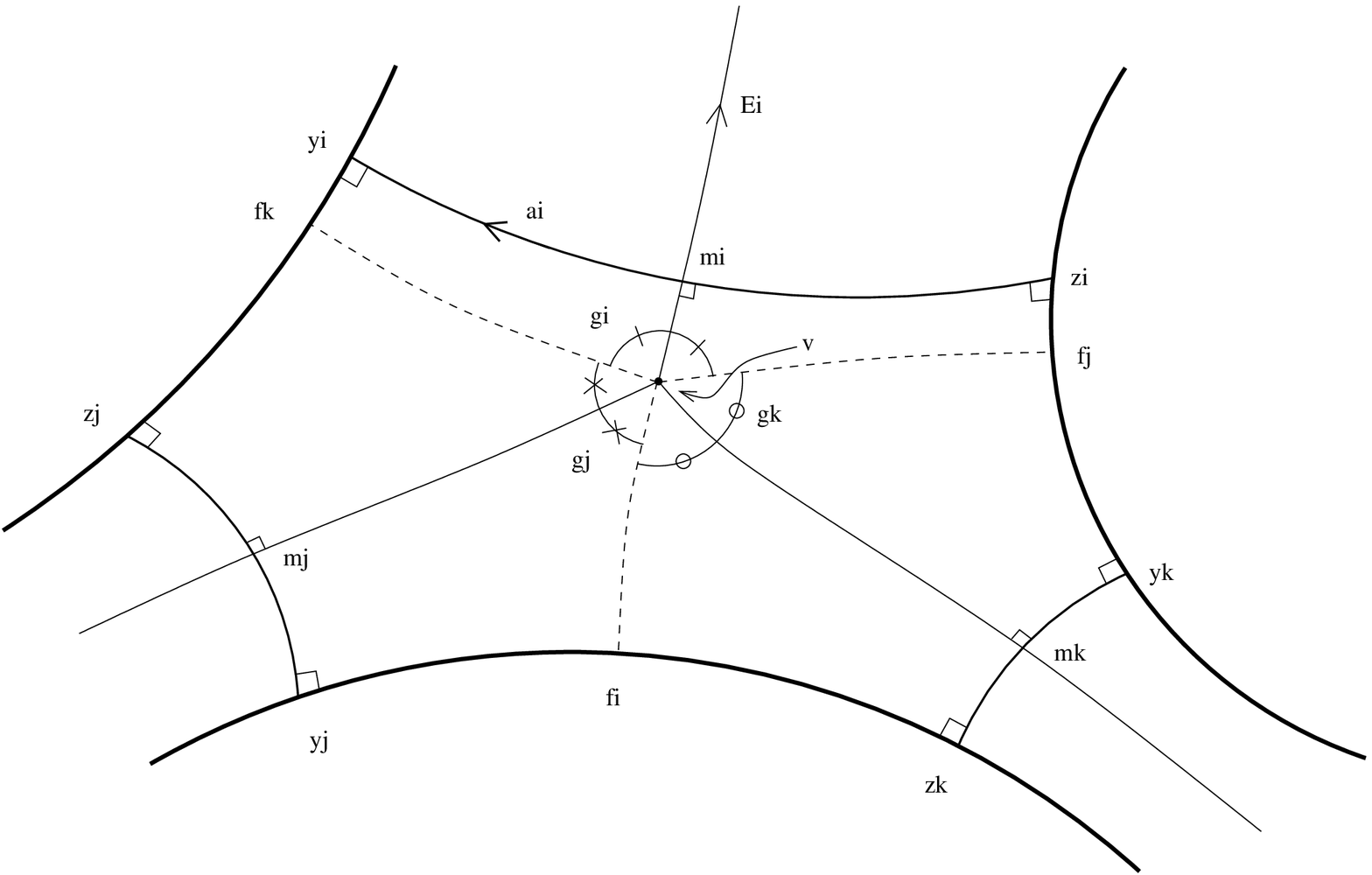}

\myCaption{Geometry of the spine close to a trivalent vertex}
\label{fig:spine}
\end{figurehere}
\end{center}
\begin{remark}
Let $H\subset\Si$ be a right-angled hexagon, bounded by
$(\ora{\a_i},\ora{\a_j},\ora{\a_k})$,
and call $\g_i=\g(\ora{\a_i})$ as in Figure~\ref{fig:spine}.
An easy computation \cite{mondello:poisson} shows that
\begin{equation}\label{eq:hex}
\sinh(w_{sp}(\ora{\a_i}))\sinh(a_i/2)=\cos(\g_i)=
\frac{s_j^2+s_k^2-s_i^2}{2s_j s_k}
\end{equation}
and so
$\dis \sinh(w_{sp}(\ora{\a_i}))=\frac{s_j^2+s_k^2-s_i^2}
{\dis 2s_j s_k \sqrt{s_i^2-1}}$
and $\dis w_{sp}(\ora{\a_i})\leq \frac{1}{2}t_{\a_i}$.
\end{remark}

\begin{theorem}[Luo \cite{luo:decomposition}]\label{thm:luo}
Given a hyperbolic surface with boundary $S$ and no cusps, the map
\[
\bm{W}\ :
\xymatrix@R=0in{
\Teich(S) \ar[rr] && |\Ao(S)|\times\R_+ \\
[f:S\rar\Si] \ar@{|->}[rr] && f^*w_{sp}
}
\]
is a $\Gamma(S)$-equivariant homeomorphism.
\end{theorem}

Notice that the construction extends to $\eTeich(S)
\setminus\eTeich(S)(0)$ but the locus $\eTeich(S)(0)$
of surfaces with $n$ cusps is problematic,
because the function ``distance from the boundary $\pa\Si$''
diverges everywhere on $\Si$.
This can be easily fixed by
considering the real blow-up $\mathrm{Bl}_0\eTeich(S)$
of $\eTeich(S)$ along $\eTeich(0)$.
The exceptional locus
can be identified to the space of {\it projectively decorated
surfaces} \cite{penner:decorated}, that is of
couples $([f:S\rar\Si],\up)$,
where $[f]$ is an $S$-marked hyperbolic surface with $n$
cusps and $\up\in\Delta^{n-1}\cong\Pr(\R_{\geq 0}^n)$
is a ray of weights (the {\it decoration}).

We have the following two simple facts.

\begin{lemma}[\cite{mondello:poisson}]\label{lemma:lambda}
For every maximal system of arcs $\bm{A}$
(of cardinality $N=6g-6+3n$), the associated $t$-lengths
extend to real-analytic map
\[
t_{\bm{A}}:
\xymatrix@R=0in{
\mathrm{Bl}_0\eTeich(S) \ar[rr] && \Delta^{N-1}\times[0,\infty)
}
\]
On the exceptional locus,
the projectivized $t$-lengths are inverses to the projectivized $\l$-lengths
(defined by Penner in \cite{penner:decorated}). Thus, $t_{\bm{A}}$
gives a system of coordinates on
$[\eTeich(S)(0)\times(\Delta^{n-1})^\circ]
\cup\Teich(S)$.
Moreover, for every $(f,\up)\in\eTeich(S)(0)\times\pa\Delta^{n-1}$
there exists an $\bm{A}$ such
that $t_{\bm{A}}$ gives a chart around $(f,\up)$.
\end{lemma}

\begin{theorem}[\cite{mondello:poisson}]\label{prop:smooth}
The map $\bm{W}$ extends to a $\G(S)$-equivariant homeomorphism
\[
\extended{\bm{W}}\ :
\xymatrix@R=0in{
\mathrm{Bl}_0\eTeich(S) \ar[rr] && |\Ao(S)|\times[0,\infty)
}
\]
On the exceptional locus, the projectivized
widths coincide with the projectived simplicial coordinates
(defined by Penner in \cite{penner:decorated})
and so there $\extended{\bm{W}}$ coincides with Penner's homeomorphism.
\end{theorem}
\end{subsection}
%
%
%
\begin{subsection}{Spines of stable surfaces}\label{sec:stable}
Notice that the spine construction extends to stable hyperbolic
surfaces $\Si$ (blowing up the locus of surfaces with $n$ cusps),
discarding the components of $\Si$ where the distance from $\pa\Si$
is infinite.
However, the weighted arc system
we can produce does not allow to reconstruct the full surface,
but just a visible portion of it.

\begin{definition}
Let $\Si$ a stable hyperbolic surface $\Si$ with boundary
(and possibly cusps) or let $(\Si,\up)$ be a stable
(projectively) decorated surface.
A component of $\Si$ is called {\it visible}
if it contains a boundary circle
or a positively weighted cusp; otherwise, it is called {\it invisible}.
Denote by $\Si_+$ the visible subsurface of $\Si$, that is the
union of the smooth points
of all visible components and by $\Si_-$ the invisible subsurface.
Two points $[f_1:S\rar\Si_1]$
and $[f_2:S\rar\Si_2]$ of $\That(S):=\mathrm{Bl}_0\Tbar^{WP}(S)$
are {\it visibly equivalent} $[f_1]\sim_{vis}[f_2]$
if there exists a third point $[f:S\rar\Si]$ and
maps $h_i:\Si\rar \Si_i$ for $i=1,2$ such that
$h_i$ restricts to an isometry $\Si_+\rar\Si_{i,+}$
and $h_i\circ f\simeq f_i$ for $i=1,2$.
\end{definition}

The spine $\Sp(\Si)$ of a stable hyperbolic surface $\Si$
with geodesic boundary (or with weighted cusps)
can only be defined inside $\ol{\Si}_+$,
so that its dual system of arcs $\bm{A}_{sp}(\Si)$
will be contained in $\ol{\Si}_+$ too.
Given a marking $[f:S\rar\Si]$, we will write
$S_+=f^{-1}(\Si_+)$ and $S_-=f^{-1}(\Si_-)$,
so that $\ol{S}_+$ will be the maximal subsurface of $S$
(unique up to isotopy),
quasi-filled by
$f^*\bm{A}_{sp}(\Si)$,
which carries positive weights $f^*w_{sp}$.

Conversely, given a system of arcs $\bm{A}\in \Af(S)$,
the {\it visible subsurface} $S_+$ associated to $\bm{A}$
is the isotopy class of maximal open subsurfaces embedded (with their
closures) in $S^\circ$ such that
$\bm{A}$ is contined in $\ol{S}_+$ as a proper
system of arcs. More concretely, $\ol{S}_+$ is the union
of a closed tubular neighbourhood of $\bm{A}$ and all components of $S\setminus\bm{A}$
which are discs or annuli that retract onto some boundary component.
If $\Si$ is obtained from
$S$ by collapsing the boundary components of $\ol{S}_+$
and the possible
resulting two-noded spheres to nodes of $\Si$,
then we obtain an isotopy class of maps $f:S\rar\Si$, which depends
only on $\bm{A}$. We will refer to this map (or just to $\Si$,
when we work in the moduli space) as the {\it topological type} of $\bm{A}$.

Given weights $w\in |\bm{A}|^\circ\times[0,\infty)$,
the components of $\ol{\Si}_+=f(\ol{S}_+)$ are quasi-filled by the arc system $f(\bm{A})$:
because of Theorem~\ref{prop:smooth}, they can be given
a hyperbolic metric such that $f(\bm{A})$ is its
spinal arc system with weights $f_*(w)$.
When no confusion is possible, we will still denote by
$[f:S\rar\Si]$ the class of visibly equivalent
$S$-marked stable surfaces determined by $f$.

This construction defines a $\G(S)$-equivariant extension
of the previous $\bm{W}^{-1}$
\[
\xymatrix@R=0in{
\wh{\bm{W}}^{-1}\ :\ |\Af(S)|\times[0,\infty) \ar[rr] &&
\That^{vis}(S)
}
\]
where $\That^{vis}(S)=\That(S)/\!\!\sim_{vis}$.

The argument above shows that $\wh{\bm{W}}^{-1}$ is bijective.
As already noticed in \cite{bowditch-epstein:natural}
and \cite{looijenga:cellular},
the map $\wh{\bm{W}}$ is not continuous if $|\Af(S)|$
is endowed with the coherent topology.

\begin{remark}
$|\Af(S)|$ is locally finite at $w$
$\iff$ $\bm{A}=\mathrm{supp}(w)$ is a proper system of arcs
$\iff$ $w$ has a countable fundamental system of
coherent neighbourhoods.
Moreover,
a sequence converges (for the coherent topology)
if and only if
it is definitely in a fixed closed simplex and there
it converges in the Euclidean topology.
\end{remark}

The discontinuity of $\wh{\bm{W}}$ at $\pa\Teich(S)$
with respect to the coherent topology
can be seen as follows. Consider a marked surface
$[f:S\rar\Si]$
with a node $f(\g)=q\in\Si$ such that not all
the boundary components of $\Si$ are cusps and
call $\bm{A}$ a maximal system of arcs of $S$ such that
$\wh{\bm{W}}(f)\in|\bm{A}|\times\R_+$.
Choose a sequence $[f_m:S\rar\Si_m]$ with $\wh{\bm{W}}(f_m)$
contained in $|\bm{A}|^\circ\times\R_+$ and such that $[f_m]\rar[f]$.
If $\tau_{\g}$ is the right Dehn twist along $\g$
and $f'_m=f_m\circ\tau_{\g}^m$, then $[f'_m]$ still converges to
$[f]$. On the other hand, the $\wh{\bm{W}}(f'_m)$'s all
belong to the interior of distinct maximal simplices of $|\Af(S)|$
and so the sequence $\wh{\bm{W}}(f'_m)$ is divergent
for the coherent topology.

The correct solution (see \cite{bowditch-epstein:natural}),
anticipated in Section~\ref{sec:laminations} and
which we will adopt without further notice,
is to equip $|\Af(S)|$
with the metric topology, whose importance
will be also clear in the proof of Lemma~\ref{lemma:isom}.

\begin{theorem}\label{prop:w}
The $\G(S)$-equivariant natural extension
\[
\xymatrix@R=0in{
\wh{\bm{W}}\ :\ \That^{vis}(S) \ar[rr] &&
|\Af(S)|\times[0,\infty)
}
\]
is a homeomorphism.
\end{theorem}
%
%
%
%
%
The following proof shares some ideas with
\cite{acgh2} (to which we refer for a more detailed discussion
of the case with $n$ cusps).
The bijectivity of $\wh{\bm{W}}$ is a direct consequences
of the work of Penner/Bowditch-Epstein and Luo.
We begin with some preparatory observations.

\begin{definition}
Let $([f:S\rar\Si],\up)\in\That(S)(0)$ be a projectively
decorated surface and let $B_i\subset\Si$ be the embedded
horoball at $x_i=f(C_i)$
with radius $p_i$.  The associated
{\it truncated surface} is
$\Si^{tr}:=\Si\setminus(B_1\cup\dots\cup B_n)$ and
the {\it reduced length} of an arc $\a\in\Aa(S)$ at $f$ is
$\tilde{\ell}_{\a}(f):=\ell(\Si^{tr}\cap f(\a))$.
\end{definition}

\begin{lemma}\label{lemma:uniform}
Let $\{f_m:S\rar\Si_m\}\subset\Teich(S)$ be a sequence that
converges to $[f_\infty:S\rar\Si_\infty]\in\That^{vis}(S)$.
\begin{enumerate}
\item[(a)]
Assume $\bm{\Ll}(f_\infty)>0$ and let $\Aa(S)=\Aa_{fin}\sqcup
\Aa_\infty$, where $\Aa_\infty$ is the subset of arcs
$\b$ such that $\ell_\b(f_\infty)=\infty$.
Then $\ell_\a(f_m)/\ell_\a(f_\infty)\rar 1$ uniformly
for all $\a\in\Aa_{fin}$. Moreover, if $\Aa_\infty\neq\emptyset$,
then there exists a diverging sequence $\{L_m\}\subset\R_+$
such that $\ell_\b(f_m)\geq L_m$ for all $\b\in\Aa_\infty$.
Hence, $t_\bullet(f_m)\rar t_\bullet(f_\infty)$ uniformly.
\item[(b)]
Assume $\Ll(f_\infty)=(\ti{\up},0)\in\Delta^{n-1}\times\{0\}$ and
let $\Aa(S)=\Aa_{fin}\sqcup \Aa_\infty$,
where $\Aa_\infty$ is the subset of arcs
$\b$ such that $\ti{\ell}_\b(f_\infty)=\infty$.
Then $\ti{\ell}_\a(f_m)/\ti{\ell}_\a(f_\infty)\rar 1$ uniformly
for all $\a\in\Aa_{fin}$. Moreover, if $\Aa_\infty\neq\emptyset$,
then there exists a diverging sequence $\{\ti{L}_m\}\subset\R_+$
such that $\ti{\ell}_\b(f_m)\geq \ti{L}_m$ for all $\b\in\Aa_\infty$.
Hence, $[t_\bullet(f_m)]\rar [t_\bullet(f_\infty)]$.
\end{enumerate}
\end{lemma}

\begin{remark}\label{rem:normalized}
A simple computation shows that a hypercycle at distance $d$
from a closed geodesic of length $\ell$ has length $\ell\cosh(d)$.
In case (b), we can assume that $\bm{\Ll}(f_m)\leq 1$ and
so we can define $\th_m\in[0,\pi/2]$ by $\sin(\th_m):=\bm{\Ll}(f_m)$.
For each boundary circle $C_{i,m}$ of $\Si_m$, let $B_{i,m}
\subset\Si_m$ be the hypercycle parallel to $C_{i,m}$ at
distance $d_m=-\log\tan(\th_m/2)$ (i.e. $\cosh(d_m)=1/\sin(\th_m)$),
which has length $\ti{p}_i(f_m):=p_i(f_m)/\sin(\th_m)\leq 1$ and so is
embedded in $\Si_m$.
Notice that the spine of $\Si_m$ coincides with the spine of its
subsurface $\Si_m^{tr}$ obtained by removing the
hyperballs bounded by the $B_{i,m}$'s: in fact, every geodesic
that meets $C_{i,m}$ orthogonally also intersects $B_{i,m}$
orthogonally. For every arc $\a$, define the {\it reduced length}
of $\a$ at $f_m$
to be $\ti{\ell}_\a(f_m):=\ell_\a(f_m)-2d_m$, namely the length
of $f_m(\a)\cap\Si_m^{tr}$. Extending a definition of Penner's
\cite{penner:decorated} (and modifying it by a factor $\sqrt{2}$),
we put $\l_\a(f_m):=\exp(\ti{\ell}_\a/2)$. Because $B_{i,m}$ limits
to a horoball of circumference $\ti{p}_i=\ti{p}_i(f_\infty)$
as $m\rar\infty$, the length
$\l_\a(f_m)\rar\l_\a(f_\infty,\up)$.
\end{remark}

\begin{notation}
In the following proof, we will denote by $S_{\infty,+}$
the open subsurface $f_\infty^{-1}(\Si_{\infty,+})$
and by $S_{\infty,+}^{tr}$
the preimage through $f_\infty$ of the analogous truncated
surface. Similar notation for $S_{\infty,-}$.
\end{notation}

\begin{proof}[Proof of Lemma~\ref{lemma:uniform}]
About (a), if $\Aa_\infty\neq\emptyset$, then
$\exists \g_1,\dots,\g_l\in\mathcal{C}(S)$ disjoint
such that $c_m=\mathrm{max}_h \ell_{\g_h}(f_m)\rar 0$
as $m\rar\infty$. Clearly, $\b\in\Aa_\infty
\iff i(\b,\g_1+\dots+\g_l)>0$.

By the collar lemma, $\ell_\b(f_m)>L_m:=T(c_m)/2$ and so
$t_\b(f_m)<T(L_m)\rar 0$ for all $\b\in\Aa_\infty$.
On the other hand,
by Proposition~\ref{prop:convergence},
we can assume that $f_m^*(g_m)\rar f_\infty^*(g_\infty)$
in $L^\infty_{loc}(S_{\infty,+})$.
Thus,
$\dis\frac{|\ell_{\a}(f_m)-\ell_{\a}(f_\infty)|}
{\ell_{\a}(f_\infty)}\rar 0$ uniformly for all $\a\in \Aa_{fin}$.

Fix $\e>0$ and let $\a_1,\dots,\a_k\in\Aa_{fin}$ be the arcs such that
$\ell_{\a_i}(f_\infty)\leq T(\e)/(1-\e)$ for $i=1,\dots,k$.
Clearly, $t_{\a_i}(f_m)\rar t_{\a_i}(f_\infty)$.
If $\a\in \Aa_{fin}$ and $\a\notin\{\a_1,\dots,\a_k\}$, then
$\ell_\a(f_m)\geq \ell_\a(f_\infty)(1-\e)>T(\e)$
and so $t_\a(f_m)<\e$ for $m$ large.
Hence, $|t_\a(f_m)-t_\a(f_\infty)|<\e$ for $m$ large and
$\a\in\Aa_{fin}\setminus\{\a_1,\dots,\a_k\}$.\\

The proof of (b) is similar.
Call $\g_1,\dots,\g_l$ the curves in the interior of $S$
that are shrunk to nodes of $\Si_\infty$
and let $J=\{j\,|\, \ti{p}_j=0\}$.
We can assume that $p_i(f_m)<\ti{p}_i(f_\infty)$.

Let $c_m=\mathrm{max}\{\ell_{\g_h}(f_m) \}$ and
$c'_m=\mathrm{max}\{p_j\,|\,j\in J \}$.
Clearly, if $\b\in\Aa_\infty$ intersects some $\g_h$, then 
$\ti{\ell}_\b(f_m)\geq T(c_m)/2\rar\infty$.
If $\b\in\Aa_\infty$ does not intersect any $\g_j$, then it starts
at some $C_j$ with $j\in J$.
Because of the collar lemma, there is a hypercycle
embedded in $\Si_m$
at distance $\d_{i,m}$ from $f_m(C_i)$,
with $p_i(f_m) \cosh(\delta_{i,m})=1$.
As $p_i(f_m)\cosh(d_m)=p_i(f_m)/\sin(\th_m)$, we get
$\cosh(\delta_{i,m})/\cosh(d_m)=\sin(\th_m)/p_i(f_m)$
and so $\delta_{j,m}-d_m\approx \log(\sin(\th_m)/p_j(f_m))
\geq \log(\sin(\th_m)/c'_m)\rar\infty$ for $j\in J$.
Hence, $\ti{\ell}_\b(f_m)\geq \ti{L}_m:=\mathrm{min}\{ T(c_m)/2,\log(\sin(\th_m)/c'_m)\}
\rar\infty$.

%
As before, we can assume that $f_m^*(g_m)$ converges uniformly over
the compact subsets of $S_{\infty,+}^{tr}$, so that
$\dis\frac{|\ti{\ell}_{\a}(f_m)-\ti{\ell}_{\a}(f_\infty)|}
{\ti{\ell}_{\a}(f_\infty)}\rar 0$ uniformly for all $\a\in \Aa_{fin}$.

Call $\a_0\in\Aa_{fin}$ the arc with smallest $\ti{\ell}_{\a_0}(f_\infty)$.
Fix $\e>0$ and let $\a_1,\dots,\a_k\in\Aa_{fin}$ be the arcs such that
$\ti{\ell}_{\a_i}(f_\infty)\leq \ti{\ell}_{\a_0}(f_\infty)-2\log(\e)$
for $i=1,\dots,k$.
Clearly, $\dis \frac{t_{\a_i}(f_m)}{t_{\a_0}(f_m)}
\rar \frac{t_{\a_i}(f_\infty)}{t_{\a_0}(f_\infty)}=
\frac{\l_{\a_0}(f_\infty)}{\l_{\a_i}(f_\infty)}$
for $i=1,\dots,k$.

If $\a\in \Aa_{fin}$ and $\a\notin\{\a_1,\dots,\a_k\}$, then
$\dis\frac{t_{\a_i}(f_m)}{t_{\a_0}(f_m)}\leq
\e+\sqrt{\exp[\ti{\ell}_{\a_0}(f_m)-\ti{\ell}_{\a_i}(f_m)]}<2\e$
for $m$ large.
Hence, $\dis\left|
\frac{t_{\a_i}(f_m)}{t_{\a_0}(f_\infty)}
\right|<2\e$ for $m$ large and
$\a\in\Aa_{fin}\setminus\{\a_1,\dots,\a_k\}$.
\end{proof}
\begin{proof}[Proof of Theorem~\ref{prop:w}]
The continuity of $\wh{\bm{W}}$ is dealt with in
Lemma~\ref{lemma:W-continuous} below.
In order to prove that $\wh{\bm{W}}$ is a homeomorphism,
it is sufficient to show that
the induced map below is.
\[
\xymatrix@R=0in{
\wh{\bm{W}}':
\Mhat^{vis}(S) \ar[rr] &&
\left(|\Af(S)|/\G(S)\right)\times[0,\infty)
}
\]
where $\Mhat^{vis}(S)=\That^{vis}(S)/\G(S)$.

In fact, we first endow $\Tbar^{WP}(S)$ with
a $\G(S)$-equivariant metric, for example pulling it
back from $\Tbar^{WP}(S)\rar\Mbar(S)$.
This way, we induce a metric on the quotient
$\Tbar^{WP}(S)\times\Delta^{n-1}/\!\!\sim_{vis}$,
where $\Delta^{n-1}$ has the Euclidean metric.
Finally, we embed
$\That^{vis}(S)$ inside $\Tbar^{WP}(S)\times
\Delta^{n-1}/\!\!\sim_{vis}$
(where the second component of the map is given by
the normalized boundary lengths), thus obtaining a
$\G(S)$-equivariant metric on $\That^{vis}(S)$.

Then, $\That^{vis}(S)$ and $|\Af(S)|$
are metric spaces
and $\G(S)$ acts on both by isometries.
Moreover, the action on $|\Af(S)|$ is simplicial on the
second baricentric subdivision, and so
its orbits are discrete.

On the other hand, the map $\wh{\bm{W}}'$ is clearly proper,
because $\That^{vis}(S)(\up)/\G(S)$ is compact for
every $\up\in\Delta^{n-1}\times[0,\infty)$.
As the image of $\wh{\bm{W}}'$ contains the dense
open subset $\dis\left(|\Ao(S)|/\G(S)\right)\times(0,\infty)$,
we have that $\wh{\bm{W}}'$ is a homeomorphism.
By Lemma~\ref{lemma:isom}(b), $\wh{\bm{W}}$ is a homeomorphism too.
\end{proof}
\begin{lemma}\label{lemma:W-continuous}
$\widehat{\bm{W}}$ is continuous.
\end{lemma}
\begin{proof}
Notice that $\That^{vis}(S)$ and $|\Af(S)|\times[0,\infty)$
have countable systems of neighbourhoods at each point.
As $\Teich(S)$ is dense in $\That^{vis}(S)$,
in order to test the continuity of $\wh{\bm{W}}$, we can consider 
a sequence $[f_m:S\rar\Si_m]\subset\Teich(S)$ converging to
a point $[f_\infty:S\rar\Si_\infty]\in\That^{vis}(S)
\setminus\That^{vis}(S)(0)$
(the case of $[f_\infty]\in\That^{vis}(S)(0)$ will be treated later).

{\bf Step 1.}
Because of Proposition~\ref{prop:convergence},
there are representatives $f_1,f_2,\dots,f_{\infty}$ such that
the hyperbolic metric $f_m^*(g_m)\rar f_{\infty}^*(g_{\infty})$
in $L^{\infty}_{loc}(S_{\infty,+})$.
Also, the distance
function $d_{f_m}(-,\pa S_{\infty,+}): S_{\infty,+}\rar \R_+$ 
with respect to the metric $f_m^* g_m$ converges in
$L^{\infty}_{loc}(S_{\infty,+})$.

{\bf Step 2.}
Let $\mathcal{E}$ be the set of edges of $\mathrm{Sp}(\Si_{\infty})$
and let $m_i$ be the midpoint of the edge $E_i\in \mathcal{E}$
in $\Si_{\infty}$.
Call $\g_{i,1}$ and $\g_{i,2}$ the shortest geodesics that join
$m_i$ to $\pa\Si_{\infty}$ and $\a_i:=f_\infty^{-1}(\g_{i,1}^{-1}\ast\g_{i,2})$
the associated arc.
Let $d(m_i,\pa\Si_{\infty})=\ell(\g_{E,1})=
\ell(\g_{E,2})$ and call
$d'(m_i,\pa\Si_{\infty})$ to be the minimum length of a geodesics
that join $m_i$ to $\pa\Si_{\infty}$ and is not homotopic to $\g_{i,1}$
or $\g_{i,2}$.
Finally, set $\e=\mathrm{min}\{d'(m_i,\pa\Si_{\infty})-
d(m_i,\pa\Si_{\infty})\,|\,E_i\in \mathcal{E}\}>0$.

Because $d_{f_m}(f_\infty^{-1}(m_i),\pa S_{\infty,+})$ and
$d'_{f_m}(f_\infty^{-1}(m_i),\pa S_{\infty,+})$ also converge
as $m\rar\infty$,
their difference is eventually positive
and so the arc $\a_i$ (up to isotopy) is still dual  
to some edge of the spines
of $(S,f_m^*(g_m))$ (which is equal to
$f_m^{-1}(\mathrm{Sp}(\Si_m))$) for $m$ large.

%
Thus, up to discarding finitely many terms
of the sequence,
we can assume that
$f_{\infty}^*\bm{A}_{sp}(\Si_{\infty})\subseteq f_m^*\bm{A}_{sp}(\Si_m)$.

{\bf Step 3.}
Let $\bm{A}_{\infty}$
be the system of arcs $f_{\infty}^*\bm{A}_{sp}(\Si_{\infty})$ on $S$.
Consider the subset $\widetilde{\mathfrak{St}}(\bm{A}_{\infty})\subset\Af(S)$
of systems $\bm{A}$ that contain $\bm{A}_{\infty}$ and such that
$f(\bm{A})$ can be represented inside $\Si_{\infty,+}$.
Let $\bm{A}_1,\dots,\bm{A}_k$ the maximal elements of
$\widetilde{\mathfrak{St}}(\bm{A}_{\infty})$
and let
$\widetilde{\mathfrak{St}}_i=\{\bm{A}_{i,r}\,|\,r\in R_i\}$
the set of maximal systems
of arcs $\bm{A}_{i,r}\supseteq\bm{A}_i$ for $i=1,\dots,k$.

{\bf Step 4.}
Clearly $\exists i_m\in\{1,\dots,k\}$ and $\exists r_m\in R_{i_m}$
such that $f_m^*\bm{A}_{sp}(\Si_m)\subseteq\bm{A}_{i_m,r_m}$
(and there are finitely many options for each $m$).
We need to show that
\[
\max
\{
|w_{\bm{A}_{i_m,r_m}}(\ora{\a},f_m)-w_{\bm{A}_{i_m,r_m}}(\ora{\a},f_\infty)|:
\a\in\bm{A}_{i_m,r_m}\supseteq f_m^*\bm{A}_{sp}(\Sigma_m)
\}
\rar 0
\]
as $m\rar\infty$.
%

{\bf Step 5.}
By Lemma~\ref{lemma:uniform}(a),
$\ell_\b(f_m)\geq L_m$ equidiverge and
$w_{\bm{A}_{i_m,r_m}}(\ora{\b},f_m)\leq t_{\b}(f_m)/2$
uniformly converge to zero, for all $\b\in\Aa_\infty$.

{\bf Step 6.}
$\bm{A}_1\cup\dots\cup\bm{A}_k$ is finite and the lengths
$\ell_\a(f_m)\rar\ell_\a(f_{\infty})<\infty$ for every
$\a$ in some $\bm{A}_i$.
Define $M(\a)=\{m\,|\,\a\in f_m^*\bm{A}_{sp}(\Si_m)\}$.
It is sufficient to prove that, if $M(\a)$ is infinite,
then $|w_{\bm{A}_{i_m,r_m}}(\ora{\a},f_m)-
w_{\bm{A}_{i_m,r_m}}(\ora{\a},f_\infty)|\rar 0$
as $m\in M(\a)$ diverges.
There are three cases.

{\bf Case 6(a).}
Let $H_m\subset S$ be a right-angled hexagon bounded
by $(\ora{\a},\ora{\a'_m},\ora{\a''_m})$, with
$\a'_m,\a''_m\in\bm{A}_{i_m}$ for suitable $i_m$.
Then
\[
\sinh
\left(
w_{\bm{A}_{i_m,r_m}}(\ora{\a},f_m)
\right)
=\frac{s_{\a'_m}(f_m)^2+
s_{\a''_m}(f_m)^2-s_\a(f_m)^2}{
\dis 2s_{\a'_m}(f_m)s_{\a''_m}(f_m)\sqrt{s_\a(f_m)^2-1}}
\]
and so
$|w_{\bm{A}_{i_m,r_m}}(\ora{\a},f_m)-w_{\bm{A}_{i_m,r_m}}(\ora{\a},f_{\infty})|\rar 0$.

{\bf Case 6(b).}
Suppose that there are hexagons $H_m\subset S$
with $\pa H_m=\ora{\a}\ast\ora{\a'_m}\ast\ora{\b_m}$,
where $\a'_m\in\bm{A}_{i_m}$ and $\b_m\in\bm{A}_{i_m,r_m}\setminus\bm{A}_{i_m}$.
We can extract a subsequence so that $\a'_m$ is a fixed
arc $\a'$. The divergence of $b_m$
and the formula
\[
\cosh(d(\ora{\a},\ola{\a'}))=
\frac{\cosh(a_m)\cosh(a'_m)+\cosh(b_m)}{\sinh(a_m)\sinh(a'_m)}
\]
(where $a_m,a'_m,b_m$ are the lengths of $\a,\a',\b_m$ at $[f_m]$)
imply that $d(\ora{\a},\ola{\a'})$
diverges, which contradicts the fact that the boundary lengths
are bounded.

{\bf Case 6(c).}
Let $H_m\subset S$ be a right-angled hexagon bounded
by $(\ora{\a},\ora{\b_m},\ora{\b'_m})$, with
$\b_m,\b'_m \in \bm{A}_{i_m,r_m}\setminus\bm{A}_{i_m}$.
Call $x_{\a,m},x_{\b,m},
x_{\b',m}$ the length
of the edges of $H_m$ opposed to $\a,\b_m,\b'_m$ and let $a_m,b_m,b'_m$
be the lengths of $\a,\b_m,\b'_m$ at $[f_m]$.
Remember that $w_{\bm{A}_{i_m,r_m}}(\ora{\b_m},f_m)$ and
$w_{\bm{A}_{i_m,r_m}}(\ora{\b'_m},f_m)$ converge to zero,
whereas $w_{\bm{A}_{i_m,r_m}}(\ora{\a},f_m)$ is bounded
(and so are $x_{\b,m}$ and $x_{\b',m}$).
%
Notice that $x_{\b,m}-w_{\bm{A}_{i_m,r_m}}(\ora{\a},f_m)$ converges
to zero
and so do the differences
$\cosh(x_{\b,m})-\cosh(w_{\bm{A}_{i_m,r_m}}(\ora{\a},f_m))$
and $\cosh(x_{\b,m})-\cosh(x_{\b',m})$.
%
But
\[
\cosh(x_{\b,m})=
\frac{\cosh(b'_m)\cosh(a_m)+\cosh(b_m)}{\sinh(b'_m)\sinh(a_m)}
=\frac{1}{\tanh(a_m)\tanh(b'_m)}+\frac{\cosh(b_m)}{\sinh(a_m)\sinh(b'_m)}
\]
and similarly for $x_{\b'_m}$, so that we obtain that
\[
\cosh(x_{\b,m})-\cosh(x_{\b',m})\approx
\frac{e^{b_m-b'_m}-e^{b'_m-b_m}}{\sinh(a_m)}
=\frac{2\sinh(b_m-b'_m)}{\sinh(a_m)}\rar 0
\]
which implies that $|b_m-b'_m|\rar 0$, because
$a_m\rar a_{\infty}\in(0,\infty)$.

Consequently, $\dis \cosh(w_{\bm{A}_{i_m,r_m}}(\ora{\a},f_m))
\rar\frac{1}{\tanh(a_\infty)}+\frac{1}{\sinh(a_\infty)}=
\frac{1}{\tanh(a_\infty/2)}$,
which gives $w_{\bm{A}_{i_m,r_m}}(\ora{\a},f_m)\rar t_{\a}(f_{\infty})/2
=w_{\bm{A}_{\infty}}(\ora{\a},f_{\infty})$.\\

{\bf\it Case of decorated surfaces.}

Suppose now that $[f_\infty,\up]\in\That^{vis}(S)(0)$.
We use the notation in Remark~\ref{rem:normalized}.
Notice that
\begin{align}\tag{*}
\l_\a(f_m) & =e^{\ell_\a/2-d_m}=e^{-d_m}
\left(s_\a(f_m)+\sqrt{s_\a(f_m)^2-1}\right)=\\
\notag
&=\tan(\vartheta_m/2)\left(s_\a(f_m)+\sqrt{s_\a(f_m)^2-1}\right)
\end{align}
The normalized widths
$\ti{w}_{\bm{A}_{i_m,r_m}}=2w_{\bm{A}_{i_m,r_m}}/\sin(\vartheta_m)$
limit to Penner's simplicial coordinates
(see below the modifications to step (6)).
So the map $\wh{\bm{W}}$ reduces to Penner's map for cusped surfaces,
in which case
we will still use the term ``normalized widths''
(instead of ``simplicial coordinates'') for brevity.

We follow the same path as before, with some modifications.

{\bf Step 1.}
As $\ti{\bm{p}}=1$, we can assume that
$p_i(f_m)<\ti{p}_i(f_m)$. Let $\Si_\infty^{tr}$ be
the truncated surface as in the proof of Lemma~\ref{lemma:uniform}(b).
By Proposition~\ref{prop:convergence}, we can assume
that $f_m(S^{tr})=\Si_m^{tr}$ and that
the metrics $f_m^*(g_m)$ converge
in $L^{\infty}_{loc}(S_{\infty,+}^{tr})$.

{\bf Step 2.}
Essentially the same, up to replacing
the distance from $\pa\Si_m$ by the distance from
$\pa\Si_m^{tr}$. 

{\bf Step 3.} Identical.

{\bf Step 4.} Now on, we have to replace
the widths by the normalized widths $\ti{w}$.
Notice that $\ti{w}_{\bm{A}_{i_m,r_m}}(\ora{\a},f_m)\leq
t_\a(f_m)/\sin(\th_m)\approx 2\mathrm{exp}(-\ti{\ell}_\a/2)$
for all $\a\in\bm{A}_{i_m,r_m}$.

{\bf Step 5.}
Similar: by Lemma~\ref{lemma:uniform}(b),
$\ti{\ell}_\b(f_m)\geq \tilde{L}_m$ equidiverge and
$\ti{w}_{\bm{A}_{i_m,r_m}}(\ora{\b},f_m)\rar 0$ uniformly,
for all $\b\in\Aa_\infty$.

{\bf Step 6.}
It follows from (1) that, as $m\rar \infty$,
for all $\a\in\bm{A}_1\cup\dots\cup\bm{A}_k$
we have $\l_\a(f_{\infty})<\infty$ and
$|\l_\a(f_m)-\l_\a(f_{\infty})|\rar 0$.

It follows from $(*)$ that $\l_\a(f_m)\sim
\vartheta_m s_\a(f_m)+O(\vartheta_m^3 s_\a(f_m)
)$.
Hence, as $m\rar\infty$,
for all these $\a$, we also have
$|\l_\a(f_m)-\vartheta_m s_\a(f_m)|\rar 0$.
On the other hand, for all $\b$ belonging
to some $\bm{A}_{i_m,r_m}\setminus\bm{A}_{i_m}$,
we have $\l_\b(f_\infty)=\infty$ and
$\l_\b(f_m)\sim \vartheta_m s_\b(f_m)$
equidiverge.

{\bf Case 6(a).}
$|\ti{w}_{\bm{A}_{i_m,r_m}}(\ora{\a},f_m)-
\ti{w}_{\bm{A}_{i_m,r_m}}(\ora{\a},f_\infty)|\rar 0$
and $\ti{w}_{\bm{A}_{i_m,r_m}}(\ora{\a},f_m)
\rar X_{\bm{A}_\infty}(\ora{\a},f_m)$,
where $\dis
X_{\bm{A}_{\infty}}(\ora{\a},f_\infty)=
\frac{\l_{\a_i}(f_\infty)^2+\l_{\a_j}(f_{\infty})^2
-\l_\a(f_\infty)^2}{\l_{\a_i}(f_\infty)\l_{\a_j}(f_\infty)
\l_\a(f_\infty)}$
is Penner's simplicial coordinate of $\ora{\a}$.

{\bf Case 6(b).}
$\b_m$ cannot cross a simple closed (nonboundary) curve of $S$
that is contracted to a node by $f_{\infty}$,
because so it would either $\a$ or $\a'$: this would
contradict the boundedness of $\ti{\ell}_\a(f_m)$ and
$\ti{\ell}_{\a'}(f_m)$.

{\bf Case 6(c).}
Because $\cosh(x_{\a,m})\approx 1+x_{\a,m}^2/2$ and
\begin{align*}
\cosh(x_{\a,m})=\frac{\cosh(b_m)\cosh(b'_m)+\cosh(a_m)}
{\sinh(b_m)\sinh(b'_m)}& \approx 1+2\exp(a_m-b_m-b'_m)+\\
& + 2\exp(-2b_m)+2\exp(-2b'_m)
\end{align*}
we obtain that $x_{\a,m}^2/\vartheta_m^2\approx
\exp(\ti{a}_m-\ti{b}_m-\ti{b}'_m)+O(\th_m^2)\rar 0$ as $m\rar \infty$.

Remember that $\ti{w}_{\bm{A}_{i_m,r_m}}(\ora{\b}_m,f_m)$
and $\ti{w}_{\bm{A}_{i_m,r_m}}(\ora{\b}'_m,f_m)$ converge to
zero uniformly and that $\ti{w}_{\bm{A}_{i_m,r_m}}(\ora{\a},f_m)$
is bounded (and so are $x_{\b,m}/\sin(\th_m)$ and $x_{\b',m}/\sin(\th_m)$).

On the other hand, 
$
\cosh(x_{\b',m})\approx 1+2\exp(b'_m-b_m-a_m)+
2\exp(-2a_m)+2\exp(-2b_m)
$
and so $x_{\b',m}^2\approx 4 \exp(b'_m-b_m-a_m)+
4\exp(-2b_m)+4\exp(-2a_m)$
and $\dis \frac{x_{\b',m}^2}{\sin^2(\th_m)}\approx
\exp(\ti{b}'_m-\ti{b}_m-\ti{a}_m)+O(\th_m^2)$.

On the other hand, $\dis \ti{w}_{\bm{A}_{i_m,r_m}}(\ora{\a},f_m)
\approx \frac{2 x_{\b,m}}{\sin(\th_m)}\approx
2\exp\left(\frac{\ti{b'}_m-\ti{b}_m-\ti{a}_m}{2}\right)+
O(\th_m)$ and an analogous estimate
holds switching the roles of $\ti{b}_m$ and $\ti{b}'_m$.
This implies that $|\ti{b}'_m-\ti{b}_m|\rar 0$.

As a consequence, $\ti{w}_{\bm{A}_{i_m,r_m}}(\ora{\a},f_m)
\approx 2\exp(-\ti{a}_m/2)\rar 2\exp(-\ti{a}_\infty/2)=
2/\l_\a(f_\infty)=X_{\bm{A}_\infty}(\ora{\a},f_\infty)$.
\end{proof}
\begin{lemma}\label{lemma:isom}
Let $f:X\rar Y$ be a $G$-equivariant
continuous map between metric spaces
on which the discrete group $G$ acts
by isometries. Assume that the $G$-orbits
on $Y$ are discrete.
\begin{itemize}
\item[(a)]
If $X/G$ is compact and $\mathrm{stab}(x)\subseteq
\mathrm{stab}(f(x))$ has finite index for all $x\in X$,
then $f$ is proper.
\item[(b)]
If $f$ is injective and
the induced map
$\ol{f}:X/G\rar Y/G$ is a homeomorphism,
then $f$ is a homeomorphism.
\end{itemize}
\end{lemma}
\begin{proof}
For (a) we argue by contradiction: let $\{x_n\}\subset X$
be a diverging subsequence such that $\{f(x_n)\}\subset Y$
is not diverging. Up to extracting a subsequence, we can
assume that $f(x_n)\rar y\in Y$ and that $[x_n]\rar [x]\in X/G$.
Thus $\exists g_n\in G$ such that
$x_n\cdot g_n\rar x$, that is $d_X(x_n,x\cdot g_n^{-1})\rar 0$.
As $\{x_n\}$ is divergent, the sequence $\{[g_n]\}\subset G/\mathrm{stab}(x)$
is divergent too, and so is it in $G/\mathrm{stab}(f(x))$.
Because $f(x_n\cdot g_n)\rar f(x)$, we have
$d_Y(f(x_n),f(x)\cdot g_n^{-1})\rar 0$ and so $\{f(x_n)\}$ is divergent,
because $f(x)\cdot G$ is discrete.

For (b), let's show first that $f$ is surjective.
Because $\ol{f}$ is bijective,
for every $y\in Y$ there exists a unique $[x]\in X/G$
such that $\ol{f}([x])=[y]$. Hence, $f(x)=y\cdot g$ for
some $g\in G$ and so $f(x\cdot g^{-1})=y$.

The injectivity of $f$ also implies that
$\mathrm{stab}(x)=\mathrm{stab}(f(x))$ for all $x\in X$.

To prove that $f^{-1}$ is continuous,
let $(x_m)\subset X$ be a sequence such that $f(x_m)\rar f(x)$
as $m\rar \infty$ for some $x\in X$.
Clearly, $[f(x_m)]\rar[f(x)]$ in $Y/G$ and so $[x_m]\rar [x]$
in $X/G$, because $\ol{f}$ is a homeomorphism.

Consider the balls $U_k=B_X(x,1/k)$ for $k>0$ and set $U_0=X$,
so that $[U_k]$ is an open neighbourhood of $[x]\in X/G$.
There exists an increasing sequence
$\{m_k\}$ such that $[x_m]\in [U_k]$,
that is $\dis x_m \in \bigcup_{g\in G} U_k \cdot g$
for all $m\geq m_k$. Consequently, there is a $g_m\in G$
such that $x_m\in U_k\cdot g_m$ for every $m_k\leq m< m_{k+1}$.
Thus, $z_m:=x_m \cdot g_m^{-1}\rar x$.
By continuity of $f$, we have $f(z_m)\lra f(x)$
and by hypothesis $f(z_m)\cdot g_m \rar f(x)$.
Moreover, $d_Y(f(x)\cdot g_m,f(x))\leq d_Y(f(x)\cdot g_m,f(x_m))+
d_Y(f(x_m),f(x))
=d_Y(f(x),f(z_m))+d_Y(f(x_m),f(x))\rar 0$ and so
$f(x)\cdot g_m\rar f(x)$. Hence,
$g_m\in\mathrm{stab}(f(x))=\mathrm{stab}(x)$ for large $m$,
because $G$ acts with discrete orbits on $Y$.
As a consequence, for $m$ large enough
$d_X(x_m,x)=d_X(z_m,x)\rar 0$ and so
$x_m\rar x$ and $f^{-1}$ is continuous at $f(x)$.
\end{proof}

\begin{remark}
In order to check that the $G$-orbits on $Y$ are discrete,
it is sufficient to show the following:
\begin{itemize}
\item[($*$)]
whenever $y\cdot g_m\rar y$
for a certain $y\in Y$ and $\{g_m\}\subset G$, the sequence
$\{g_m\}$ definitely belongs to $\mathrm{stab}_G(y)$.
\end{itemize}
Assuming ($*$), 
there is an $\e>0$ and
a ball $B=B(z,\e)$ such that $B\cap z\cdot G=\{z\}$.
Given a sequence $\{g_m\}\subset G$ such that
$y\cdot g_m\rar z\in Y$,
then $d(z\cdot g_j^{-1} g_i,z)\leq d(z\cdot g_j^{-1},y)+
d(y,z\cdot g_i^{-1})=d(z,y\cdot g_j)+d(y\cdot g_i,z)<\e$
for $i,j\geq N_\e$.
Thus, $g_j^{-1}g_i\in \mathrm{stab}_G(z)$
and $d(y\cdot g_j,z)=d(y\cdot g_j g_j^{-1} g_i,z)=
d(y\cdot g_i,z)$ for all $i,j\geq N_\e$.
Hence, $y\cdot g_i=z$ for all $i\geq N_\e$ and so the orbit is discrete.
\end{remark}

%
%
\end{subsection}
%
%
\begin{subsection}{The bordification of arcs}
Define a map
\[
\Phi:|\Af(S)|\times[0,\infty]\lra\Tbar^a(S)
\]
in the following way:
\[
\Phi(w,\bm{p})=\begin{cases}
([\l^{-1}_\bullet(\wh{\bm{W}}^{-1}(w,0))],0) & \text{if $\bm{p}=0$}\\
j(\wh{\bm{W}}^{-1}(w,\bm{p})) & \text{if $0<\bm{p}<\infty$} \\
([w],\infty) & \text{if $\bm{p}=\infty$.}
\end{cases}
\]

The situation is thus as in the following diagram.
\[
\xymatrix{
|\Af(S)|\times[0,\infty)
\ar@{<-}[rr]^{\qquad\wh{\bm{W}}}_{\qquad\cong}
\ar@{^(->}[d] &&
\That^{vis}(S) \ar@{^(->}[d]_{\hat{j}} \\
|\Af(S)|\times[0,\infty] \ar[rr]^{\qquad\Phi} &&
\Tbar^a(S)
}
\]
\begin{theorem}\label{thm:phi}
$\Phi$ is a $\G(S)$-equivariant
homeomorphism. Thus, $\Mbar^a(S)=\Tbar^a(S)/\G(S)$ is compact.
\end{theorem}
For homogeneity of notation, we will call
$\ol{\bm{W}}^a:=\Phi^{-1}:\Tbar^a(S)\rar|\Af(S)|\times[0,\infty]$.

In order to prove
Theorem~\ref{thm:phi},
we need a few preliminary results.
\begin{proposition}\label{prop:j-continuous}
The map $\Teich(S)\hra
\Tbar^a(S)$ extends to a continuous
$\hat{j}:\That^{vis}(S)\hra\Tbar^a(S)$.
\end{proposition}
\begin{proof}
The continuity of $\hat{j}$ follows from Lemma~\ref{lemma:uniform}.
Moreover, Lemma~\ref{lemma:t} and Lemma~\ref{lemma:lambda}
assure that
the $t$-lengths separate the points
of $\That^{vis}(S)$ and so $\hat{j}$ is injective.
\end{proof}

\begin{lemma}
Let $[f_m:S\rar\Si_m]$ be a sequence in $\Teich(S)$.
\begin{itemize}
\item[(a)]
$\|t_\bullet(f_m)\|_\infty\rar 0$ if and only if $\bm{\Ll}(f_m)\rar 0$
\item[(b)]
$\|t_\bullet(f_m)\|_\infty\rar\infty$ if and only if
$\bm{\Ll}(f_m)\rar\infty$
\item[(c)]
$\exists M>0$ such that $1/M\leq \|t_\bullet(f_m)\|_\infty\leq M$
if and only if $\exists M'>0$ such that
$1/M'\leq \bm{\Ll}(f_m)\leq M'$.
\end{itemize}
\end{lemma}
\begin{proof}
Because $w_{sp}(\a,f_m)\leq t_\a(f_m)$ for $\a\in\bm{A}_{sp}(f_m)$,
we conclude
\[
(6g-6+3n)\|t_\bullet(f_m)\|_\infty\geq 2\bm{\Ll}(f_m)
\]
By the collar lemma, $\ell_\a(f_m)\geq T(\bm{\Ll}(f_m))/2$
for all $\a\in\Aa(S)$ and so $\|t_\bullet(f_m)\|\leq T(T(\bm{\Ll}(f_m))/2)$.
\end{proof}

\begin{lemma}
The map $\Phi$ is continuous and injective.
\end{lemma}
\begin{proof}
Injectivity of $\Phi$ is immediate.

As we already know that $\hat{j}$ is continuous,
consider a sequence
$\{f_m:S\rar\Si_m\}\subset\Teich(S)$ such that $\bm{W}(f_m)\rar
w\in|\Af(S)|\times\{\infty\}$,
where $\bm{A}:=\mathrm{supp}(w)=\{\a_0,\dots,\a_k\}$, and
assume that $w(\a_0)\geq
w(\a_i)$ for every $1\leq i\leq k$.
Thus,
\[
\sup_{\b\in\bm{A}_{sp}(f_m)\setminus\bm{A}}
\frac{w_{sp}(\b,f_m)}{w_{sp}(\a_0,f_m)} \rar 0
\]

We want to show that $j(f_m)\rar ([w],\infty)$; equivalently,
that for every subsequence of $(f_m)$ (which we will
still denote by $(f_m)$) we can extract
a further subsequence that converges to $([w],\infty)$.

Because of Equation~\ref{eq:hex}
(applied to any maximal system of arcs $\bm{A}'$
that contains $\bm{A}$), $\ell_{\a_i}(f_m)\rar 0$
for all $\a_i\in\bm{A}$.

The collar lemma ensures that $\exists \d>0$ such that
two simple closed geodesics of length $\leq \d$
in a closed hyperbolic surface cannot intersect each other.
Thus, $\bm{A}\subseteq \bm{A}_{sp}(f_m)$ for $m$ large.

{\it Claim: for all $\b\notin\bm{A}$, the ratio
$t_\b(f_m)/t_{\a_0}(f_m)\rar 0$ uniformly.}

By contradiction, suppose $\exists \eta>0$
and $\{\b_m\}\subset\Aa(S)\setminus\bm{A}$
such that $t_{\b_m}(f_m)/t_{\a_0}(f_m)\geq \eta$.
Thus, $\ell_{\b_m}(f_m)\rar 0$ and $\b_m\in\bm{A}_{sp}(f_m)$.
By Equation~\ref{eq:hex},
\[
\sinh(w_{sp}(\ora{\b_m},f_m)) =
\frac{s_x^2+s_y^2-s_{\b_m}^2}{2s_x s_y\sqrt{s_{\b_m}^2-1}}
\approx \frac{s_x^2+s_y^2-1}{s_x s_y \ell_{\b_m}}\geq
\frac{1}{2\ell_{\b_m}}
\]
Thus, asymptotically 
$w_{sp}(\b_m,f_m)\geq 2\log(1/\ell_{\b_m}(f_m))
\approx t_{\b_m}(f_m)$.
As $w_{sp}(\a_0,f_m)\approx t_{\a_0}(f_m)$,
we conclude that $w_{sp}(\b_m,f_m)/w_{sp}(\a_0,f_m)\geq \eta/2$
for $m$ large. This contradiction proves the claim.

Given a small $\e>0$, we pick a $\d>0$ as above such that
$\cosh(\d/2)^2<1+2\e$.
If $\liminf \ell_\b(f_m)<\d$ for $\b\notin\bm{A}$,
then we can extract a subsequence
of $(f_m)$ such that $\ell_\b(f_m)<\d$ for large $m$.
Again, we can assume that
$\b$ belongs to $\bm{A}_{sp}(f_m)$ for large $m$
and so also to $\bm{A}$.
Clearly, we can only add a finite number of $\b$'s to
$\bm{A}$ and so we extract a subsequence
only a finite number of times.

Again up to subsequences, we can assume that for every
$\a\in\bm{A}$ either: $\ell_\a(f_m)\in(\d',\d)$
or $\ell_\a(f_m)\rar 0$.
If $\ell_\a(f_m)\in(\d',\d)$,
then clearly
$t_{\a}(f_m)/t_{\a_0}(f_m)\rar 0$.
Instead, if $\lim \ell_\a(f_m)\rar 0$,
then Equation~\ref{eq:hex} gives
\[
\cos(\g(\ora{\a},f_m))\geq 1-s_\a^2/2\geq 1/2-\e
\]
for large $m$. This implies that
\[
w_{sp}(\a,f_m)\approx\log(16\cos(\g(\ora{\a},f_m))
\cos(\g(\ola{\a},f_m)))-2\log(\ell_{\a}(f_m))
\]
and so $\dis
\frac{t_{\a}(f_m)}{t_{\a_0}(f_m)}\approx
\frac{\log(\ell_{\a}(f_m))}{\log(\ell_{\a_0}(f_m))}
\approx
\frac{w_{sp}(\a,f_m)}{w_{sp}(\a_0,f_m)}
\rar
\frac{w(\a)}{w(\a_0)}$.
\end{proof}

\begin{proposition}
$\G(S)$ acts on $\Tbar^a(S)$ by isometries
and with discrete orbits. Hence,
$\Mbar^a(S)=\Tbar^a(S)/\G(S)$ is Hausdorff.
\end{proposition}
\begin{proof}
Suppose $t\cdot g_m\rar t$, with $t\in\Tbar^a(S)$
and $g_m\in\G(S)$. Consider a sequence
$[f_m:S\rar\Si_m]$ such that $j(f_m)\rar t$
in $\Tbar^a(S)$.

{\bf Case 1: $\|t\|_\infty<\infty$.}
Passing to a subsequence, $[f_m]\cdot h_m
\rar [f_\infty:S\rar\Si_\infty]\in\That^{vis}(S)$
for suitable $h_m\in \G(S)$.
Thus, $\hat{j}(f_\infty)\cdot h_m^{-1}\rar t$.

Assume first that $\|t\|_\infty>0$ and so $\Si_\infty$
does not have $n$ cusps.

Let $\bm{A}=\bm{A}_{sp}(f_\infty)$, so that it
is supported on $f_\infty^{-1}(\Si_{\infty,+})$ and
$\ell_\a(f_\infty)<\infty$
for all $\a\in\bm{A}$.
Because the length
spectrum of finite arcs in $\Si_\infty$ is discrete
(with finite multiplicities)
and $\hat{j}(f_\infty)\cdot h_m^{-1}$ is a Cauchy sequence,
there exists an integer $M$ such that
(up to subsequences)
$h_m^{-1}$ fixes $\bm{A}$ for all $m\geq M$.
Thus, $h_m$ is the composition of a diffeomorphism
of $f_\infty^{-1}(\Si_{\infty,-})$
and an isometry of $f_\infty^{-1}(\Si_{\infty,+})$
(with the pull-back metric) for $m\geq M$.
Hence, $t=\hat{j}(f_\infty)\cdot h_m^{-1}$ for $m\geq M$
and so $t=\hat{j}(\hat{f}_\infty)$ for some
$\hat{f}_\infty:S\rar\Si_\infty$.
Similarly, $\hat{j}(\hat{f}_\infty)\cdot g_m\rar
\hat{j}(\hat{f}_\infty)$ and so $g_m$ is
the composition of a diffeomorphism
of $\hat{f}_\infty^{-1}(\Si_{\infty,-})$
and an isometry of $\hat{f}_\infty^{-1}(\Si_{\infty,+})$
for large $m$. Hence, $t\cdot g_m$ cannot accumulate at $t$.

Assume now that $\|t\|_\infty=0$, so that $\Si_\infty$
has $n$ cusps.

It follows from the classical case that
the spectrum of the finite reduced lengths
(and so of the finite $\l$-lengths) of $(\Si_\infty,\up)$
is discrete and with finite multiplicities.
Because $[t_\bullet(f_\infty)]=[\l_\bullet(f_\infty)]$,
we can conclude as in the previous case.

{\bf Case 2: $\|t\|_\infty=\infty$.}\\
Let $w^{m}=\bm{W}(f_m)$. Up to subsequences,
$w^{m}\cdot h_m\rar w^{\infty}$ in $|\Af(S)|\times[0,\infty]$
for suitable $h_m\in\G(S)$ and $w^{\infty}\in |\Af(S)|\times
\{\infty\}$.
As before, $\Phi(w^{\infty})\cdot h_m^{-1}\rar t$ in $\Tbar^a(S)$.
Because $w^{\infty}$ has finite support,
$t=\Phi(w^{\infty})\cdot h_m^{-1}$ for $m\geq M$
and so $t=\Phi(\hat{w}^\infty)$, for some $\hat{w}^\infty
\in|\Af(S)|\times\{\infty\}$.
Thus, $\Phi(\hat{w}^\infty)\cdot g_m \rar\Phi(\hat{w}^\infty)$
and $g_m$ is the composition of a diffeomorphism
of $S_-$ and an isometry of $S_+$,
where $S_+$ (resp. $S_-$) is the $\hat{w}^{\infty}$-visible
(resp. invisible) subsurface of $S$.
Hence, $t\cdot g_m$ cannot accumulate at $t$.
\end{proof}
\begin{proof}[Proof of Theorem~\ref{thm:phi}]
In order to apply Lemma~\ref{lemma:isom}(b), we only
need to prove that
$\Phi':|\Af(S)|/\G(S)\times[0,\infty]\rar \Mbar^a(S)$
is a homeomorphism.

We already know that $\Phi'$ is continuous, injective.
Moreover, its image contains $\M(S)$ which is dense
in $\Mbar^a(S)$.
As $|\Af(S)|/\G(S)$ is compact and $\Mbar^a(S)$
is Hausdorff, the map $\Phi'$ is closed and so
it is also surjective. Hence, $\Phi'$ is a
homeomorphism.
\end{proof}

\begin{corollary}
$\hat{j}$ is a homeomorphism onto the finite
part of $\Tbar^a(S)$.
\end{corollary}
\end{subsection}
%
%
\begin{subsection}{The extended Teichm\"uller space}
We define the {\it extended Teichm\"uller space}
$\Ttil(S)$ to be
\[
\Ttil(S):=\Tbar^{WP}(S)\cup |\Af(S)|_\infty
\]
where $|\Af(S)|_\infty$ is just a copy of $|\Af(S)|$.

Clearly, there is map $\mathrm{Bl}_0 \Ttil(S)\rar
\Tbar^a(S)$, which identifies visibly equivalent
surfaces of $\That(S)\subset\mathrm{Bl}_0\Ttil(S)$.

We define a topology on $\Ttil(S)$ by
requiring that $\Tbar^{WP}(S)\hra \Ttil(S)$ and
$|\Af(S)|_\infty\hra\Ttil(S)$ are homeomorphisms
onto their images, that $\Tbar^{WP}(S)\subset\Ttil(S)$
is open and
we declare that a sequence
$\{f_m\}\subset \Tbar^{WP}(S)$ is converging to
$w\in|\Af(S)|_\infty$ if and only if
$\bm{W}(f_m)\rar (w,\infty)$ in $|\Af(S)|\times(0,\infty]$.

Notice that $\Mtil(S):=\Ttil(S)/\G(S)$
is an orbifold with corners, which acquires some
singularities at infinity. In fact, $\Mtil(S)$ is homeomorphic
to $\Mbar^{WP}(R,x)\times\Delta^{n-1}\times[0,\infty]/\!\!\sim$, where
$(R,x)$ is a closed $x$-marked surface such that $S\simeq R\setminus x$
and $(R',\up',t')\sim(R'',\up'',t'')$ $\iff$ $t'=t''=\infty$ and
$(R',\up')$ is visibly equivalent to $(R'',\up'')$.
\end{subsection}
%
%
\end{section}
%
%
\begin{section}{Weil-Petersson form and circle actions}\label{sec:wp}
%
%
\begin{subsection}{Circle actions on moduli spaces}
Let $S$ be a compact surface of genus $g$ with boundary
components $C_1,\dots,C_n$ (assume as usual that $2g-2+n>0$).
Let $v_i$ be a point of $C_i$ and set $v=(v_1,\dots,v_n)$.

We denote by $\mathrm{Diff}_+(S,v)$ the group of orientation-preserving
diffeomorphisms of $S$ that fix $v$ pointwise
and by $\mathrm{Diff}_0(S,v)$ its connected component of the identity.

The Teichm\"uller space $\Teich(S,v)$ is the space of
hyperbolic metrics on $S$ up to action of $\mathrm{Diff}_0(S,v)$
and the mapping class
group of $(S,v)$ is
$\G(S,v)=\mathrm{Diff}_+(S,v)/\mathrm{Diff}_0(S,v)$.
Thus, $\M(S,v)=\Teich(S,v)/\G(S,v)$ is the resulting moduli space.

Clearly, $\R^n$ acts on $\Teich(S,v)$ by Fenchel-Nielsen twist
(with unit angular speed)
around the boundary components and $\Teich(S,v)/\R^n=\Teich(S)$.
Similarly, the torus $\To^n=(\R/2\pi\Z)^n$
acts on $\M(S,v)$ and the quotient is $\M(S,v)/\To^n=\M(S)$.

Mimicking what done for $\Tbar^{WP}(S)$, we can define
an augmented Teichm\"uller space $\Tbar^{WP}(S,v)$
and an action of $\R^n$ on it.
However, we want to be a little more careful
and require that a marking $[f:S\rar\Si]\in\Tbar^{WP}(S,v)$
that shrinks $C_i$ to a cusp $y_i\in\Si$ is
smooth with $\mathrm{rk}(df)=1$ at $C_i$, so that
$f$ identifies $C_i$ to the sphere tangent bundle
to $ST_{\Si,y_i}$ and $v_i$ to a point in $ST_{\Si,y_i}$.

Thus, $\Tbar^{WP}(S,v)\rar\Tbar^{WP}(S)$ is an $\R^n$-bundle
and $\Mbar(S,v)\rar \Mbar(S)$ is a $\To^n$-bundle,
which is a product $L_1\times\dots\times L_n$
of circle bundles $L_i$ associated to $v_i\in C_i$.

If one wish, one can certainly lift the action to
$\That(S,v)=\mathrm{Bl}_0\Tbar^{WP}(S,v)$.

As for the definition of
$\That^{vis}(S,x)=\That(S,x)/\!\!\sim_{vis}$,
we declare $[f_1:S\rar\Si_1],
[f_2:S\rar\Si_2]\in\That(S,v)$ to be
visibly equivalent
if $\exists\, [f:S\rar\Si]\in\That(S,v)$ and
maps $h_i:\Si\rar \Si_i$ for $i=1,2$ such that
$h_i$ restricts to an isometry $\Si_+\rar\Si_{i,+}$
and $h_i\circ f\simeq f_i$ (for $i=1,2$) through
homotopies that fix $f^{-1}(\ol{\Si}_+)\cap v$ (but
not necessarily $f^{-1}(\ol{\Si}_-)\cap v$).

This means that, if $[f:S\rar\Si]\in\That^{vis}(S,v)$
has $f(C_i)\subset\Si_-$, then $[f]$ does not record
the exact position of the point $v_i\in C_i$.
In other words, the $i$-th component of $\R^n$
acts trivially on $[f]$.
\end{subsection}
%
%
\begin{subsection}{The arc complex of $(S,v)$}
Let $\Aa(S,v)$ to be the set of nontrivial isotopy
classes of simple arcs in $S$
that start and end at $\pa S\setminus v$ and let
$\b_i$ be a (fixed) arc from $C_i$ to $C_i$ that
separates $v_i$ from the rest of the surface.

A subset $\bm{A}=\{\b_1,\dots,\b_n,\a_1,\dots,\a_k\}\subset\Aa(S,v)$ is
a $k$-system of arcs on $(S,v)$ if $\b_1,\dots,\b_n,\a_1,\dots,\a_k$ admit
disjoint representatives. The arc complex $\Af(S,v)$ 
is the set of systems of arcs on $(S,v)$.
A point in $\Af(S,v)$ can be represented as a sum
$\sum_j w_j \a_j$, provided we remember the $\b_i$'s
(that is, as $\sum_j w_j\a_j+\sum_i 0\b_i$)
or as a function $w:\Aa(S,v)\rar\R$.

We can define $\Ao(S,v)\subset\Af(S,v)$ to be the subset
of simplices representing systems of arcs that cut $S$
into a disjoint union of discs and annuli homotopic to
some boundary component.

Remark that there is a natural map $\Af(S,v)\rar\Af(S)$,
induced by the inclusion $S\setminus v\hra S$ and
that forgets the $\b_i$'s,
and so a simplicial map $|\Af(S,v)|\rar|\Af(S)|$.

We can also define a suitable map
$\wh{\bm{W}}_v$ for the pointed surface $(S,v)$
in such a way that the following diagram commutes.
\[
\xymatrix{
\That^{vis}(S,v) \ar[rr]^{\wh{\bm{W}}_v\quad} \ar[d]
&& |\Af(S,v)|\times [0,\infty) \ar[d] \\
\That^{vis}(S)   \ar[rr]^{\wh{\bm{W}}\quad}
&& |\Af(S)|  \times [0,\infty)
}
\]
Let $[f:S\rar\Si]\in\That^{vis}(S,v)$.
If we consider it as a point of $\That^{vis}(S)$,
then $\wh{\bm{W}}(f)$ is a system of arcs in $S$.

For every $i=1,\dots,n$ such that $f(C_i)\in\Si_+$,
consider the geodesic $\rho_i\subset\Si$ coming
out from $f(v_i)$ and perpendicular to $f(C_i)$
(if $f(C_i)$ is a cusp, let $\rho_i$ be the geodesic
originating at $f(C_i)$ in direction $f(v_i)$).
Call $z_i$ the point where $\rho_i$ first meets
the spine of $\Si$ and $e_i$ an infinitesimal portion
of $\rho_i$ starting at $z_i$ and going towards $f(C_i)$.

Define $\mathrm{Sp}(\Si,f(v))$ to be the
one-dimensional CW-complex obtained from
$\mathrm{Sp}(\Si)$
by adding the vertices $z_i$ (in case $z_i$ was
not already a vertex) and the infinitesimal edges $e_i$.
Consequently, we have a well-defined
system of arcs $\bm{A}_{sp}(\Si,f(v))$ dual to
$\mathrm{Sp}(\Si,f(v))$ and widths $w_{sp,f(v)}$, in which
the arc dual to $e_i$ plays the role of $f(\b_i)$
(which thus has zero weight).

We set $\wh{\bm{W}}_v(f)=f^* w_{sp,f(v)}$.

The following is an immediare consequence
of Theorem~\ref{prop:w}.

\begin{proposition}
The map $\wh{\bm{W}}_v$
is a $\G(S,v)$-equivariant homeomorphism.
\end{proposition}

We can make $\R^n$ act on $|\Af(S,v)|$ via $\wh{\bm{W}}_v$
and so on $|\Af(S,v)|\times[0,\infty]$.
Thus, the action also prolongs to the extended Teichm\"uller space
$\Ttil(S,v):=\Tbar^{WP}(S,v)\cup|\Af(S,v)|_{\infty}$.
\end{subsection}
%
%
%
\begin{subsection}{Weil-Petersson form}
Chosen a maximal set of simple closed curves
$\bm{\g}=\{\g_1,\dots,\g_{6g-6+2n},C_1,\dots,C_n\}$ on $S$,
we can define a symplectic form $\omega_v$ on $\Teich(S,v)$
by setting
\[
\omega_v=\sum_{i=1}^{6g-6+2n} d\ell_i\wedge d\tau_i+
\sum_{j=1}^n dp_j\wedge d t_j
\]
where $t_j=p_j\th_j/2\pi$ is the twist parameter at $C_j$.
As usual, $\omega_v$ does not depend on the choice
of $\bm{\g}$ and it descends to $\M(S,v)$.
Its independence of the particular
Fenchel-Nielsen coordinates permits to extend
$\omega_v$ to a symplectic form on $\Mbar(S,v)$.

Moreover, as $\omega_v(dp_j,-)=\pa/\pa t_j$,
the twist flow on $\Mbar(S,v)$ is Hamiltonian
and the associated moment map is exactly
$\mu=(p_1^2/2,\dots,p_n^2/2)$.
Thus, the leaves $(\Mbar(S)(\up),\omega_{\up})$
are exactly the symplectic
reductions of $(\Mbar(S,v),\omega_v)$ with respect to
the $\To^n$-action. 

As remarked by Mirzakhani \cite{mirzakhani:volumes},
it follows by standards results of symplectic geometry
that there is a symplectomorphism
$\Mbar(S)(\up)\rar \Mbar(S)(0)$ which pulls
$[\omega_0]+\sum_{i=1}^n \frac{p_i^2}{2} c_1(L_i)$
back to $[\omega_{\up}]$.

Penner has provided a beautiful formula for $\omega_0$
in term of the $\tilde{a}$-coordinates.

\begin{theorem}[\cite{penner:wp}]
Let $\bm{A}$ be a maximal system of arcs on $S$.
If $\pi:\Teich(S)(0)\times\R_+^n\rar\Teich(S)(0)$
is the projection onto the first factor,
then
\[
\pi^*\omega_0=-\frac{1}{2}\sum_{t\in T}
\left(
d\tilde{a}_{t_1}\wedge d\tilde{a}_{t_2}+
d\tilde{a}_{t_2}\wedge d\tilde{a}_{t_3}+
d\tilde{a}_{t_3}\wedge d\tilde{a}_{t_1}
\right)
\]
where $T$ is the set of ideal triangles in $S\setminus\bm{A}$,
the triangle $t\in T$ is bounded by the (cyclically ordered) arcs
$(\a_{t_1},\a_{t_2},\a_{t_3})$.
\end{theorem}

The whole $\Teich(S)$ is naturally a Poisson manifold
with the Weil-Petersson pairing $\eta$ on the cotangent
bundle, whose symplectic leaves
are the $\Teich(S)(\up)$.
A general formula expressing $\eta$ in term of lengths of
arcs and widths is given by the following.

\begin{theorem}[\cite{mondello:poisson}]\label{thm:poisson}
Let $\bm{A}$ be a maximal system of arcs on $S$.
Then
\[
\eta=\frac{1}{4}\sum_{k=1}^n
\sum_{\substack{y_i\in\a_i\cap C_k \\ y_j\in\a_j\cap C_k}}
\frac{\sinh(p_k/2-d_{C_k}(y_i,y_j))}{\sinh(p_k/2)}
\frac{\pa}{\pa a_i}\wedge\frac{\pa}{\pa a_j}
\]
where $d_{C_k}(y_i,y_j)$ is the length
of the geodesic running from $y_i$ to $y_j$ along $C_k$ in the positive
direction.
\end{theorem}

In order to understand the limit for large $\up$, it makes
sense to rescale the main quantities as $\tilde{w}_i=(\bm{\Ll}/2)^{-1}w_i$,
$\tilde{\omega}=(1+\bm{\Ll}/2)^{-2}\omega$
and $\tilde{\eta}=\left(1+\bm{\Ll}/2\right)^2\eta$.

\begin{lemma}[\cite{kontsevich:intersection}]
The class $[\tilde{\omega}_\infty]\in H^2_{\G(S)}(|\Af(S)|)$
is represented by a piecewise linear $2$-form on $|\Af(S)|$
whose dual can be written (on the maximal simplices) as
\[
\tilde{H}=\frac{1}{2}\sum_{r}\left(
\frac{\pa}{\pa \tilde{w}_{r_1}}\wedge\frac{\pa}{\pa\tilde{w}_{r_2}}+
\frac{\pa}{\pa \tilde{w}_{r_2}}\wedge\frac{\pa}{\pa\tilde{w}_{r_3}}+
\frac{\pa}{\pa \tilde{w}_{r_3}}\wedge\frac{\pa}{\pa\tilde{w}_{r_1}}
\right)
\]
where $r$ ranges over all (trivalent) vertices of the ribbon graph
represented by a point in $|\Ao(S)|$ and $(r_1,r_2,r_3)$ is the (cyclically)
ordered triple of edges incident at $r$.
\end{lemma}

The above result admits a pointwise sharpening as follows.

\begin{theorem}[\cite{mondello:poisson}]
The bivector field
$\ti{\eta}$ extends over $\Ttil(S)$ and,
on the maximal simplices of $|\Af(S)|_\infty$, we have
\[
\tilde{\eta}_\infty=\tilde{H}
\]
pointwise.
\end{theorem}

Thus, we have a description of the degeneration of $\eta$
when the boundary lengths of the hyperbolic surface become very large.
\end{subsection}
%
%
\end{section}
%
%
\begin{section}{From surfaces with boundary to
pointed surfaces}\label{sec:grafting}
%
%
\begin{subsection}{Ribbon graphs}
Let $S$ be a compact oriented surface of genus $g$ with boundary
components $C_1,\dots,C_n$ and assume that $\chi(S)=2-2g-n<0$.
Let $\bm{A}=\{\a_0,\dots,\a_k\}\in\Af(S)$
be a system of arcs in $S$ and
$S_+$ the corresponding visible subsurface of $S$.

If $\ora{\a}$ is an oriented arc supported on $\a$,
then we will refer
to the symbol $\ora{\a}^*$ as to the
{\it oriented edge dual to $\ora{\a}$}.
\begin{remark}
If $S$ carries a hyperbolic metric and $\bm{A}$ is its spinal
system of arcs, then $\ora{\a}^*$ must be considered
the edge of the spine dual to $\a$ and oriented in such a way
that, at the point $\ora{\a}^*\cap\ora{\a}$ (unique, up to prolonging
$\ora{\a}^*$) the tangent
vectors $\langle v_{\ora{\a}^*}, v_{\ora{\a}}\rangle $ form a positive basis
of $T_p S$.
\end{remark}
Let $E(\bm{A}):=\{\ora{\a}^*,\ola{\a}^*\,|\,\a\in\bm{A}\}$
and define the following
operators $\s_0,\s_1,\s_{\infty}$ on $E(\bm{A})$:
\begin{itemize}
\item[(1)]
$\s_1$ reverses the orientation of each arc
(i.e. $\s_1(\ora{\a}^*)=\ola{\a}^*$)
\item[($\infty$)]
if $\ora{\a}$ ends at $x_{\a}\in C_i$,
then $\s_{\infty}(\ora{\a}^*)$
is dual to
the oriented arc $\ora{\b}$ that ends at $x_\b\in C_i$,
where $x_\b$ comes just {\it before} $x_\a$ according to the
orientation induced on $C_i$ by $S$
\item[(0)]
$\s_0$ is defined by $\s_0=\s_1\s_{\infty}^{-1}$.
\end{itemize}
If we call $E_i(\bm{A})$ the orbits of $E(\bm{A})$ under the
action of $\s_i$, then
\begin{itemize}
\item[(1)]
$E_1(\bm{A})$ can be identified with $\bm{A}$
\item[($\infty$)]
$E_{\infty}(\bm{A})$ can be identified with
the subset of the boundary components of $S$ that
belong to $S_+$
\item[(0)]
$E_0(\bm{A})$ can be identified to the set of connected
components of $\dis S_+\setminus\bm{A}$.
\end{itemize}
\end{subsection}
%
%
\begin{subsection}{Flat tiles and Jenkins-Strebel
differentials}\label{sec:flat-tiles}
Keeping the notation as before,
let $f:S\rar\hat{S}$ be the topological type of $\bm{A}$
(see Section~\ref{sec:stable}).

For every system of weights $w$ supported on
$\bm{A}$, the surface $\hat{S}_+$ can be endowed with a flat metric
(with conical singularities) in the following way.

Every component $\hat{S}_{i,+}$ of $\hat{S}_+$ is quasi-filled
by the arc system $f(\bm{A})\cap\hat{S}_{i,+}$. As we can
carry on the construction componentwise, we can assume
that $\bm{A}$ quasi-fills $S$.

In this case, we consider the flat {\it tile}
$\dis T=[0,1]\times[0,\infty]\Big/[0,1]\times\{\infty\}$
and we call {\it point at infinity} the class
$[0,1]\times\{\infty\}$.
Moreover, we define $\Si:=\bigcup_{\ora{e}\in E(\bm{A})}T_{\ora{e}}/\!\!\sim$,
where $T_{\ora{e}}:=T\times\{\ora{e}\}$ and
\begin{itemize}
\item
$(u,0,\ora{e})\sim(1-u,0,\ola{e})$ for all $\ora{e}\in E(\bm{A})$
and $u\in[0,1]$
\item
$(1,v,\ora{e})\sim(0,v,\s_\infty(\ora{e}))$
for all $\ora{e}\in E(\bm{A})$ and $v\in[0,\infty]$.
\end{itemize}
We can also define an embedded graph $G\subset\Si$
by gluing the segments $[0,1]\times\{0\}\subset T$
contained in each tile. Thus, we can identify $\a^*$
with an edge of $G$ for every $\a\in\bm{A}$.

It is easy to check that there is a homeomorphism
$\hat{S}\rar\Si$, well-defined up to isotopy,
that takes boundary components to points at infinity
or to vertices.

Moreover, for every $\ora{\a}^*\in E(\bm{A})$,
we can endow $T_{\ora{\a}^*}$ with
the quadratic differential $dz^2$, where
$z=w(\a)u+iv$.
These quadratic differentials
glue to give a global $\varphi$ (and so a conformal
structure on the whole $\Si$), which 
has double poles with negative quadratic residue
at the points at infinity and
is holomorphic elsewhere,
with zeroes of order $k-2$ at the $k$-valent vertices of $G$.
Furthermore, $\a^*$ has length $w(\a)$
with respect to the induced flat metric $|\varphi|$.

Finally, the {\it horizontal trajectories} of $\varphi$
(that is, the curves along which $\varphi$ is
positive-definite) are: either closed circles that wind
around some point at infinity, or edges of $G$.

Thus, $\varphi$ is a
{\it Jenkins-Strebel quadratic differential} and $G$
is its {\it critical graph}, i.e. the union of all
horizontal trajectories that hit some zero or some
pole of $\varphi$.

If $\bm{A}$ does not quasi-fill $S$, then we will define
the Jenkins-Strebel differential componentwise,
by setting it to zero on the invisible components.

See \cite{harer:virtual}, \cite{kontsevich:intersection},
\cite{looijenga:cellular} and \cite{mondello:survey}
for more details.
\end{subsection}
%
%
\begin{subsection}{HMT construction}
We begin by recalling the following result of Strebel.
\begin{theorem}[\cite{strebel:67}]\label{thm:strebel}
Let $R'$ be a compact Riemann surface of genus $g$
with $x'=(x'_1,\dots,x'_n)$ distinct points such that
$n\geq 1$ and $2g-2+n>0$. For every $(p_1,\dots,p_n)\in
\R^n_{\geq 0}$ (but not all zero),
there exists a unique (nonzero) quadratic differential
$\varphi$ on $R'$ such that
\begin{itemize}
\item
$\varphi$ is holomorphic on $R'\setminus x'$
\item
horizontal trajectories of $\varphi$ are either circles
that wind around some $x'_i$ or closed arcs between
critical points
\item
the critical graph $G$ of $\varphi$ cuts $R'$ into semi-infinite
flat cylinders (according to the metric $|\varphi|$),
whose circumferences are closed trajectories
\item
if $p_i=0$, then $x'_i$ belongs to the critical graph
\item
if $p_i>0$, then
the cylinder around $x'_i$ has circumference length $p_i$.
\end{itemize}
\end{theorem}

Notice that the graph $G$ plays a role analogous to
the spine of a hyperbolic surface.
In fact, given a point $[f:R\rar R']\in\Teich(R,x)$
and $(p_1,\dots,p_n)\in\Delta^{n-1}$,
we can consider the unique $\varphi$
given by the theorem above and
the system of arcs $\bm{A}\in\Ao(R,x)$
such that $f(\bm{A})$ is dual to the critical graph $G$ of $\varphi$,
and we can
define the width $w(\a)$ to be the $|\varphi|$-length
of the edge $\a^*$ of $G$ dual to $\a\in\bm{A}$.

\begin{theorem}[Harer-Mumford-Thurston \cite{harer:virtual}]
The map $\dis\Teich(R,x)\times\Delta^{n-1}
\lra |\Ao(R,x)|$ just constructed is a
$\G(R,x)$-equivariant homeomorphism.
\end{theorem}

Clearly, if $R'$ is a stable Riemann surface, then
the theorem can be applied on every visible component
of $R'$ (i.e. on every component that contains some $x'_i$
with $p_i>0$) and $\varphi$ can be extended by zero
on the remaining part of $R'$. Hence, we can extend
the previous map to
\[
\bm{W}_{HMT}:\Tbar^{vis}(R,x)\times\Delta^{n-1}
\lra |\Af(R,x)|
\]
which is also a $\G(R,x)$-equivariant homeomorphism
(see, for instance, \cite{looijenga:cellular}
and \cite{mondello:survey}).

The purpose of the following sections is to relate
this $\bm{W}_{HMT}$ to the spine construction.
\end{subsection}
%
%
%
\begin{subsection}{The grafting map}
Given a hyperbolic surface $\Si$ with boundary components
$C_1,\dots,C_n$, we can
{\it graft a semi-infinite flat cylinder at each $C_i$} of circumference
$p_i=\ell(C_i)$. The result is a surface $\gr8(\Si)$ with
a $C^{1,1}$-metric, called the {\it Thurston metric}
(see \cite{scannell-wolf:grafting} for the case
of a general lamination, or
\cite{kulkarni-pinkall:canonical} a
higher dimensional analogues).
If $\Si$ has cusps,
we do not glue any cylinder at the cusps of $\Si$.
Notice that $\gr8(\Si)$ has the conformal type of a punctured
Riemann surface and it will be sometimes regarded as
a closed Riemann surface with marked points.
\begin{notation}
Choose a closed surface $R$ with distinct marked points
$x=(x_1,\dots,x_n)\subset R$ and an identification
$R\setminus x\cong\gr8(S)$ such that $x_i$ corresponds to $C_i$.
Clearly, we can identify $\Af(S)\cong\Af(R,x)$ and
$\G(S)\cong\G(R,x)$.
\end{notation}
We use the grafting construction to define a map
\[
(\gr8,\Ll):\Ttil(S) \lra \Tbar^{WP}(R,x)\times\Delta^{n-1}\times
[0,\infty]/\!\!\sim
\]
where $\sim$ identifies $([f_1],\up,\infty)$ and $([f_2],\up,\infty)$
if $([f_1],\up)$ and $([f_2],\up)$ are visibly equivalent.

We set $\gr8(f:S\rar\Si):=[\gr8(f):R\rar\gr8(\Si)]$,
on the bounded part $\Tbar^{WP}(S)\subset\Ttil(S)$.
On the other hand, if $\tilde{w}\in|\Af(S)|_\infty$
represents
a point at infinity of $\Ttil(S)$, then we define
$(\gr8,\Ll)(\tilde{w}):=
(\bm{W}_{HMT}^{-1}(\tilde{w}),\infty)$.

\begin{theorem}\label{thm:grafting}
The map $(\gr8,\Ll)$ is a $\G(S)$-equivariant
homeomorphism
that preserves the topological types and whose
restriction to each topological stratum of the finite part
and to each simplex of $|\Af(S)|_\infty$ is a real-analytic
diffeomorphism.
\end{theorem}

\begin{corollary}
(a) The induced $\Mtil(S)\lra \Mbar(R,x)
\times\Delta^{n-1}\times[0,\infty]/\!\!\sim$ is a homeomorphism,
which is real-analytic on $\Mhat(S)$ and piecewise real-analytic
on $|\Af(S)|_\infty/\G(S)$.

(b)
Let $\That^{vis}(R,x)$ (resp. $\Mhat^{vis}(R,x)$)
be obtained from $\Tbar^{WP}(R,x)\times\Delta^{n-1}$
(resp. $\Mbar^{WP}(R,x)\times\Delta^{n-1}$) by identifying
visibly equivalent surfaces.
Then, the induced $\Tbar^a(S)\lra \That^{vis}(R,x)
\times[0,\infty]$
and $\Mbar^a(S)\lra\Mhat^{vis}(R,x)\times[0,\infty]$
are homeomorphisms.
\end{corollary}

We can summarize our results in the following commutative
diagram
\[
\xymatrix{
\That^{vis}(R,x)\times[0,\infty]
\ar[rrd]_{\Psi} &&
\Tbar^a(S) \ar[ll]_{\qquad\qquad (\gr8,\Ll)}
\ar[d]^{\ol{\bm{W}}^a} \\
&& |\Af(R,x)|\times[0,\infty]
}
\]
in which $\Psi=\ol{\bm{W}}^a\circ(\gr8,\Ll)^{-1}$ and
all maps are $\G(R,x)$-equivariant homeomorphisms.
For every $t\in[0,\infty]$, call $\Psi_t:\That^{vis}(R,x)
\rar|\Af(R,x)|$
the restriction of $\Psi$ to $\That^{vis}(R,x)\times\{t\}$
followed by the projection onto $|\Af(R,x)|$.

\begin{corollary}\label{cor:grafting}
$\Psi_t$ is a continuous family of $\G(R,x)$-equivariant
triangulations of $\That^{vis}(R,x)$,
whose extremal cases are Penner/Bowditch-Epstein's
for $t=0$ and Harer-Mumford-Thurston's for $t=\infty$.
\end{corollary}

In order to prove Theorem~\ref{thm:grafting}, we will show first that
$(\gr8,\Ll)$ is continuous.
Lemma~\ref{lemma:isom}(a) will ensure that it is proper.
Finally, we will prove that the restriction of $(\gr8,\Ll)$
to each stratum is bijective onto its image, and so
that $(\gr8,\Ll)$ is bijective.
\end{subsection}
%
%
\begin{subsection}{Continuity of $(\gr8,\Ll)$}
To test the continuity of $(\gr8,\Ll)$
at $q\in\Ttil(S)$, we split the problem
into two distinct cases:
\begin{enumerate}
\item
$\Ll(q)$ bounded and so $q=[f:S\rar\Si]$
\item
$\Ll(q)$ not bounded and so
$q=\tilde{w}\in|\Af(S)|_\infty$.
\end{enumerate}
%
%
\begin{subsubsection}{$\Ll(q)$ bounded.}\label{sec:bounded}
Let $\{f_m:S\rar\Si_m\}\subset\Teich(S)$
be a sequence
that converges to $[f]$, so that $\Ll(f_m)\rar\Ll(f)$.

Condition (2) of Proposition~\ref{prop:convergence}
ensures that there are
maps $\tilde{f}_m:S\rar\Si_m$ which have a standard behavior
on a neighbourhood of the thin part of $\Si_m$ and
such that the metric $\tilde{f}_m^*(g_m)$
converges to $f^*(g)$ uniformly on the complement.
Fixed some $\gr8(f):R\rar\gr8(\Si)$, define
$\hat{f}_m:R\rar\gr8(\Si_m)$ in such a way that
$F_m=\gr8(f)\circ\hat{f}_m^{-1}:\gr8(\Si_m)\rar\gr8(\Si)$
has the following properties.

If $\ell_{C_i}(f)>0$, then
let $\phi_m^i:\pa^i \Si_m\rar \pa^i\Si$ be restriction
of $f\circ\tilde{f}_m^{-1}$ to the $i$-th boundary component.
Moreover, we can give orthonormal
coordinates $(x,y)$ (with $y\geq 0$ and $x\in [0,\ell_{C_i})$)
on the $i$-th cylinder
in such a way that $C_i=\{y=0\}$ and $S$ induced on $C_i$
the orientation along which $x$ {\it decreases}.

For every $i$ such that $\ell_{C_i}(f)>0$,
we define $F_m$
to be $(x,y)\mapsto (\phi_m^i(x),y)$ on the $i$-th cylinder.

For every $i$ such that $\ell_{C_i}(f)=0$ and $\ell_{C_i}(f_m)>0$,
we can assume that $\ell_{C_i}(f_m)<1/2$
and we can consider a hypercycle $H_i\subset\Si_m\subset\gr8(\Si_m)$
parallel to the $i$-th boundary component and of length $2\ell_{C_i}(f_m)$.

We define $F_m$ to agree with
$f\circ \tilde{f}_m^{-1}$ on the portion of $\gr8(\Si_m)$
which is hyperbolic and bounded by the possible hypercycles $H_i$'s.

Finally, we extend $F_m$
outside the possible hypercycles $H_i$'s by a diffeomorphism.

Clearly, condition (5) of Proposition~\ref{prop:convergence}
for the sequence $\{\hat{f}_m\}$
and $\gr8(f)$ is verified
and so $[f_m]=[\hat{f}_m]\rar[\gr8(f)]$ in $\That(R,x)$.
\end{subsubsection}
%
%
%
\begin{subsubsection}{$\Ll(q)$ not bounded.}
Let $S_+\subset S$ be the visible subsurface determined
(up to isotopy) by $\bm{A}=\mathrm{supp}{\tilde{w}}$
and let $\{\bm{A}_i\,|\,i\in I\}$ be the set of all maximal
system of arcs of $S$ that contain $\bm{A}$.

Consider a sequence $[f_m:S\rar\Si_m]\in\Teich(S)$
that converges to $q=\tilde{w}\in|\Af(S)|_\infty\subset\Ttil(S)$
and such that
$\widehat{\bm{W}}(f_m)\in |\bm{A}_{i_m}|^\circ\times\R_+$,
with $i_m\in I$.

We must show that $\gr8(f_m)\rar
\gr8(q)=[f:R\rar\Si]$
in $\Tbar(R,x)\times\Delta^{n-1}\times[0,\infty]/\!\!\sim$,
where $f$ and $\Si$ are constructed as in Section~\ref{sec:flat-tiles}.

It will be convenient to denote by $\tilde{w}(\a,f)$ the
weight $\tilde{w}(\a)$ for every $\a\in\bm{A}$.
Moreover, we will use the notation $\tilde{w}(-,f_m)$ to denote
the normalized quantity $2w_{sp}(-,f_m)/\bm{\Ll}(f_m)$
and $\tilde{w}_m(\ora{\a},f)$
to denote
$\dis\frac{\tilde{w}(\a,f)}{\tilde{w}(\a,f_m)}\tilde{w}(\ora{\a},f_m)$.

\begin{remark}
Using Proposition~\ref{prop:convergence},
it is sufficient to show that condition (5)
for the sequence $\{\gr8(f_m)\}$ and $\gr8(q)$
is satisfied on the positive components of $\Sigma$.
As usual, we will define a sequence of homeomorphisms
$\hat{f}_m:R\rar\gr8(\Si_m)$ (that satisfies
condition (5)) by describing $F_m=f\circ\hat{f}_m^{-1}:
\gr8(\Si_m)\rar\Si$.
\end{remark}

For every $m$ and every small $\e>0$,
define the following regions of $\gr8(\Si_m)$ and of $\Si$.
\begin{center}
\begin{figurehere}
\psfrag{a}{$\ora{\a}$}
\psfrag{e}{$\ora{\a}^*$}
\psfrag{v}{$v$}
\psfrag{p}{$s$}
\psfrag{P-}{$P_-$}
\psfrag{Re-}{$R^\e_-$}
\psfrag{Qe-}{$Q^\e_-$}
\psfrag{Ue}{$\hat{U}^\e$}
\psfrag{Rh}{$\hat{R}$}
\psfrag{Rh-}{$\hat{R}_-$}
\psfrag{Sin}{$\Si_m$}
\psfrag{Sp}{$\mathrm{Sp}(\Si_m)$}
\includegraphics[width=\textwidth]{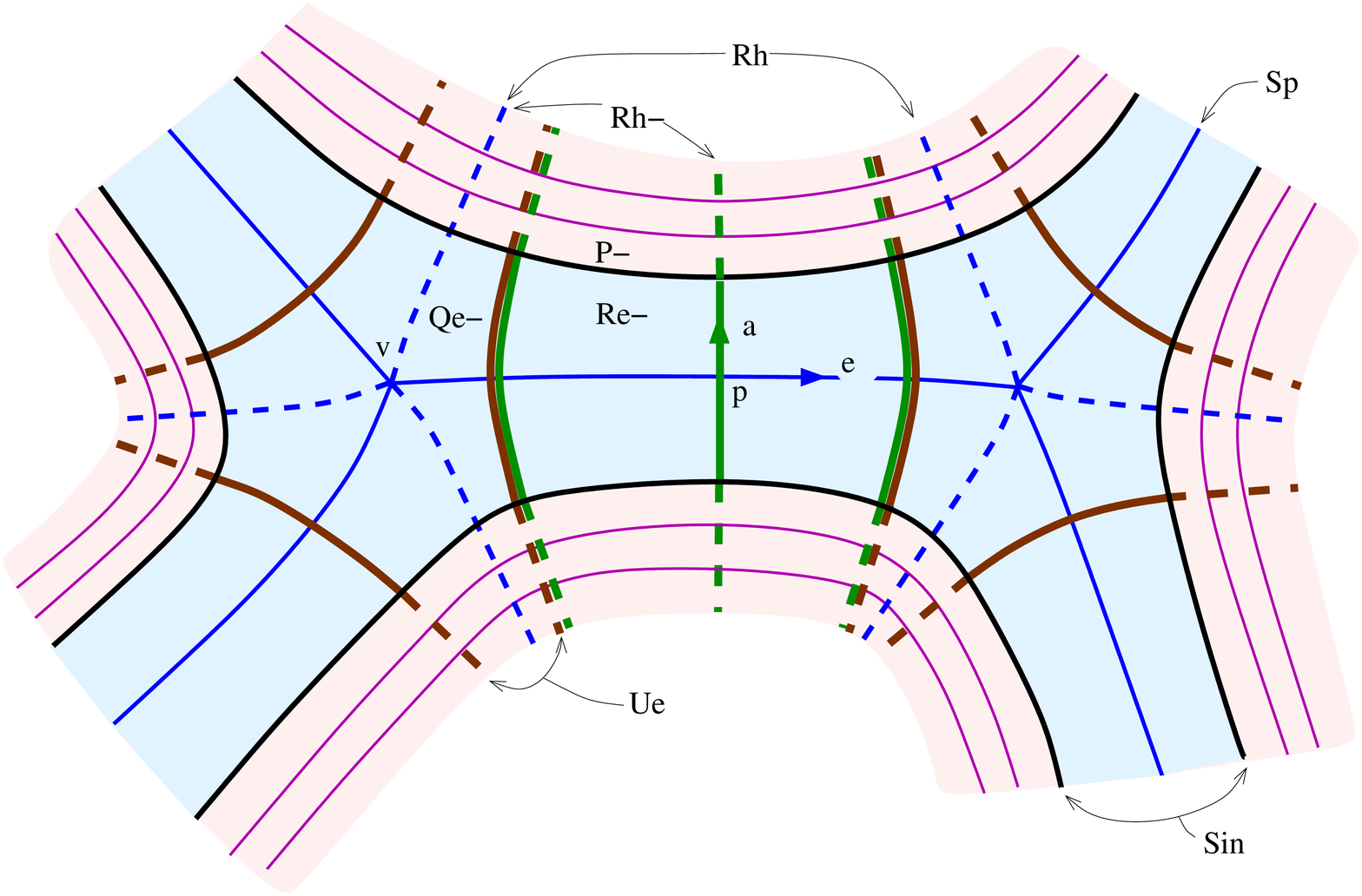}

\myCaption{Regions of $\gr8(\Si_m)$:
$\hat{U}^\e$ refers to $v$ and $P,Q,R,U$ to $\ora{\a}^*$.}
\label{fig:regions}
\end{figurehere}
\end{center}
\begin{itemize}
\item
Let $\ora{\a}$ be an oriented arc on $S$ with support $\a\in\bm{A}_{i_m}$.
Call $P(\ora{\a}^*,\Si_m)$ the projection of the geodesic edge
$f_m(\a)^*$ of $\mathrm{Sp}(\Si_m)$
to the boundary component of $\Si_m$ pointed by $f(\ora{\a})$
and orient $P(\ora{\a}^*,\Si_m)$ coherently with $\ora{\a}^*$
(and so reversing the orientation induced by $\Si_m$).

For every $b\in P(\ora{\a}^*,\Si_m)$, call $g_b$ the geodesic arc
that leaves $P(\ora{\a}^*,\Si_m)$ perpendicularly at $b$
and ends at $\a^*$ and define the quadrilateral
\[
R(\ora{\a}^*,\Si_m):=\bigcup_{b\in P(\ora{\a}^*,\Si_m)} g_b
\]
and let $\hat{R}(\ora{\a}^*,\Si_m)$ be the union of $R(\ora{\a}^*,\Si_m)$
and the flat rectangle of
$\gr8(\Si_m)$ of infinite height with basis $P(\ora{\a}^*,\Si_m)$.
\item
Assume now $\a\in \bm{A}\subset \bm{A}_{i_m}$ and
let $(\ora{\a},\ora{\b_1},\ora{\b_2})$ bound a hexagon
of $S\setminus \bm{A}_{i_m}$.
The formula
\[
\sinh(a/2)\sinh(w(\ora{\a}))=\frac{s_{\b_1}^2+s_{\b_2}^2-s_\a^2}
{2s_{\b_1} s_{\b_2}}
\]
shows that $w_{sp}(\ora{\a},f_m)>0$ for $m$ large enough,
because $a(f_m)=\ell_{\a}(f_m)\rar 0$. Thus, we can assume
that $w_{sp}(\ora{\a},f_m),w_{sp}(\ola{\a},f_m)>0$ for all $\a\in\bm{A}$
and all $m$.

Call $x$ the arc-length coordinate on $P(\ora{\a}^*,\Si_m)$ that
is zero at the projection of $s:=\a^*\cap\a$ and let
$P_-=P\cap\{x\leq 0\}$. Define
\[
R^\e(\ora{\a}^*,\Si_m):=\{ r_x\,|\,
x\in[-(1-\e)w_{sp}(\ora{\a},f_m), (1-\e)w_{sp}(\ola{\a},f_m)]
\}
\]
where $r_x$ is the hypercyclic arc parallel to $f(\a)$ that joins
$x\in P(\ora{\a}^*,\Si_m)$ and $f(\a)^*$, and let
$\hat{R}^\e(\ora{\a}^*,\Si_m)$ be the union of $R^\e(\ora{\a}^*,\Si_m)$
and the flat rectangle of $\gr8(\Si_m)$ that leans on it.

We can clearly put coordinates $(x,y)$
on $R^\e(\ora{\a}^*,\Si_m)\cup(\hat{R}(\ora{\a}^*,\Si_m)\setminus\Si_m)$
such that
\begin{itemize}
\item
$x$ extends the arc-length coordinate of $P$
\item
$(x,y)$ are orthonormal on the flat
part $\hat{R}(\ora{\a}^*,\Si_m)\setminus\Si_m$, which corresponds
to $[-w(\ora{\a},f_m),w(\ola{\a},f_m)]\times[0,\infty)$
\item
$(x,y)$ are orthogonal on the hyperbolic part
$R^\e(\ora{\a}^*,\Si_m)$, which corresponds to
$[-(1-\e)w(\ora{\a},f_m),(1-\e)w(\ola{\a},f_m)]\times[-a(f_m)/2,0]$;
moreover,
$\{x=const\}$ is a hypercycle parallel to $f(\a)$ and $\{y=const\}$
is a geodesic that crosses $f(\a)$ perpendicularly.
\end{itemize}
Finally, we set $R_-^\e:=R^\e\cap\{x\leq 0\}$ and
$\hat{R}_-^\e:=\hat{R}^\e\cap\{x\leq 0\}$, and we let $\hat{Q}_-^\e:=
\hat{R}_-\setminus\hat{R}_-^\e$.

Define analogously the regions for $\Si$, some of which will
depend on $m$. First, we call $\hat{R}(\ora{\a}^*,\Si,m):=T_{\ora{\a}^*}
\subset\Si$ and we put coordinates $\ti{x}=-\ti{w}_m(\ora{\a},f)+
\ti{w}(\a,f)u$ (which depends on $m$)
and $\ti{y}=\ti{w}(\a,f)v$ on it,
so that the Jenkins-Strebel differential $\varphi$ on $\Si$
restricts to $(d\tilde{x}+id\ti{y})^2$ on $\hat{R}(\ora{\a}^*,\Si,m)$.
Then, we define $\hat{R}_-(\ora{\a}^*,\Si,m):=\hat{R}(\ora{\a},\Si,m)
\cap\{\ti{x}\leq 0\}$ and $\hat{R}_-^\e(\ora{\a}^*,\Si,m):=
\hat{R}(\ora{\a},\Si,m)\cap\{-(1-\e)\ti{w}_m(\ora{\a}^*,f)
\leq \ti{x}\leq 0\}$
and finally $\hat{Q}_-^\e:=\hat{R}_-\setminus\hat{R}_-^\e$.
Define similarly the regions with $\ti{x}\geq 0$.
%
%
%
\item
If $v$ is a vertex of $G\subset\Si$, then let $f(\ora{\b_1})^*,\dots,
f(\ora{\b_j})^*$ be the (cyclically ordered) set of edges of $G$
outgoing from $v$, where $\b_h\in\bm{A}$ (the indices of the $\b$'s
are taken in $\Z/j\Z$). For every $m$ and $h$ there is an $l_h\geq 1$
such that $\ora{\b_h},\s_\infty^{-1}(\ora{\b_h}),
\s_{\infty}^{-2}(\ora{\b_h}),\dots,\s_\infty^{-l_h}(\ora{\b_h})=
\ola{\b_{h+1}}$ are distinct. Call
\[
\hat{U}^\e(\ora{\b_h},\Si_m):=
\hat{Q}^\e_-(\ora{\b_h},\Si_m)\cup \hat{Q}^\e_+(\ora{\b_{h+1}},\Si_m)
\cup\bigcup_{i=1}^{l_h-1}\hat{R}(\s_\infty^{-i}(\ora{\b_h}),\Si_m)
\]
and let
$\dis
\ti{w}(v,f_m)=\sum_h \sum_{i=1}^{l_h-1}
\ti{w}(\s_\infty^{-i}(\ora{\b_h}),\Si_m)
$
be the total (normalized) weight of the edges
$\{\eta_k\}$ of $G_m\subset\Si_m$
{\it that shrink to $v$}, that is
the edges supporting $\s_\infty^{-i}\ora{\b_h}$
with $h=1,\dots,j$ and $i=1,\dots,l_h-1$.

Set $U^\e:=\hat{U}^\e\cap\Si_m$ and,
similarly, $\hat{U}^\e(\ora{\b_h}^*,\Si,m):=
\hat{Q}_-^\e(\ora{\b_h}^*,\Si,m)\cup
\hat{Q}_+^\e(\ora{\b_{h+1}}^*,\Si,m)$
and
$\hat{U}^\e(v,\Si,m):=
\bigcup_{h=1}^j \hat{U}^\e(\ora{\b_h}^*,\Si,m)$.
\item
If $v$ is a nonmarked (smooth or singular) vertex of $G\subset\Si$, then 
we simply set
$\hat{U}^\e(v,\Si_m):=
\bigcup_{h=1}^j \hat{U}^\e(\ora{\b_h}^*,\Si_m)$.
\item
If $v$ is a smooth vertex of $\Si$ marked by $x_i$, then
we set $\hat{U}^\e(v,\Si_m):=\{x_i\}\cup\tilde{C}_i\cup
\bigcup_{h=1}^j \hat{U}^\e(\ora{\b_h}^*,\Si_m)$,
where $\tilde{C}_i$ is the flat cylinder corresponding
to $x_i$.
\end{itemize}
%
%
%
We choose $\dis\e_m=\mathrm{max}\{1/\bm{\Ll}(f_m),
1-\sum_{\a\in\bm{A}}\ti{w}(\a,f_m) \}^{1/2}$,
so that $\e_m\rar 0$, $\e_m \bm{\Ll}(f_m)\rar\infty$
and $(1-\sum_{\a\in\bm{A}} \tilde{w}(\a,f_m))/\e_m\rar 0$.
Moreover, we set $\d_m=\exp(-\e_m w(\a_0,f_m)/4)\rar 0$,
where $\a_0\in\bm{A}$
and $\tilde{w}(\a_0,f)=\mathrm{min}
\{\tilde{w}(\a,f)>0\,|\,\a\in\bm{A}\}$,
so that $a_i(f_m)<\d_m$ for $m$ large.\\

Define $F_m:\gr8(\Si_m)\rar\Si$ according
to the following prescriptions.

{\bf Edges.} For every $\a\in\bm{A}$ and every orientation $\ora{\a}$,
$F_m$ continuously
maps $\hat{R}_+^{\e_m}(\ora{\a}^*,\Si_m)$
onto $\hat{R}_+^{\e_m}(\ora{\a}^*,\Si,m)$ in such a way that
$\dis F_m(x,y)=\frac{2}{\bm{\Ll}(f_m)}\left(
\frac{\tilde{w}_m(\ora{\a},f)}{\tilde{w}(\ora{\a},f_m)}x,y\right)$
for $y\geq \d_m$ and the vertical arcs
$\{x\}\times[-a/2,\d_m]$ (whose length is $\d_m+a\cosh(x)/2$)
are homothetically
mapped to vertical trajectories
$\{\tilde{x}'\}\times\left[0,2\d_m/\bm{\Ll}(f_m)\right]$. 
Thus, the differential of $F_m$ (from the $xy$-coordinates on $\Si_m$
to the $\tilde{x}\tilde{y}$-coordinates on $\Si$) is
\[
dF_m=
\begin{cases}
\dis\frac{2}{\bm{\Ll}(f_m)}
\left(\begin{array}{cc}
\dis\frac{\ti{w}_m(\ora{\a},f)}{\ti{w}(\ora{\a},f_m)} & 0 \\
0 & 1
\end{array}\right)
& \text{if $y\geq\d_m$}\\
\dis\frac{2}{\bm{\Ll}(f_m)}
\left(\begin{array}{cc}
\dis\frac{\tilde{w}_m(\ora{\a},f)}{\tilde{w}(\ora{\a},f_m)} & 0 \\
\dis -\frac{ya\sinh(x)}{2\d_m\left(1+\frac{a}{2\d_m}\cosh(x)\right)^{2}}
& \dis\left(1+\frac{a}{2\d_m}\cosh(x)\right)^{-1}
\end{array}
\right)
& \text{if $0\leq y\leq\d_m$} \\
\dis\frac{2}{\bm{\Ll}(f_m)}
\left(\begin{array}{cc}
\dis \frac{\tilde{w}_m(\ora{\a},f)}{\tilde{w}(\ora{\a},f_m)} & 0 \\
\dis \frac{y\sinh(x)}{\left(1+\frac{a}{2\d_m}\cosh(x)\right)^{2}}
& \dis\frac{\cosh(x)}{1+\frac{a}{2\d_m}\cosh(x)}
\end{array}
\right)
& \text{if $y\leq 0$}
\end{cases}
\]
Because the metric $g_m$ on $\hat{R}_+^{\e_m}(\ora{\a}^*,\Si_m)$
in the $xy$-coordinates is
\[
g_m=
\begin{cases}
\left(\begin{array}{cc}
1 & 0 \\
0 & 1
\end{array}
\right)
& \text{if $y\geq 0$} \\
\left(\begin{array}{cc}
1 & 0 \\
0 & \cosh(x)^2
\end{array}
\right)
& \text{if $y\leq 0$}
\end{cases}
\]
we obtain $(F_m^{-1})^*(g_m)=M^t\, M$ (with respect to the
$\tilde{x}\tilde{y}$-coordinates), where
\[
M=\sqrt{g_m}dF_m^{-1}=
\begin{cases}
\dis\frac{\bm{\Ll}(f_m)}{2}
\left(
\begin{array}{cc}
\dis\frac{\ti{w}(\ora{\a},f_m)}{\ti{w}_m(\ora{\a},f)} & 0 \\
0 & 1
\end{array}\right)
\\
\dis\frac{\bm{\Ll}(f_m)}{2}
\left(\begin{array}{cc}
\dis\frac{\tilde{w}(\ora{\a},f_m)}{\tilde{w}_m(\ora{\a},f)} & 0 \\
\dis \frac{ya\sinh(x)\ti{w}(\ora{\a},f_m)}
{2\d_m\ti{w}_m(\ora{\a},f)(1+\frac{a}{2\d_m}\cosh(x))}
& 1+\frac{a}{2\d_m}\cosh(x)
\end{array}
\right)
\\
\dis\frac{\bm{\Ll}(f_m)}{2}
\left(\begin{array}{cc}
\dis \frac{\tilde{w}(\ora{\a},f_m)}{\tilde{w}_m(\ora{\a},f)} & 0 \\
\dis -\frac{y\sinh(x)\ti{w}(\ora{\a},f_m)}
{\ti{w}_m(\ora{\a},f)(1+\frac{a}{2\d_m}\cosh(x))}
& 1+\frac{a}{2\d_m}\cosh(x)
\end{array}
\right)
\end{cases}
\]
in the three different regions.

If $w(\ora{\a},f_m)\geq w(\a,f_m)/2$, then
$\dis\frac{a}{2}\sinh(w(\ora{\a},f_m))\leq 1$ implies
$\dis\frac{a}{2}\sinh(x)\leq\frac{a}{2}\sinh[(1-\e_m)w(\ora{\a},f_m)]
\lessapprox \exp(-\e_m w(\ora{\a},f_m))$
and so $\dis\frac{ya\sinh(x)}{2\d_m}\lessapprox \exp(-\e_m w(\a_0,f_m)/2)$
and $\dis\frac{a\cosh(x)}{2\d_m}\lessapprox \exp(-\e_m w(\a_0,f_m)/4)$.
Hence, on the region where we have defined $F_m$,
the distortion is bounded by
\[
\mathrm{max}
\left\{
\frac{\tilde{w}_m(\ora{\a},f)}{\tilde{w}(\ora{\a},f_m)},\,
\frac{\tilde{w}(\ora{\a},f_m)}{\tilde{w}_m(\ora{\a},f)}\,\Big|\,
\a\in\bm{A}
\right\}
\cdot
\left(
1+\exp[-\e_m w(\a_0,f_m)/4]
\right)\rar 1\, .
\]

{\bf Around the vertices.}
For every vertex $v\in G\subset\Si$ with outgoing edges
$\{\ora{\b_h}^*\}$, we require
$F_m$ to map $\hat{U}^{\e_m}(v,\Si_m)\cap\{y\geq\d_m\}$
onto $\hat{U}^{\e_m}(v,\Si,m)\cap\{\tilde{y}\geq 2\d_m/\bm{\Ll}(f_m)\}$
with differential (from $(x,y)$ to $(\tilde{x},\tilde{y})$)
constantly equal to
\[
\frac{2}{\bm{\Ll}(f_m)}
\left(\begin{array}{cc}
c & 0 \\
0 & 1
\end{array}\right),\
\text{where}\
c=\frac{\dis \e_m \sum_h \tilde{w}_m(\ora{\b_h},f)}
{\dis\e_m\sum_h \tilde{w}(\ora{\b_h},f_m)+\tilde{w}(v,f_m)}
\]
Notice that
\[
c-1=\frac{\dis \sum_h \left(\tilde{w}_m(\ora{\b_h},f)-
\tilde{w}(\ora{\b_h},f_m)\right)
-\e_m^{-1}\ti{w}(v,f_m)}
{\dis\sum_h \tilde{w}(\ora{\b_h},f_m)+\e_m^{-1}\tilde{w}(v,f_m)
}\rar 0
\]
because $\dis \e_m^{-1}\tilde{w}(v,f_m)\leq\e_m^{-1}
(1-\sum_{\a\in\bm{A}}\tilde{w}(\a,f_m))\rar 0$.
Hence, the distortion of $F_m$ goes to $1$.

{\bf Neighbourhoods of the vertices.}
If $v\in G\subset\Si$ is smooth, then define $F_m$ to be
a diffeomorphism between
$\hat{U}^{\e_m}(v,\Si_m)\setminus\{y\geq\d_m\}$
and $\hat{U}^{\e_m}(v,\Si,m)
\setminus\{\tilde{y}\geq 2\d_m/\bm{\Ll}(f_m)\}$.
If $v$ is also marked, then we can require $F_m$ to preserve
the marking.

If $v$ is a node between two visible components, then
$F_m$ maps $\hat{U}^{\e_m}(v,\Si_m)\setminus\{y\geq \d_m\}$
onto $\hat{U}^{\e_m}(v,\Si,m)\setminus\{\tilde{y}\geq 2\d_m/\bm{\Ll}(f_m)\}$
shrinking the edges $\{\eta_i^*\}$ to $v$ and as a diffeomorphism
elsewhere.

If $\Si'\subset\Si$ is an invisible component and
$v_1,\dots,v_l$ are vertices of $G\subset\Si_+$ and nodes of $\Si'$, then
let $\{\eta_i\}$ be the sub-arc-system $\bm{A}_{i_m}\cap f^{-1}(\Si')$.
We require $F_m$ to map
\[
\left(
\bigcup_i \left(\hat{R}^{\e_m}(\ora{\eta_i}^*,\Si_m)
\cup \hat{R}^{\e_m}(\ola{\eta_i}^*,\Si_m)\right)
\cup \bigcup_h \hat{U}^{\e_m}(v_h,\Si_m)
\right)
\setminus\{y\geq\d_m\}
\]
onto
$\dis\Si'\cup \left(\bigcup_h \hat{U}^{\e_m}(v_h,\Si,m)\setminus
\{\tilde{y}\geq 2\d_m/\bm{\Ll}(f_m)\}\right)$
by shrinking $\hat{U}^{\e_m}(v_h,\Si_m)\cap\{y=\d_m/2\}$ to $v_h$
and as a diffeomorphism elsewhere.
\end{subsubsection}
%
%
\end{subsection}
%
%
\begin{subsection}{Bijectivity of $(\gr8,\Ll)$}
The bijectivity at infinity (namely, for $\bm{\Ll}=\infty$) follows
from Theorem~\ref{thm:strebel}. Thus, let's select a (possibly
empty) system of curves $\bm{\g}=\{\g_1,\dots,\g_k\}$ on $S$ and let's
consider the stratum $\mathcal{S}(\bm{\g})\subset\That(S)$
in which $\bm{\Ll}<\infty$ and $\ell_{\g_i}=0$.

To show that $(\gr8,\Ll)$ gives a bijection of $\mathcal{S}(\bm{\g})$
onto its image, it is sufficient to work separately on each component
of $S\setminus\bm{\g}$. Thus, we can reduce to the case in which
$\g_i=C_i\subset\pa S$ and $\gr8$ glues a cylinder
at the boundary components $C_{k+1},\dots,C_n$.
Hence, we are reduced to show that the grafting map
\[
\gr8':\Teich(S)(\up)\lra \Teich(S)(0)
\]
is bijective for every $p_{k+1},\dots,p_n\in\R_+$
, where $p_1=\dots=p_k=0$.
We already know that $\gr8'$ is continuous and proper:
we will show that it is a local homeomorphism by adapting
the argument of \cite{scannell-wolf:grafting}.
Here we describe what considerations are needed to make their
proof work in our case.
\begin{remark}
Here we are using the notation $\Teich(S)(0)$ instead
of $\Teich(R,x)$ because we want to stress that
we are regarding $\gr8(\Si)$ as a hyperbolic surface,
with the metric coming from the uniformization.
\end{remark}
The grafted metrics are $C^{1,1}$ but the map $\gr8'$ is real-analytic.
In fact, given a real-analytic arc $[f_t:S\rar\Si_t]$ in $\Teich(S)(\up)$
and chosen representatives $f_t$ so that $f_0\circ f_t^{-1}:\Si_t\rar\Si_0$
is an isometry on the boundary components $C_{k+1,t},\dots,C_{n,t}$ and
harmonic in the interior with respect to the hyperbolic
metrics (so that the hyperbolic
metrics pull back to a real-analytic family $\s_t$ on $S$),
we can choose
the grafted maps $\gr8'(f_t):S\rar\Si_t$
so that $\gr8'(f_0)\circ\gr8'(f_t)^{-1}$
extend $f_0\circ f_t^{-1}$ as isometries on the cylinders
$\tilde{C}_{i,t}:=C_{i,t}\times[0,\infty)$.
Hence, the family of metrics $\gr8'(\s_t)$ on $S$, obtained
by pulling the Thurston metric back thourgh $\gr8'(f_t)$,
is real-analytic in $t$ and so the arc $[\gr8'(f_t)]$ in $\Teich(S)(0)$
is real-analytic.

Thus, it is sufficient to show that the differential
$d\gr8'$ is injective at every point of $\Teich(S)(\up)$.\\

Given a real-analytic one-parameter family $f_t:S\rar\Si_t$
corresponding to a tangent vector $v\in T_{[f_0]}\Teich(S)(\up)$,
assume that the grafted family $[\gr8'(f_t):S\rar\gr8'(\Si_t)]$ 
defined above
determines the zero tangent vector in $T_{[\gr8'(f_0)]}\Teich(S)(0)$.

Call $\gr8'(\sigma_t)$ the pull-back via $f_t$ of the hyperbolic
metric of $\Si_t$
and construct the harmonic representative
$F_t:(S,\gr8'(\s_t))\rar(S,\gr8'(\s_0))$
in the class of the identity as follows.

Give orthonormal coordinates $(x,y)$ to the
cylinder $\tilde{C}_{i,t}\cong C_{i,t}\times[0,\infty)$
that is glued at the boundary component
$C_{i,t}\subset\Si_t$ for $i=k+1,\dots,n$,
in such a way that $x$ is the arc-length parameter of the circumferences and
$y\in[0,\infty)$.
\begin{remark}\label{rem:coordinates}
The $(x,y)$ coordinates can be extended to
an orthogonal system in a small
hyperbolic collar of $C_{i,t}$ in such a way that $y$ is the arc-length
parameter along the geodesics $\{x=const\}$.
Thus, for $y\in(-\e,0)$, the metric 
looks like $\cosh(y)^2 dx^2+dy^2=dx^2+dy^2+O(\e^2)$.
\end{remark}
Call {\it $M$-ends} of $\gr8'(\Si_t)$ the subcylinders
$\tilde{C}_{i,t}\times[M,\infty)$ for
$i=k+1,\dots,n$ and use the same terminology
for their images in $S$ via $\gr8'(f_t)^{-1}$.

For every $t$, let $\mathfrak{F}_t:=\bigcup_{M\geq 0}\mathfrak{F}_t(M)$
where $\mathfrak{F}_t(M)$ is the set of $C^{1,1}$ diffeomorphisms
$g_t:(S,\gr8(\s_t))\rar(S,\gr8(\s_0))$ homotopic to the identity,
such that $g_t$ isometrically preserves the $M$-ends.
Clearly, $\mathfrak{F}_t(M)\subseteq\mathfrak{F}_t(M')$, if $M\leq M'$.

Let $e(g_t)=\frac{1}{2}\|\nabla g_t\|^2$ be the energy density of $g_t$,
$\dis \Hh(g_t)=\|dg_t(\pa_z)\|^2\,\frac{dz\,d\bar{z}}{\gr8'(\s_t)}$,
where $z$ is a local conformal coordinate on $(S,\gr8'(\s_t))$,
and $\Jac(g_t)$ the Jacobian determinant of $g_t$,
so that $e(g_t)=2\Hh(g_t)-\Jac(g_t)$.
Notice that, if $g_t$ is an oriented
diffeomorphism, then $0<\Jac(g_t)\leq\Hh(g_t)\leq e(g_t)$
at each point.

Define also the reduced quantities $\tilde{e}(g_t)=e(g_t)-1$, $\tilde{\Hh}(g_t)
=\Hh(g_t)-1$ and $\tilde{\Jac}(g_t)=\Jac(g_t)-1$, so that the reduced energy
\[
\tilde{E}(g_t):=\int_{S}\tilde{e}(g_t)\gr8'(\s_t)
\]
is well-defined for every $g_t\in\mathfrak{F}_t$.
For instance, the identity map on $S$ belongs to $\mathfrak{F}_t(0)$
and its reduced energy is $E(f_0\circ f_t^{-1})-2\pi\chi(S)$.

As $\gr8'(\s_0)$ is nonpositively curved, the map $F_{t,M}$
of least energy in
$\mathfrak{F}_t(M)$ is harmonic away from the $M$-ends and so is an
oriented diffeomorphism. Thus,
\[
0=\int_{S}\tilde{\Jac}(F_{t,M})\gr8'(\s_t)\leq
\int_{S}\tilde{\Hh}(F_{t,M})\gr8'(\s_t)\leq
\tilde{E}(F_{t,M})
\]
Thus, the map $F_t$ of least (reduced) energy in $\mathfrak{F}_t$ can be
obtained as a limit of the $F_{t,M}$'s and it is clearly unique.
Call $\tilde{\Hh}_t:=\tilde{\Hh}(F_t)$ and
similarly $\tilde{e}_t=\tilde{e}(F_t)$.

Following Scannell-Wolf (but noticing that the roles of $x$ and $y$
here are exchanged compared to their paper), one can show that
\begin{itemize}
\item
the family $\{F_t\}$ is real-analytic in $t$
\item
for every small $t$, the map $F_t$ is (locally) $C^{2,\a}$ on $S$;
so is the vector field $\dot{F}:=\dot{F}_0$ (hence, the analyticity
of $F_t$ implies that $\tilde{\Hh}_t$ and $\tilde{e}_t$
are real-analytic in $t$ too)
\item
the function $\dot{\tilde{\Hh}}:=\dot{\tilde{\Hh}}_0$
is locally Lipschitz and it is
harmonic on the flat cylinders
\item
along every $C_{k+1},\dots,C_n$, we have
\[
V=-\frac{1}{2}
\left(
(\pa_y\dot{\tilde{\Hh}})_+
-(\pa_y\dot{\tilde{\Hh}})_-
\right)
\]
where $V(x,y)$ is a harmonic function defined on the
cylinders $\tilde{C}_{i,0}$
(and on first-order thickenings of $C_{i,0}$) that can be identified
to the $y$-component of $\dot{F}$ and
$w(x,0)_+$ simply means $\dis\lim_{y\rar 0^+}w(x,y)$.
\item
$\dis V_y=\frac{1}{2}\dot{\Hh}+c_i$ on each
$\tilde{C}_i$, where $c_i$ is a constant that may depend on
the cylinder.
\end{itemize}

Now on, let all line integrals be with respect to the arc-length
parameter $dx$ and all surface integrals with respect to $\gr8'(\s_0)$.
Notice that
\[
\int_S \dot{\Hh}=
\int_S \dot{\tilde{\Hh}}=
\lim_{t\rar 0} \frac{1}{t}
\int_S \tilde{\Hh}_t
\]
because $\ti{\Hh}_0=0$.
As the integral on the right is a real-analytic function of $t$
which vanishes at $t=0$, we conclude that $\dot{\ti{\Hh}}$
is integrable.
By the same argument, so is $\dot{\tilde{e}}$.

On the other hand,
$\dis\dot{\tilde{e}}=
\frac{1}{2}\frac{\pa}{\pa t}\|\nabla F_t\|^2 \Big|_{t=0}
\geq |V_y|$
and so $V_y$ is integrable too and all constants $c_i=0$.
Thus, $V$ and $\dot{\tilde{\Hh}}$ decay at least
as $\exp(-2\pi y/p_i)$ on $\tilde{C}_{i,0}$ and we can write
\[
0=\int_{\tilde{C}_{i,0}}V\Delta V =
-\int_{\tilde{C}_{i,0}}\|\nabla V\|^2
+\int_{C_{i,0}}V\pa_n V
\]
Moreover,
\begin{align}\label{eq:first}
0\leq\int_{\tilde{C}_{i,0}}\|\nabla V\|^2=
\int_{C_{i,0}}V_y V=
\frac{1}{2}\int_{C_{i,0}}\dot{\tilde{\Hh}}V
\end{align}
On the other hand, multiplying by $\dot{\tilde{\Hh}}=
\dot{\Hh}$ and integrating by parts the
linearized equation 
\[
(\Delta_{\gr8'(\s_0)}+2K_0)\dot{\Hh}=0
\]
where $K_0$ is the curvature of $\gr8'(\s_0)$, we obtain
\[
\begin{cases}
\dis 0\leq\int_{\tilde{C}_{i,0}}\|\nabla\dot{\Hh}\|^2=
\int_{C_{i,0}}\dot{\Hh}(\pa_n\dot{\Hh})_+\\
\dis 0\leq\int_{S_{hyp}}\|\nabla\dot{\Hh}\|^2+
2|\dot{\Hh}|^2=
-\sum_{i=k+1}^n \int_{-C_{i,0}}
\dot{\Hh}(\pa_n\dot{\Hh})_-
\end{cases}
\]
where $S_{hyp}$ is the $\gr8'(\s_0)$-hyperbolic part of $S$.

From $\dis 0\leq \sum_{i=k+1}^n \int_{C_{i,0}}
\dot{\Hh}\left(
(\pa_y\dot{\Hh})_-
-(\pa_y\dot{\Hh})_+
\right)$, we
finally get
\begin{equation}\label{eq:second}
0\leq 2\int_S \|\nabla\dot{\Hh}\|^2-
2K\|\dot{\Hh}\|^2=
-\sum_{i=k+1}^n \int_{C_{i,0}} \dot{\Hh}V
\end{equation}
Combining Equation~\ref{eq:first} and Equation~\ref{eq:second},
we obtain
\[
\int_{C_{i,0}}\dot{\Hh}V=0\qquad
\forall\ i=k+1,\dots,n
\]
and so $\dot{\Hh}=0$ on $S$.

Hence, $F_t$ is a $(1+o(t))$-isometry between $\gr8'(\s_t)$ and $\gr8'(\s_0)$
and one can easily conclude that $\s_t$ and $\s_0$ are
$(1+o(t))$-isometric too.
\end{subsection}
%
%
%
\begin{subsection}{More on infinitely grafted structures}
\begin{subsubsection}{Projective structures.}
Consider a compact Riemann surface $R$ without boundary and
of genus at least $2$. A {\it projective structure}
on a marked surface $[f:R\rar R']$ is an equivalence class of
holomorphic atlases $\mathfrak{U}=\{f_i:U_i\rar\C\Pr^1\,|\,R'\supset U_i
\ \text{open} \}$
for $R'$ such that the transition
functions belong to
$\mathrm{Aut}(\C\Pr^1)\cong\mathrm{PSL}(2,\C)$, that is
$f_i\Big|_{U_i\cap U_j}$ and $f_j\Big|_{U_i\cap U_j}$ are
{\it projectively equivalent}.

Given two projective structures, represented by
maximal atlases $\mathfrak{U}$ and $\mathfrak{V}$, on the same
$[f:R\rar R']\in\Teich(R)$
and a point $p\in R'$, we want to measure how charts of $\mathfrak{U}$
are not projectively equivalent to charts in $\mathfrak{V}$ around $p$.
So, let $f:U\rar\C\Pr^1$ be a chart in $\mathfrak{U}$ and
$g:U\rar\C\Pr^1$ a chart in $\mathfrak{V}$, with $U\subset R'$.
There exists a unique $\sigma\in\mathrm{PSL}(2,\C)$ such that
$f$ and $\s\circ g$ agree up to second order at $p$.
Then, $(f-\s\circ g)''':T_p U\rar T_{f(p)}\C\Pr^1$
is a homogeneous cubic map and
$f'(p)^{-1}\circ (f-\s\circ g)'''$ is a homogeneous cubic
endomorphism of $T_p U$, and so an element $\Sch(f,g)(p)$
of $(T^*_p U)^{\otimes 2}$. The holomorphic quadratic differential
$\Sch(\mathfrak{U},\mathfrak{V})$ on $R'$
is called {\it Schwarzian derivative}.
It is known that, given a $\mathfrak{U}$ and a holomorphic
quadratic differential $\varphi\in\Qq_{R'}$,
there exists a unique projective structure $\mathfrak{V}$ on $R'$
such that $\Sch(\mathfrak{U},\mathfrak{V})=\varphi$.

Thus,
the natural projection $\pi:\Pj(R)\rar\Teich(R)$
from the set $\Pj(R)$ of projective structures on $R$
(up to isotopy) to the Teichm\"uller space of $R$ is 
a principal $\Qq$-bundle, where
$\Qq\rar\Teich(R)$ is the bundle of holomorphic
quadratic differentials.

On the other hand, the grafting map $\mathrm{gr}:\Teich(R)\times\ML(R)
\rar\Teich(R)$ admits a lifting
\[
\mathrm{Gr}:\Teich(R)\times\ML(R)\arr{\sim}{\lra}\Pj(R)
\]
which is a homeomorphism (Thurston) and such that $\mathrm{Gr}(-,0)$
corresponds to the Poincar\'e structure.
We recall that a surface with projective
structure comes endowed with a Thurston $C^{1,1}$ metric:
in particular, if $\l=c_1\g_1+\dots+c_n\g_n$ is a multi-curve
on $R$, then $\mathrm{Gr}(R',\l)$ is made of a hyperbolic
piece, isometric to $R'\setminus\mathrm{supp}(\l)$ and
$n$ flat cylinders $F_1,\dots,F_n$, with $F_i$ homotopic to $\g_i$
and of height $c_i$.

It is a general fact that $\mathrm{Gr}(-,\l)$
is a real-analytic section of $\pi$ for all $\l\in\ML$.
\end{subsubsection}
\begin{subsubsection}{A compactification of $\Pj(R)$.}
The homeomorphism $\Teich(R)\times\ML(R)\cong\Pj(R)$
shows that
sequences $([f_m:R\rar R'_m],\l_m)$
in $\Pj(R)$ can diverge in two ``directions''.

Dumas \cite{dumas:grafting} provides a {\it grafting
compactification}
of $\Pj(R)$ by
separately compactifying $\Teich(R)$
and $\ML(R)$. In particular, he defines $\ol{\Pj}(R)
:=\Tbar^{Th}(R)\times\ol{\ML}(R)$, where $\Tbar^{Th}(R)=
\Teich(R)\cup\Pr\ML(R)$ is Thurston's compactification and
$\ol{\ML}(R)=\ML(R)\cup\Pr\ML(R)$ is the natural projective
compactification of $\ML(R)$.
In particular, the locus $\Tbar^{Th}(R)\times\Pr\ML(R)$ corresponds
to ``infinitely grafted surfaces''.

In order to describe the asymptotic properties of $\ol{\Pj}(R)$,
we recall the following well-known result.

\begin{theorem}[\cite{hubbard-masur:foliations}]
The map
\[
\Lambda:\Qq\rar \Teich(R)\times\ML(R)
\]
defined as $([f:R\rar R'],\varphi)\mapsto
([f:R\rar R'],f^*\Lambda_{R'}(\varphi))$
is a homeomorphism, where $\Lambda_{R'}(\varphi)$ is the
measured lamination on $R'$
obtained by straightening the (measured) horizontal
foliation of $\varphi$.
\end{theorem}

The {\it antipodal map} is the homeomorphism
$\i:\Teich(R)\times\ML(R)\rar\Teich(R)\times\ML(R)$
given by $\i([f:R\rar R'],\l)=\Lambda(-\Lambda^{-1}([f],\l))$.
The following result
shows that the restriction
$\i_{f}:\ML(R)\rar\ML(R)$
of $i$ to a certain point $[f]\in\Teich(R)$
controls the asymptotic behavior of $\pi^{-1}(f)$.

\begin{theorem}[\cite{dumas:grafting},
\cite{dumas:schwarzian}]\label{thm:dumas}
Let $\{([f_m:R\rar R_m],\l_m)\}\subset\Teich(R)\times\ML(R)$
be a diverging sequence such that
$\pi\circ\mathrm{Gr}_{\l_m}(f_m)=[f:R\rar R']$.
The following are equivalent:
\begin{enumerate}
\item
$\l_m\rar [\l]$ in $\ol{\ML}(R)$, where $[\l]\in\Pr\ML(R)$
\item
$[f_m]\rar [\i_{f}(\l)]$ in $\Tbar^{Th}(R)$
\item
$\Lambda_{f}(-\Sch(\mathrm{Gr}_{\l_m}(f_m)))\rar [\l]$
in $\ol{\ML}(R)$
\item
$\Lambda_{f}
(\Sch
(\mathrm{Gr}_{\l_m}(f_m)))
\rar
[i_{f}(\l)]$
in $\ol{\ML}(R)$,
where the Schwarzian derivative is considered with respect
to the Poincar\'e structure.
\end{enumerate}
When this happens, we also have $[\Ho(\kappa_m)]\rar
[\Lambda^{-1}(R',\l)]$ in $\Pr L^1(R',K^{\otimes 2})$, where
$\Ho(\kappa_m)$ is the Hopf differential of the
collapsing map $\kappa_m:R'\rar R_m$.
\end{theorem}
We recall that, if $\l_m$ is a multi-curve $c_1\g_1+\dots+c_n\g_n$,
then $\kappa_m$ collapses the $n$ grafted cylinders onto the respective
geodesics and is the identity elsewhere.
Thus, if the $j$-th flat cylinder
is isometric to $[0,\ell_j]\times[0,c_j]/ (0,y)\sim(\ell_j,y)$, then
$\Ho(\kappa_m)$ restricts to $dz^2$ on the grafted cylinders
and vanishes on the remaining hyperbolic portion of $R'$.
\begin{remark}
The theorem implies that the boundary
of $\pi^{-1}(f)\subset\ol{\Pj}(R)$ is
exactly the graph of the projectivization of $i_{f}$.
\end{remark}
\end{subsubsection}
\begin{subsubsection}{Surfaces with infinitely grafted ends.}
We can adapt Theorem~\ref{thm:grafting} to our situation,
when we restrict our attention to smooth hyperbolic surfaces with
large boundary.

Let $S$ be a compact oriented surface of genus $g$ with boundary components
$C_1,\dots,C_n$ (and $\chi(S)=2-2g-n<0$) and let $dS$ be its double.
%
%
\begin{theorem}\label{thm:mydumas}
Let $\{f_m:S\rar \Si_m\}\subset\Teich(S)$
be a sequence such that
$(\mathrm{gr}_\infty,\Ll)(f_m)=([f:(R,x)\rar (R',x')],\up_m)\in
\Teich(R,x)
\times\R_+^n$.
The following are equivalent:
\begin{enumerate}
\item
$\up_m\rar(\up,\infty)$ in $\Delta^{n-1}\times(0,\infty]$
\item
$[f_m]\rar \tilde{w}$ in $\Tbar^a(S)$,
where $\tilde{w}$ is the projective multi-arc associated to the
vertical foliation of the Jenkins-Strebel differential $\varphi_{JS}$
on $(R',x')$ with weights $\up$ (see Theorem~\ref{thm:strebel}).
\end{enumerate}
When this happens, we also have
\begin{itemize}
\item[(a)]
$4\bm{p_m}\!\!\!^{-2}\Ho(\kappa_m)\rar
\varphi_{JS}$ in $L^1_{loc}(R',K(x')^{\otimes 2})$, where
$\Ho(\kappa_m)$ is the Hopf differential of the
collapsing map $\kappa_m:R'\rar \Si_m$
\item[(b)]
with respect to the Poincar\'e projective structure,
$2\bm{p_m}\!\!\!^{-2}\Sch(\mathrm{Gr}_\infty(f_m))\rar -\varphi_{JS}$
in $H^0(R',K(x')^{\otimes 2})$.
\end{itemize}
\end{theorem}
\begin{remark}
We have denoted by $\mathrm{Gr}_\infty(f_m:S\rar\Si_m)$
the ($S$-marked)
surface with projective structure obtained from $\Si_m$
by grafting cylinders of infinite length at its ends.
This is a somewhat ``very exotic''
projective structure, whose developing map
wraps infinitely many times around $\mathbb{CP}^1$.
Its Schwarzian with respect to the Poincar\'e structure
has double poles at the cusps.
\end{remark}

We have already shown that (1) and (2) are equivalent
to each other.

\begin{lemma}
Given an increasing sequence $\{m_k\}\subset\N$ and a diverging
sequence $\{t_k\}\subset\R_+$,
we have $\bigcup_k R'_k=\dot{R}'$, where we consider
$R'_k:=\mathrm{gr}_{t_k\bm{p_{m_k}}\pa\Si_{m_k}}(\Si_{m_k})$
as embedded inside $\dot{R}'$.
\end{lemma}
\begin{proof}
Let $z\in R'\setminus\bigcup_k R'_k$ and notice that,
for each $k$, $R'\setminus R'_k$ is a disjoint union of $n$ discs.
Up to extracting a subsequence, we can assume that $z$ belongs
to the $j$-th disc, together with $x'_j$.
The image of $C_j\subset S$ inside $R'$ separates
$\{x_j,z\}$ from the rest of the surface.
Because $t_k\rar \infty$, the extremal length
$Ext_{C_j}(R'_k)\rar 0$ as $k\rar\infty$.
This implies $z=x_j$.
\end{proof}

The proof of (a) follows
\cite{dumas:grafting} (see also \cite{dumas:grafting2})
with minor modifications:
\begin{itemize}
\item
because of the previous lemma,
for every compact $K\subset \dot{R}'$,
there exist $t_0>0$ such that
$K\subset K_m^t:=\mathrm{gr}_{t\bm{p_m}\pa\Si_m}(\Si_m)\subset
\dot{R}'$ for
every $t\geq t_0$
\item
let $h_m:\dot{R}'\rar\Si_m$ be the harmonic map homotopic to
$\kappa_m$, that is the limit as $s\rar\infty$ of
the harmonic maps $h_m^s:\mathrm{gr}_{s\pa\Si_m}(\Si_m)\rar\Si_m$
that restrict to isometries at the boundary: we clearly
have
\[
\|\Ho(h_m)-\Ho(\kappa_m)\|_{L^1(K)}\leq
\|\Ho(h_m)-\Ho(\kappa_m)\|_{L^1(K_m^t)}
\]
and $E_{K_m^t}(h_m)<E_{K_m^t}(\kappa_m)=2\pi|\chi(S)|+
t\bm{p_m}\!\!\!^2/2
\leq E_{K_m^t}(h_m)+2\pi|\chi(S)|$, where $E_{K_m^t}$ is the
integral of the energy density on $K_m^t$
\item
the statement that $[\Ho(h_m)]\rar[\varphi]$ as $m\rar\infty$
is basically proven by Wolf in \cite{wolf:harmonic}; in fact,
the considerations involved in his argument do not require
the integrability of $\Ho(h_m)$ or $\varphi$ over the whole $\dot{R}'$:
rescaling the Hopf differential in order to have the right
boundary lengths, one obtains
\[
4\bm{p_m}\!\!\!^{-2}\Ho(h_m)\rar\varphi
\qquad \text{in $L^1_{loc}(\dot{R}')$}
\]
\item
the local estimate
\[
\|\Ho(h_m)-\Ho(\kappa_m) \|_{L^1(K)}
\leq \sqrt{2(E_K(h_m)-E_K(\kappa_m))}
\left(\sqrt{E_K(h_m)}+\sqrt{E_K(\kappa_m)} \right)
\]
is obtained in the proof of Proposition~2.6.3
of \cite{korevaar-schoen:sobolev}
\item
one easily concludes, because
$\|\Ho(h_m)-\Ho(\kappa_m) \|_{L^1(K_m^t)}
=O(\bm{p_m}\sqrt{t})$ and $\|\Ho(h_m)\|_{L^1(K_m^t)}=
O(\bm{p_m}\!\!\!^2 t)$.
\end{itemize}

Assertion (b) is also basically proven in \cite{dumas:schwarzian}
up to minor considerations.
\begin{itemize}
\item
Call $\rho$ the hyperbolic metric on $\dot{R}'$ and $\rho_m$
the Thurston metric on $\mathrm{Gr}_\infty(\Si_m)\cong\dot{R}'$;
moreover, let $\beta_m$ be the Schwarzian tensor
$\beta(\rho,\rho_m)=[\mathrm{Hess}_{\rho}(\sigma_m)-d\s_m\otimes d\s_m]^{2,0}$,
where $\s_m=\s(\rho,\rho_m)=\log(\rho_m/\rho)$.
\item
The decomposition (\cite{dumas:schwarzian}, Theorem~7.1)
\[
\bm{S}(\mathrm{Gr}_\infty(\Si_m))=2\beta_m-2\Ho(\kappa_m)
\]
(where $\bm{S}$ is with respect to the Poincar\'e structure on $\dot{R}'$)
still holds, because it relies on local considerations.
\item
Let $K$ be the compact subsurface of $\dot{R}'$
obtained by removing all $n$ horoballs of circumference
$1/4$ at $x'$.
Moreover, let $\rho_{\up}$ be the Thurston metric on $\Si$
obtained by grafting infinite
flat cylinders at the boundary of
$(\mathrm{gr}_\infty,\Ll)^{-1}([f],\up)$ and call
$\hat{\rho}_{\up}:=(1+\bm{p}^2)\rho_{\up}$ the normalized
metrics.
The set $\mathcal{N}=\{\hat{\rho}_{\up}\,|\,
\up\in\Delta^{n-1}\times[0,\infty]\}$
is compact in $L^\infty(K)$.
Thus, $\|\hat{\rho}_{\up}/\rho\|_{L^\infty(K)}<c$
and all restrictions to $K$ of metrics in $\mathcal{N}$
are pairwise H\"older equivalent with factor and exponent
dependent on $\dot{R}'$ only (same proof as in Theorem 9.2
of \cite{dumas:schwarzian}).
\item
The same estimates of \cite{dumas:schwarzian} give
(Theorem~11.4)
\[
\|\beta_m\|_{L^1(D_{\d/4},\rho)}\leq c
\]
where $c$ depends on $\dot{R}'$ and $\d$.
\item
All norms are equivalent on $H^0(R',K(x')^{\otimes 2})$,
so we consider the $L^1$ norm on $K\subset\dot{R}'$
and we observe that $\|\psi\|_{L^1(D_{\d/4},\rho)}\leq c'
\|\psi\|_{L^1(K)}$
for any $\rho$-ball of radius $\d/4$ embedded in $K$.
\item
There exists $t_0$ (dependent only on $\dot{R}'$) such that
$K\subset K_m^t$ for all $m$.
Thus,
\[
\begin{array}{l}
\qquad
\|2\bm{S}(\mathrm{Gr}_\infty(\Si_m))+
\varphi_{JS}\|_{L^1(K)}\leq \\
\leq c_1\|2\bm{S}(\mathrm{Gr}_\infty(\Si_m))+
4\Ho(\kappa_m)\|_{L^1(D_{\d/4},\rho)}
+\|4\Ho(\kappa_m)-\bm{p_m}\!\!\!^2\varphi_{JS}\|_{L^1(K_m^{t_0})}
\leq \\
\leq 4c_1\|\b_m \|_{L^1(D_{\d/4},\rho)}+c_2(1+\bm{p_m}\sqrt{t_0})
\leq c_3 (1+\bm{p_m}\sqrt{t_0})
\end{array}
\]
where $c_3$ depends on $\dot{R}'$ only. We conclude as in (a).
\end{itemize}
\end{subsubsection}
\end{subsection}
%

%
\end{section}
%
%
\bibliographystyle{amsalpha}
\bibliography{bib-tri}
\end{document}